\documentclass[reqno]{amsart}
\usepackage{amsmath} 
\usepackage{amssymb}
\usepackage{color}

\usepackage{graphicx}
\usepackage{tikz}
\usepackage{amsfonts,amsmath, amssymb}

\numberwithin{equation}{section}

\usepackage{graphics}
\usepackage{epsfig}
\newtheorem{claim}{\bf \t}[part]

\newtheorem{theorem}{Theorem}[section]

\newtheorem{lemma}[theorem]{Lemma}

\newtheorem{remark}[theorem]{Remark}

\def\v{\varepsilon}

\def\x{\xi}
\def\t{\theta}
\def\k{\kappa}
\def\n{\nu}
\def\m{\mu}
\def\a{\alpha}
\def\b{\beta}
\def\g{\gamma}
\def\d{\delta}
\def\l{\lambda}

\def\r{\rho}
\def\s{\sigma}
\def\z{\zeta}

\def\f{\frac}
\begin{document}

\title[Global Well-posedness of the  Relativistic  Boltzmann Equation]{Global Well-posedness of the  Relativistic  Boltzmann Equation}

\author[Y. Wang]{Yong Wang}
\address[Y. Wang]{Institute of Applied Mathematics, Academy of Mathematics and Systems Science, Chinese Academy of Sciences, Beijing 100190, China, and University of Chinese Academy of Sciences}
\email{yongwang@amss.ac.cn}

\begin{abstract}
In this paper, we prove the global existence and uniqueness of mild solution to the relativistic Boltzmann equation both in the whole space and in torus for a class of initial data with bounded velocity-weighted $L^\infty$-norm and some smallness on $L^1_xL^\infty_p$-norm as well as on  defect mass, energy and entropy. Moreover, the asymptotic stability of the solutions is also investigated in the case of torus. 
\end{abstract}

\date{\today}

\maketitle


\section{Introduction}

The relativistic  Boltzmann equation is written as
\begin{equation}\label{1.1}
p^{\mu}\partial_\mu F=\mathcal{C}(F,F),
\end{equation}
where  the collision operator $\mathcal{C}(F,F)$ takes the bilinear form
\begin{align}\label{1.2}
	 \mathcal{C}(F_1,F_2)=\int_{\mathbb{R}^3}\int_{\mathbb{R}^3}\int_{\mathbb{R}^3} W(p,q|p',q') [F_1(p')F_2(q')-F_1(p)F_2(q)]\f{dp'}{p_0'}\f{dq'}{q_0'}\f{dq}{q_0},
\end{align}
Here the translation rate $W(p,q|p',q')$ is given by
\begin{equation}\label{1.3}
W(p,q|p',q')=\f{c}{2}s\sigma(g,\t) \d^{(4)}(p^\m+q^\m-{p^\m}'-{q^\m}'),
\end{equation}
where $\sigma(g,\t)$ is the scattering kernel measuring the interactions between particles and Dirac function $\d^{(4)}$ is the delta function of four variables. The constant $c>0$ is the light speed.  The relativistic momentum of a particle is denoted by $p^\m$, $\mu=0,1,2,3$. We raise and lower the indices with the Minkowski metric $p^\mu=g_{\mu\nu}p_\n$ where $g_{\mu\nu}=\mbox{diag}=(-1,1,1,1)$.  The signature of the metric is $(-,+,+,+)$.  For $p\in\mathbb{R}^3$, we write $p^\m=(p_0,p)$ where $p_0\doteq\sqrt{|p|^2+c^2}$ is the energy of a relativistic particle with velocity $p$.   We use the Einstein convection of implicit summation over repeated indices, then the Lorentz inner product is given by
\begin{equation*}\label{1.4}
p^\m q_\mu=-p_0q_0+\sum_{i=1}^3p_iq_i.
\end{equation*}
It is noted that the momentum of each particle is restricted to the mass shell $p^\m p_\m=-c^2$ with $p_0>0$.  We refer the interesting readers to \cite{Cercignani1,Groot,Dijkstra,Glassey1,Stewart} for background of the relativistic kinetic theory.

The streaming term of relativistic Boltzmann equation is given by
$$p^\m \partial_\m= p_0\partial_t+cp\cdot\nabla_x.$$
Then we write the relativistic Boltzmann equation \eqref{1.1}  as
\begin{align}\label{1.5}
\partial_t F+\hat{p}\cdot\nabla_xF=Q(F,F),
\end{align}
where $F(t,x,p)$ is a distribution function for fast moving particles at time $t>0$, position $x\in\Omega=\mathbb{R}^3$ or $\mathbb{T}^3$ and particle velocity $p\in\mathbb{R}^3$.
The collision operator  $Q(F,F)\doteq \f{1}{p_0}\mathcal{C}(F,F)$ and the normalized particle velocity $\hat{p}$ is given by
\begin{equation*}\label{1.6}
\hat{p}:= c\f{p}{p_0}\equiv \f{p}{\sqrt{1+|p|^2/c^2}}.
\end{equation*}
We impose the relativistic Boltzmann equation \eqref{1.5} with the following initial data
\begin{equation}\label{1.5-1}
F_0(t,x,p)|_{t=0}=F_0(x,p).
\end{equation}

Now we define the quantity $s$, which is the square of the energy in the "center of momentum"  system, $p+q=0$, as
\begin{equation*}\label{1.7}
s=s(p^\m,q^\m):= -(p^\m+q^\m)(p_\m+q_\m)=2(-p^\m q_\m+c^2)\geq0.
\end{equation*}
The relative momentum $g\geq 0$ is defined as
\begin{equation*}\label{1.8}
g^2=g^2(p^\m,q^\m):=(p^\m-q^\m)(p_\m-q_\m)=2(-p^\m q_\m-c^2)\geq0.
\end{equation*}
A direct calculation shows that $s=g^2+4c^2$. Conversation of momentum and energy for elastic collisions is described as
\begin{equation}\label{1.9}
\begin{cases}
p+q=p'+q',\\
p_0+q_0={p_0}'+{q_0}'.
\end{cases}
\end{equation}
The scattering angle $\t$ is defined by
\begin{equation*}\label{1.10}
\cos\t=\f{(p^\m-q^\m)({p_\m}'-{q_\m}')}{g^2}.
\end{equation*}
This angle is well defined under \eqref{1.9}, see \cite{Glassey}.

The steady solutions of this model are the well known J\"{u}ttner solution, also known as the relativistic Maxwellian, i.e.,
\begin{equation*}
J(p)=\f{e^{-\f{cp_0}{k_B T}}}{4\pi c k_B T K_2(c^2/(k_BT))},
\end{equation*}
where $K_2(\cdot)$ is the Bessel function $K_2(z)=\f{z^2}{2}\int_1^{\infty}e^{-zt}(t^2-1)^{\f32}dt$, $T$ is the temperature and $k_B$ is the Boltzmann constant. Throughout this paper, we normalize all the physical constants to be one, including the speed of light. Then the normalized relativistic Maxwellian becomes
\begin{equation*}\label{1.11}
J(p)=\f{1}{4\pi}e^{-p_0}.
\end{equation*}

\

Using the Lorentz transformations as described in \cite{Groot,Strain3}, one can carry out  the {\it center-of-momentum} expression to  reduce the delta functions and obtain
\begin{align}\label{1.12}
Q(F_1,F_2)&=\int_{\mathbb{R}^3}\int_{\mathbb{S}^2} v_{\phi} \s(g,\t)[F_1(p')F_2(q')-F_1(p)F_2(q)]d\omega dq\nonumber\\
&:= Q_+(F_1,F_2)-Q_-(F_1,F_2),
\end{align}
where $v_{\phi}=v_{\phi}(p,q)$ is the M{$\phi$}ller velocity
\begin{equation}\label{1.13}
v_{\phi}=v_{\phi}(p,q)\doteq\sqrt{\Big|\f{p}{p_0}-\f{q}{q_0}\Big|^2-\Big|\f{p}{p_0}\times\f{q}{q_0}\Big|^2}=\f{g\sqrt{s}}{2p_0q_0}.
\end{equation}
The post-collisional momentum in the expression \eqref{1.12} satisfies
\begin{align}
\begin{cases}
p'=\f{1}{2}(p+q)+\f{1}{2}g\Big(\omega+(\tilde{\g}-1)(p+q)\f{(p+q)\cdot\omega}{|p+q|^2}\Big),\\[2mm]
q'=\f{1}{2}(p+q)-\f{1}{2}g\Big(\omega+(\tilde{\g}-1)(p+q)\f{(p+q)\cdot\omega}{|p+q|^2}\Big).
\end{cases}\nonumber
\end{align}
where $\tilde{\g}=(p_0+q_0)/\sqrt{s}$.  And the energies are given by
\begin{align}
\begin{cases}
p_0'=\f{1}{2}(p_0+q_0)+\f{g}{2\sqrt{s}}(p+q)\cdot\omega,\\[2mm]
q_0'=\f{1}{2}(p_0+q_0)-\f{g}{2\sqrt{s}}(p+q)\cdot\omega.
\end{cases}\nonumber
\end{align}
For other representation of the collision operator, we refer to \cite{Andreasson,Glassey1,Glassey3}.

For  functions $h(p),g(p)$ with sufficient decay  at infinity, the collision operator satisfies
\begin{equation*}\label{1.16}
\int_{\mathbb{R}^3}Q(h,g)dp=\int_{\mathbb{R}^3}pQ(h,g)dp=\int_{\mathbb{R}^3}p_0Q(h,g)dp\equiv0.
\end{equation*}
Let $F$ be a solution of the relativistic Boltzmann equation \eqref{1.5}, formally,  $F$ satisfies the conservations of mass, momentum and energy
\begin{align}
\int_{\Omega}\int_{\mathbb{R}^3}[F(t,x,p)-J(p)]dpdx&=\int_{\Omega}\int_{\mathbb{R}^3}[F_0-J(p)]dpdx=M_0,\label{1.17}\\
\int_{\Omega}\int_{\mathbb{R}^3}p[F(t,x,p)-J(p)]dpdx&=\int_{\Omega}\int_{\mathbb{R}^3}p[F_0-J(p)]dpdx=\tilde{M}_0,\label{1.18}\\
\int_{\Omega}\int_{\mathbb{R}^3}p_0[F(t,x,p)-J(p)]dpdx&=\int_{\Omega}\int_{\mathbb{R}^3}p_0[F_0-J(p)]dpdx=E_0,\label{1.19}
\end{align}
as well as the additional entropy inequality
\begin{equation}
 \int_{\Omega}\int_{\mathbb{R}^3}[F(t)\ln{F(t)}-J\ln{J}]dpdx\leq \int_{\Omega}\int_{\mathbb{R}^3}[F_0\ln{F_0}-J\ln{J}]dpdx.\label{1.20}
\end{equation}
For any function satisfying \eqref{1.17}, \eqref{1.19} and \eqref{1.20}, a standard Taylor expansion shows that 
\begin{equation}\label{1.20-1}
\mathcal{E}(F(t)):= \int_{\Omega}\int_{\mathbb{R}^3}\Big\{F(t)\ln F(t)-J\ln J\Big\}dpdx
+[\ln(4\pi)-1]M_0+E_0\geq 0,
\end{equation}
see \eqref{2.53} for more details.

In 1940 Lichnerowicz-Marrot \cite{L-Marrot} derived the relativistic Boltzmann equation which  is a fundamental model for relativistic particles whose speed is comparable to the speed of light.   The local existence and uniqueness were firstly investigated by Bichteler \cite{Bichteler} in the $L^\infty$ framework under smallness conditions on the initial data. Dudy\'{n}ski and Ekiel-Je\.{z}ewska \cite{Dudy,Dudy-E4} studied  the linearized relativistic Boltzmann equation.  It is well known that the global existence of renormalized solution to the Newtonian Boltzmann equation was proved by DiPerna and Lions \cite{D-Lion} for large initial data, the uniqueness of such solution, however, is unknown.  In 1992, Dudy\'{n}ski and Ekiel-Je\.{z}ewska \cite{Dudy-E5}  obtained  the global existence of the DiPerna-Lions renormalized solution of the relativistic Boltzmann equation by using their results \cite{Dudy-E2,Dudy-E3}. For other interesting works,  see  \cite{Andreasson,Jiang,Jiang1,Wennberg} and the references therein.

On the other hand, when the amplitude of  initial data is small, there are lots of results on the existence and uniqueness of global solutions to the relativistic Boltzmann equation. 1n 1993 Glassey and Strauss \cite{Glassey1} proved the global existence of  smooth solution  on the torus for the relativistic Boltzmann equation, the exponential decay rate was also obtained. It is noted that they \cite{Glassey1} considered only  the hard potential cases. 1995, they \cite{Glassey2} extended that results to  the Cauchy problem. In 2006, Hsiao and Yu \cite{Hsiao-Yu} relaxed the restriction on the cross-section of \cite{Glassey1}, but
is still restricted to the hard potential.  In 2010, Strain \cite{Strain} proved the unique solution of the relativistic Boltzmann equation exists for all time and decay with any polynomial rate towards the relativistic Maxwellian on torus for the soft potentials.  Recently, Jang \cite{Jang} investigated the global classical solutions to the relativistic Boltzmann equation without angular cut-off, which extended the result of  Newtonian Boltzmann equation \cite{Strain2}. For other interesting works, we refer to \cite{Glassey4,Ha} for the case near vacuum,  \cite{Yang-Yu,Strain-Guo1} for Landau system,  \cite{Strain5,Xiaoqh,Liu-Zhao,Yang-Yu1,Yu} for Landau-Maxwell system, \cite{Rein} for Vlasov-Maxwell system and \cite{Guo3} for relativistic Vlasov-Maxwell-Boltzmann equation and the references therein. Along this direction, the very interesting paper \cite{Strain}  is in the frontier of this topic. We would like to mention that based on some new observations, the results of this paper significantly improve the paper \cite{Strain} .

We would like to mention some results on the Newtonian Boltzmann equation. Under a uniform bound assumption in a strong Sobolev space, Desvillettes-Villani \cite{Desvillettes-V} obtained an almost exponential decay rate of large amplitude solutions to the global Maxwellian. The result has been recently improved by Gualdani, Mischler and Mouhot \cite{Gualdani} to a sharp exponential time decay rate. On the other hand, there are many studies on the global existence of small perturbation solutions to the Boltzmann equation, for instance, \cite{Guo1,Liu-Yang-Yu} by using the energy method, \cite{Guo2,Guo,Ukai-Yang} by using $L^2\cap L^\infty$ approach, and \cite{Xu-Yang,Strain2} for non-cutoff Boltzmann equation. For other interesting results, see \cite{B-Guo,Guo4,Guo5,Strain-Guo,Strain-Guo2,Guo6,Villani1} and the references therein. Finally, we  mention some results on  the Newtonian limit of the relativistic Boltzmann equation, see \cite{Calogero,Strain4,Speck} and the references therein.

It is noted that the initial data in \cite{Strain} are required to have small amplitude perturbation in $L^{\infty}_{x,v}$-norm around the global Maxwellian. Recently, the authors \cite{DHWY} developed a new $L^\infty_xL^1_v\cap L^\infty_{x,v}$ approach, and proved the global existence and uniqueness of mild solutions to the Boltzmann equation in the whole space and torus for a class of initial data with bounded velocity-weighted $L^\infty$-norm under some smallness conditions on $L^1_xL^\infty_v$-norm as well as defect mass, energy and entropy. The purpose of this paper is to extend \cite{DHWY} to the relativistic Boltzmann equation, i.e. we consider the global existence and uniqueness of mild solution to the relativistic Boltzmann equation with bounded $L^\infty$-norm and some smallness conditions on $L^1_xL^\infty_p$-norm as well as on defect mass, energy and entropy. The main difficulty is that the collision kernel of the relativistic Boltzmann equation is much more complicated than the non-relativistic case. 

Now we begin to  formulate our main results. Define a weight function
\begin{equation}\label{WF}
w_\b(p):= (1+|p|^2)^{\f\beta2},
\end{equation}
and the perturbation
\begin{equation}\label{1.21}
f(t,x,p):=\f{F(t,x,p)-J(p)}{\sqrt{J(p)}},
\end{equation}
then the relativistic Boltzmann equation \eqref{1.5} is rewritten as
\begin{align}\label{1.22}
f_t+\hat{p}\cdot\nabla_xf+\nu(p)f-Kf=\Gamma(f,f),
\end{align}
where the linearized operator of the Boltzmann equation is
\begin{align}\label{1.23}
Lf=\nu(p) f-Kf=-\f1{\sqrt{J}}\Big\{Q(J,\sqrt{J}f)+Q(\sqrt{J}f, J)\Big\},
\end{align}
the collisional frequency $\nu(p)$ is defined by
\begin{equation}\label{n1.20}
\nu(p)=\int_{\mathbb{R}^3}\int_{\mathbb{S}^2}v_{\phi} \s(g,\t)J(q)d\omega dq,
\end{equation}
and the operator  $K:= K_2-K_1$ are defined as in \cite{Strain}:
\begin{align}
(K_1f)(p)&:=\int_{\mathbb{R}^3}\int_{\mathbb{S}^2}v_{\phi} \s(g,\t)\sqrt{J(p)J(q)}f(q)d\omega dq,\label{1.24}\\
(K_2f)(p)&:=\f1{\sqrt{J}}\Big\{Q_+(J,\sqrt{J}f)+Q_+(\sqrt{J}f,J)\Big\}\nonumber\\
&=\int_{\mathbb{R}^3}\int_{\mathbb{S}^2}v_{\phi} \s(g,\t)\sqrt{J(q)J(q')}f(p')d\omega dq+\int_{\mathbb{R}^3}\int_{\mathbb{S}^2}v_{\phi} \s(g,\t)\sqrt{J(q)J(p')}f(q')d\omega dq,\label{1.25}
\end{align}
and
\begin{align}\label{1.26}
\Gamma(f,f)\equiv\f1{\sqrt{J}}Q(\sqrt{J}f,\sqrt{J}f)
&=\f1{\sqrt{J}}Q_+(\sqrt{J}f,\sqrt{J}f)-\f1{\sqrt{J}}Q_-(\sqrt{J}f,\sqrt{J}f)\nonumber\\
&:= \Gamma_+(f,f)-\Gamma_-(f,f).
\end{align}
Then, for any $(t,x,p)$, the mild form of the relativistic Boltzmann equation \eqref{1.22} is given by
\begin{align}\label{1.27}
f(t,x,p)&=e^{-\nu(p)t}f_0(x-\hat{p}t,p)+\int_0^te^{-\nu(p)(t-s)} (Kf)(s,x-\hat{p}(t-s),p)ds\nonumber\\
&~~~~+\int_0^te^{-\nu(p)(t-s)} \Gamma(f,f)(s,x-\hat{p}(t-s),p)ds,
\end{align}
with initial condition
\begin{equation}\label{1.27-1}
f_0(x,p)=\f{F_0(x,p)-J(p)}{\sqrt{J(p)}}.
\end{equation}

\

To consider the global well-posedness of the relativistic Boltzmann equation, we need the following hypothesis on $\s$:\\[1mm]
{\it
\noindent{\bf H).}  For soft potentials, we assume that the collision kernel of  \eqref{1.5} satisfies
\begin{align}\label{1.30}
\f{g}{\sqrt{s}} g^{-b}\s_0(\t)\lesssim \s(g,\t)\lesssim g^{-b}\s_0(\t),
\end{align}
where $b,\g$ satisfy $b\in(0,2), \g>-\min\{\f43,4-2b\}$.
In addition, we assume that   $\s_0(\t)\lesssim \sin^{\g}\t$ and $\s_0(\t)$ is
  non-zero on a set of positive measure. \\[2mm]
For hard potentials, we assume
\begin{equation}\label{1.31}
\f{g}{\sqrt{s}}g^a\s_0(\t)\lesssim \s(g,\t)\lesssim (g^a+g^{-b})\s_0(\t).
\end{equation}
where  $\g>-\f43, a\in[0,2]\cap[0,\min\{2+\g,4+3\g\}),b\in[0,2)$.}\\[2mm]
We point out that the short range interactions collision kernel is included in the hard potentials above, and the Newtonian limit of the relativistic Boltzmann equation in this case is the hard-sphere Boltzmann equation.\\[2mm]

The first result of this paper is:
\begin{theorem}[Global Existence]\label{thm1.1}
Let  $\Omega=\mathbb{T}^3~\mbox{or}~\mathbb{R}^3$, and {\bf H)} hold. For any given $\b>14,~\bar{M}\geq 1$,  suppose that the initial data $F_0$  satisfies $F_0(x,p)=J(p)+\sqrt{J(p)}f_0(x,p)\geq 0$ and $\|w_\b f_0\|_{L^\infty}\leq \bar{M}$. There is a  small constant $\epsilon_0>0 $  depending on $a,b, \g,\b,\bar{M}$ such that if
	\begin{equation} \label{1.35}
	\mathcal{E}(F_0)+\|f_0\|_{L^1_xL^\infty_p}\leq \epsilon_0,
	\end{equation}
	the Boltzmann equation \eqref{1.5},  \eqref{1.5-1} has a global unique  mild solution $F(t,x,p)=J(p)+\sqrt{J(p)}f(t,x,p)\geq0$  satisfying \eqref{1.17}-\eqref{1.20} and
	\begin{align} \label{1.36}
	\|w_\b f(t)\|_{L^\infty}\leq \tilde{C}_1\bar{M}^2,
	\end{align}
	where the positive constant $\tilde{C}_1$  depends only on $a,b,\g,\b$. Moreover, if the initial data $f_0$ is continuous in $(x,p)\in\Omega\times\mathbb{R}^3$, then the solution $f(t,x,p)$ is continuous in $[0,\infty)\times\Omega\times\mathbb{R}^3$.
\end{theorem}
\begin{remark}
It is noted that there exists a lot of initial data satisfying \eqref{1.35}. For example, we take 
\begin{equation}\nonumber
F_0(x,p)=\rho_0(x) J(p), \quad (x,p)\in\Omega\times\mathbb{R}^3,
\end{equation}
with $\r_0(x)\geq0$,  $\r_0\in L^\infty_x$, $\r_0-1\in L^1$ and $\r_0\ln\r_0-\r_0+1\in L^1_x$.   Then, it is direct to check that  
\begin{equation}
\mathcal{E}(F_0)+\|f_0\|_{L^1_xL^\infty_p}\leq \|\r_0\ln\r_0-\r_0+1\|_{L^1}+C\|\r_0-1\|_{L^1}.
\end{equation} 
Therefore if $\|\r_0\ln\r_0-\r_0+1\|_{L^1}+C\|\r_0-1\|_{L^1}$ is small, then \eqref{1.35} holds, and the initial data is allowed to have large oscillations in $L^\infty_{x,p}$, see the following figure:
\begin{figure}[h!]
	\centering
	\includegraphics[width=8cm, height=4cm]{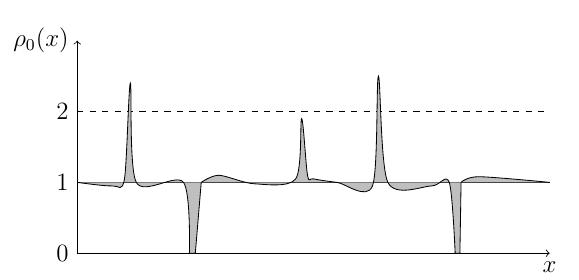}
		\caption{}\label{figure1.1}
\end{figure} 
\end{remark}

\begin{remark}
As pointed out in \cite{Strain,Dudy-E4}, the full ranges should be $\g>-2,a\in[0,2+\g],b\in[0,\min\{4,4+\g\})$ for hard potentials, and $\g>-2, b\in(0,\min\{4,4+\g\})$ for soft potentials. In this paper, due to some technique difficulties, we need the restrictions on $\g, a, b$ as in {\bf H)}. It is an interesting open problem to consider the remaining cases.  Indeed it is not known how to construct the local solution to the relativistic Boltzmann equation with $L^\infty$ bounded  initial data for hard potentials with $2<a\leq 2+\g,\g>0$, see Theorem \ref{thm7.1} below.
\end{remark}

Furthermore,  one can  obtain the following decay estimates for the solutions obtained in Theorem \ref{thm1.1} in the case of torus $\Omega=\mathbb{T}^3$.
\begin{theorem}[Decay Estimate for Hard Potentials]\label{thm1.2}
For hard potentials, let $\Omega=\mathbb{T}^3$, $\b>14$ and  $\g>-\f43, a\in[0,2]\cap[0,\min\{2+\g,4+3\g\}),b\in[0,2)$. Assume  $(M_0,\tilde{M}_0,E_0)=(0,0,0)$,  and  $\epsilon_0>0$ sufficiently small, then there exists a positive constant $\l_0>0$ such that the  solution $f(t,x,p)$ obtained in Theorem \ref{thm1.1} satisfies
\begin{align}\label{1.29}
\|w_\b f(t)\|_{L^\infty}\leq \tilde{C_2}e^{-\l_0 t},
\end{align}
where $\tilde{C}_2>0$ is a positive constant  depending only $a,b,\g,\b$ and $\bar{M}$.
\end{theorem}

\begin{theorem}[Decay Estimate for Soft Potentials]\label{thm1.3}
For soft potentials, let $\Omega=\mathbb{T}^3$, $\b>14$ and  $b\in(0,2), \g>-\min\{\f43,4-2b\}$. Assume  $(M_0,\tilde{M}_0,E_0)=(0,0,0)$,  and  $\epsilon_0>0$ sufficiently small, then  the  solution $f(t,x,p)$ obtained in Theorem \ref{thm1.1} satisfies
\begin{align}\label{1.38}
\|f(t)\|_{L^\infty}\leq \tilde{C}_3(1+t)^{-1-\f{\xi_1}{b}}.
\end{align}
where the positive constant $\xi_1>0$ is defined in Lemma \ref{lem2.4} below,   and $\tilde{C}_3>0$ depends only on $a,b,\g,\b, \bar{M}$.
\end{theorem}

\begin{remark}
From \eqref{1.38} and \eqref{1.36}, we have $\|w_{\f\b2}f\|_{L^\infty}\leq C(1+t)^{-\f12}$, which yields  $\|w_{\f\b2}f\|_{L^\infty}\ll1$ when $t\gg1$.Then, one can apply the iteration method in Section 6, 7 of \cite{Strain} to improve the decay rate to any polynomial when $\b$ is  large enough. However, we shall not discuss it in this paper since the main aim of this paper is the existence of global solution with uniqueness for the relativistic Boltzmann equation.
\end{remark}

Now we explain the strategy of the proof of Theorem \ref{thm1.1}. As mentioned previously, the  only global existence of large solutions to the relativistic Boltzmann equation is due to Dudy\'{n}ski and Ekiel-Je\.{z}ewska \cite{Dudy-E5}, the uniqueness  of these renormalized solutions, however, is completely open due to the lack of $L^\infty$ estimates.  Indeed, it is difficult to establish the global $L^\infty$ bound for the solutions of relativistic Boltzmann equation due to the nonlinear term $\Gamma(f,f)$.  In the previous references \cite{Strain,Glassey1}, one usually  bounds the nonlinear term in the following way
\begin{align}\label{1.39}
	|w_\b(p)\Gamma(f,f)(t)|\leq C\nu(p)\|w_\b f(t)\|^2_{L^\infty},
\end{align}
then the smallness  assumption on the $L^\infty$-norm is needed. Indeed,  it is hard to prove  even the local existence of solution to the relativistic Boltzmann equation with general bounded $L^\infty$-norm initial data by using \eqref{1.39}  for hard potentials. In this paper, we firstly establish a new bound  on the gain term $\Gamma_+(f,f)$(see \eqref{7.62-1} below), i.e.,
\begin{align}\label{1.32}
|w_\b(p) \Gamma_+(f,f)(t)|\leq C\|w_\b f(t)\|^2_{L^\infty},~\mbox{for suitably large } \b>0,
\end{align}
which enable us to obtain the local existence of $L^\infty$ solution to the relativistic  Boltzmann equation without any smallness assumption on the $L^\infty$-norm of initial data, see Theorem \ref{thm7.1} below.

Although we have obtained the local solution with general bounded $L^\infty_{x,p}$ initial data, but  it is very difficult to extend such local solution to a  global one due to the difficulty of quadratic  term $\Gamma(f,f)$.  To avoid the smallness assumption on the $L^\infty$-norm, motivated by \cite{DHWY}, we firstly establish the following estimate for the  nonlinear term $\Gamma(f,f)$ of relativistic Boltzmann equation(see Lemma \ref{lem4.1} below), i.e., for $\b\geq1$,
\begin{align}\label{1.40}
&\Big|w_\b(p)\Gamma(f,f)(t,x,p)\Big|
\leq  C\nu(p)\|w_{\b}f(t)\|^{2-\vartheta}_{L^\infty}\cdot\Big(\int_{\mathbb{R}^3}|f(t,x,q)|dq\Big)^{\vartheta},
\end{align}
for some $0<\vartheta<1$. We remark that one should be very careful to establish the above two inequalities \eqref{1.32} and \eqref{1.40} due to the complexity of  cross-sections and the Lorentz transformation for the relativistic Boltzmann equation. Indeed, we need Lemmas \ref{lemA.2} and \ref{lemA.3}(see appendix), which refine the corresponding lemmas in \cite{Glassey1}. 

Finally, based on the above preparation and  under the initial condition \eqref{1.35}, we prove that $\int_{\mathbb{R}^3}|f(t,x,q)|dq$ should be small after some positive time due to the hyperbolicity of  relativistic Boltzmann equation,  even though $\int_{\mathbb{R}^3}|f_0(x,q)|dq$ may be  large initially.   Then  we can finally establish following uniform estimate
\begin{equation*}\label{1.41}
\sup_{0\leq s\leq t}\|w_\b f(s)\|_{L^\infty}\leq C\bar{M}^2.
\end{equation*}
through careful analysis. It is noted that  the smallness of $\mathcal{E}(F_0)+\|f_0\|_{L^1_xL^\infty_p}$  implies that the initial data may have large oscillations.

\

\noindent{\bf Organization of the paper.}   In section 2, we give some useful estimates which will be used frequently. Section 3 is devoted to the local existence of unique solution to the relativistic Boltzmann equation with arbitrary $L^\infty$ data. In section 4, we first establish a key inequality Lemma \ref{lem4.1}, then give the details of proof of Theorem \ref{thm1.1}. Section 5 is devoted to the decay estimates in the case of torus. 

\

\noindent{\bf Notations.}  Throughout this paper, we will use the $L^2$ norms
\begin{align}
\|h\|_{L^2}:=\left(\int_{\Omega}\int_{\mathbb{R}^3}|h(x,p)|^2dpdx\right)^{\f12},~~~
|h|_{L^2}:=\left(\int_{\mathbb{R}^3}|h(p)|^2dp\right)^{\f12},\nonumber
\end{align}
and $\|\cdot\|_{L^\infty}$ denotes the $L^\infty(\Omega\times\mathbb{R}^3_p)$-norm.
The $L^2(\mathbb{R}^3_p)$ inner product is denoted $\langle\cdot,\cdot\rangle$. We also need to measure the dissipation of the linearized operator
\begin{align}\nonumber
\|h\|_{\nu}:=\left(\int_{\Omega}\int_{\mathbb{R}^3}\nu(p)|h(x,p)|^2dpdx\right)^{\f12},~~~
|h|_{\nu}:=\left(\int_{\mathbb{R}^3}\nu(p)|h(p)|^2dp\right)^{\f12}.
\end{align}
We will further use $A\lesssim B$ to mean that there exists a positive constant $C>0$ such that $A\leq CB$ holds uniformly over the range of parameters which are present in the inequality and  the precise magnitude of the constant is not
important. The notation $B\gtrsim A$ is equivalent to $A\lesssim B$, and $A\approx B$ means that $A\lesssim B$ and $B\lesssim A$. We also use
$C>0$ to denote a generic positive constant  which may depend on $\g,\b$ and  vary from line to line, and $c>0$ to denote a small constant.  $C_\vartheta,\cdots$ denote the generic positive constants depending on $\vartheta,\cdots$, respectively, which also may vary from line to line.


\section{Preliminaries}

Define
\begin{align}\label{A.4}
l:=\f{p_0+q_0}{2},~~j:= \f{|p\times q|}{g},
\end{align}
From \cite{Groot,Dudy,Strain1}, we know that
\begin{align}\label{2.1}
(K_if)(p)=\int_{\mathbb{R}^3}k_i(p,q)f(q)dq,~~i=1,2,
\end{align}
with the symmetric kernels
\begin{align}
k_1(p,q)&=c_1\f{g\sqrt{s}}{p_0q_0}e^{-l}\int_0^{\pi}\s(g,\t)\sin\t d\t,\label{2.3}\\
k_2(p,q)&=c_2\f{s^{\f32}}{gp_0q_0}\int_0^\infty \f{y[1+\sqrt{1+y^2}]}{\sqrt{1+y^2}}\s\Big(\f{g}{\sin\f{\psi}{2}},\psi\Big) e^{-l\sqrt{1+y^2}}I_0(jy)dy,\label{2.4}
\end{align}
where $c_1>0$, $c_2>0$ are positive constants and
\begin{equation*}\label{2.6}
0\leq \sin\f{\psi}{2}=\f{\sqrt2 g}{[g^2-4+(g^2+4)\sqrt{1+y^2}]^{\f12}}.
\end{equation*}
The modified Bessel function $I_0(x)$ of imaginary function is defined as
\begin{equation*}\label{A.1}
I_0(z):=\f{1}{2\pi}\int_{0}^{2\pi} e^{z\cos\varphi}d\varphi.
\end{equation*}

\begin{lemma}[Glassey\& Strauss \cite{Glassey1}]\label{lem2.1}
It holds that
\begin{align}
\f{[|p\times q|^2+|p-q|^2]^{\f12}}{\sqrt{p_0q_0}}&\leq g\leq |p-q|~~\mbox{and}~~g\leq 2\sqrt{p_0q_0},\label{2.7}\\
v_{\phi}&=\f{g\sqrt{s}}{p_0q_0}\lesssim  1,\label{2.8}\\
l^2-j^2=\f{g^2+4}{4g^2}&|p-q|^2\geq 1+\f14|p-q|^2,\label{2.9}\\
\f{1}{\sqrt2}g(1+y^2)^{\f14}\leq&\f{g}{\sin\f{\psi}{2}}\leq \sqrt{s}(1+y^2)^{\f14},\label{2.10}\\
\f{y}{2\sqrt{1+y^2}}\leq &\cos\f{\psi}{2}\leq 1.\label{2.11}
\end{align}
\end{lemma}

\

We define
\begin{equation}\label{2.15-1}
\s_a(g,\psi):=g^{a} \sin^{\g}\psi~\mbox{and}~\s_b(g,\psi):= g^{-b} \sin^{\g}\psi,
\end{equation}
and
\begin{align}
k_{2a}(p,q):=\f{s^{\f32}}{gp_0q_0}\int_0^\infty y e^{-l\sqrt{1+y^2}}I_0(jy)\cdot\s_a\Big(\f{g}{\sin\f\psi2},\psi\Big) dy,\label{2.15}\\
k_{2b}(p,q):=\f{s^{\f32}}{gp_0q_0}\int_0^\infty  y e^{-l\sqrt{1+y^2}}I_0(jy)\cdot\s_b\Big(\f{g}{\sin\f\psi2},\psi\Big)dy.\label{2.16}
\end{align}
Then  the following estimates hold:
\begin{lemma}\label{lem2.3}
Under the assumptions of \eqref{1.30} and \eqref{1.31}, it holds that
\begin{align}
0\leq k_1(p,q)\lesssim \Big\{1+|p-q|^{1-b}\Big\}e^{-\f14(p_0+q_0)}.\label{2.11-3}
\end{align}
and
\begin{align}\label{2.12}
k_2(p,q)\lesssim
\begin{cases}
k_{2a}(p,q)+k_{2b}(p,q),~~\mbox{for hard potentials},\\
k_{2b}(p,q),~~\mbox{for soft potentials},
\end{cases}
\end{align}
where $k_{2a}(p,q)$ and $k_{2b}(p,q)$ satisfy
\begin{align}\label{2.13}
k_{2a}(p,q)\lesssim
\begin{cases}
(p_0q_0)^{\f{a-\g-2}{4}}e^{-\f18|p-q|},~~\mbox{for}~~a\geq\g\geq1,\\[1.5mm]
(p_0q_0)^{\f{a-\g-2}{4}}\f{(p_0q_0)^{\f{1-\g}2}}{[|p\times q|+|p-q|]^{1-\g}}e^{-\f18|p-q|},~~\mbox{for}~~a\geq\g\geq0,~\gamma<1,\\[1.5mm]
(p_0q_0)^{-\f12+\f{\zeta_1}{4}}e^{-\f18|p-q|},~~\mbox{for}~~1\leq a<\g,\\[1.5mm]
(p_0q_0)^{-\f12+\f{\zeta_1}{4}}\f{(p_0q_0)^{\f{1-a}2}}{[|p\times q|+|p-q|]^{1-a}}e^{-\f18|p-q|},~~\mbox{for}~~a<\g,~a<1,\\[1.5mm]
(p_0q_0)^{\f{a+|\g|-2}4}\f{(p_0q_0)^{\f{1+|\g|}2}}{[|p\times q|+|p-q|]^{1+|\g|}}  e^{-\f18|p-q|},~\mbox{for}~-2<\gamma< 0,~a\leq 2+\g,
\end{cases}
\end{align}
and
\begin{align}\label{2.14}
k_{2b}(p,q)\lesssim
\begin{cases}
(p_0q_0)^{-\f12-\f{\zeta_2}{4}}\f{(p_0q_0)^{\f{1+b}2}}{[|p\times q|+|p-q|]^{1+b}}e^{-\f18|p-q|},~\mbox{for}~\gamma\geq 0,~b<2,\\[1.5mm]
(p_0q_0)^{\f14(|\g|-b-2)}\f{(p_0q_0)^{\f{1+b}2}}{[|p\times q|+|p-q|]^{1+b}}e^{-\f18|p-q|},~\mbox{for}~-b<\gamma< 0,~|\g|< b<2,\\[1.5mm]
(p_0q_0)^{\f14(|\g|-b-2)}\f{(p_0q_0)^{\f{1+|\g|}2}}{[|p\times q|+|p-q|]^{1+|\g|}}e^{-\f18|p-q|},~\mbox{for}~-2<\gamma< 0,~|\g|\geq b,
\end{cases}
\end{align}
where $\zeta_1=\max\{-2,a-\g\}$, $\zeta_2=\min\{2,b+\g\}$.
\begin{remark}
Here we assume $0\leq b<2$ to guarantee the integrability of $k_b(p,q)$,  i.e., $\int_{\mathbb{R}^3}k_{2b}(p,q)dq<\infty$.
\end{remark}
\end{lemma}
\noindent{\bf Proof.} It is noted that \eqref{2.11-3} follows easily from \eqref{2.3}. Now we focus on the estimation of $k_2(p,q)$ which is much more complicated.
It is noted that \eqref{2.12} follows easily from \eqref{1.30} and \eqref{1.31}. In the following, we try to estimate \eqref{2.13} and \eqref{2.14}.\\[2mm]
\noindent{\underline{Estimation of $k_{2a}(p,q)$}:}  Noting
\begin{equation}\label{2.17}
\sin\psi=2\sin\f\psi2 \cos\f\psi2,
\end{equation}
which, together with  \eqref{2.15}, yields that
\begin{align}\label{2.18}
k_{2a}(p,q)\lesssim\f{g^{\g-1}s^{\f32}}{p_0q_0}\int_0^\infty y e^{-l\sqrt{1+y^2}}I_0(jy)\Big(\f{g}{\sin\f\psi2}\Big)^{a-\g} \cos^{\g}\f\psi2dy,
\end{align}
We divide the proof into the following cases. \\[1.5mm]
\noindent{\it Case 1:}  For $\g\geq0$. \\[1.5mm]
\noindent{1)} For  $a\geq\g\geq0$, it follows from \eqref{2.9}-\eqref{2.11}, \eqref{2.18} and Lemma \ref{lemA.2} that
\begin{align}
k_{2a}(p,q)&\lesssim\f{g^{\g-1}s^{\f32+\f{a-\g}2}}{p_0q_0}\int_0^\infty y e^{-l\sqrt{1+y^2}}I_0(jy) (1+y^2)^{\f{a-\g}4} dy\nonumber\\
&\lesssim \f{g^{\g-1}s^{\f32+\f{a-\g}2}}{p_0q_0} \f{l^{1+\f{a-\g}2}}{(l^2-j^2)^{1+\f{a-\g}{4}}}e^{-\sqrt{l^2-j^2}}\label{2.19-1}\\
&\lesssim \f{g^{\g-1}s^{\f32+\f{a-\g}2}}{p_0q_0} l^{1+\f{a-\g}2}e^{-\f12|p-q|}
\lesssim\f{g^{\g-1}}{p_0q_0} l^{1+\f{a-\g}2}e^{-\f38|p-q|},\label{2.19}
\end{align}
where we have used the facts $0\leq a-\g\leq 2$ and  $s= 4+g^2\leq 4+|p-q|^2$. If $\g\geq1$, then it follows from \eqref{2.7}, \eqref{2.19} and  \eqref{A.9} that
\begin{align}\label{2.20}
k_{2a}(p,q)&
\lesssim\f{ l^{1+\f{a-\g}2}}{p_0q_0}e^{-\f14|p-q|}\lesssim (p_0q_0)^{\f{a-\g-2}4}e^{-\f18|p-q|},~~\mbox{for}~a\geq\g\geq1.
\end{align}
If $0\leq \g<1$, it follows from \eqref{2.19}, \eqref{2.7} and \eqref{A.9} that
\begin{align}\label{2.21}
k_{2a}(p,q)&\lesssim\f{(p_0q_0)^{\f{1-\g}2}}{[|p\times q|+|p-q|]^{1-\g}}\f{l^{1+\f{a-\g}2}}{p_0q_0} e^{-\f38|p-q|}\nonumber\\
&\lesssim \f{(p_0q_0)^{\f{1-\g}2}}{[|p\times q|+|p-q|]^{1-\g}} (p_0q_0)^{\f14(a-\g-2)} e^{-\f18|p-q|},~~\mbox{for}~a\geq \g,~0\leq\g<1.
\end{align}

\

\noindent{2)} For  $0\leq a<\g$, it follows from \eqref{2.9}-\eqref{2.11}, \eqref{2.18} and Lemma \ref{lemA.2} that
\begin{align}
k_{2a}(p,q)&\lesssim\f{g^{\g-1}s^{\f32}}{p_0q_0}\int_0^\infty y e^{-l\sqrt{1+y^2}}I_0(jy) [g(1+y^2)^{\f14}]^{a-\g} dy\nonumber\\
&\lesssim \f{g^{a-1}s^{\f32}}{p_0q_0} \f{l^{1+\f{\z_1}2}}{(l^2-j^2)^{1+\f{\z_1}{4}}}e^{-\sqrt{l^2-j^2}}\label{2.22-1}\\
&\lesssim \f{g^{a-1}s^{\f32}}{p_0q_0} l^{1+\f{\z_1}2}e^{-\f12|p-q|}
\lesssim \f{g^{a-1}}{p_0q_0} l^{1+\f{\z_1}2}e^{-\f38|p-q|}, \label{2.22}
\end{align}
where $\zeta_1=\max\{-2,a-\g\}\leq0$. If $a\geq1$, it follows from \eqref{2.7}, \eqref{2.22} and \eqref{A.9}  that
\begin{align}\label{2.23}
k_{2a}(p,q)&\lesssim \f{l^{1+\f{\z_1}2}}{p_0q_0} e^{-\f14|p-q|}
\lesssim   (p_0q_0)^{-\f12+\f{\z_1}4} e^{-\f18|p-q|},~\mbox{for}~1\leq a<\g.
\end{align}
On the other hand, if $0\leq a<1$, it follows from \eqref{2.7}, \eqref{2.22} and \eqref{A.9}  that
\begin{align}\label{2.24}
k_{2a}(p,q)&\lesssim \f{(p_0q_0)^{\f{1-a}2}}{[|p\times q|+|p-q|]^{1-a}}\f{l^{1+\f{\z_1}2}}{p_0q_0} e^{-\f38|p-q|}\nonumber\\
&\lesssim  (p_0q_0)^{-\f12+\f{\z_1}4} \f{(p_0q_0)^{\f{1-a}2}}{[|p\times q|+|p-q|]^{1-a}}  e^{-\f18|p-q|},
~\mbox{for}~0\leq a<\g,~a<1.
\end{align}

\noindent{\it Case 2:}  For $-2<\g<0$, it follows from \eqref{2.18} and \eqref{2.9}-\eqref{2.11} that
\begin{align}
k_{2a}(p,q)&\lesssim\f{g^{\g-1}s^{\f32}}{p_0q_0}\int_0^\infty y e^{-l\sqrt{1+y^2}}I_0(jy) \Big(\sqrt{s}(1+y^2)^{\f14}\Big)^{a+|\g|}\cdot\Big(\f{y}{\sqrt{1+y^2}}\Big)^{-|\g|} dy\nonumber\\
&\lesssim \f{g^{\g-1}s^{\f32+\f{a+|\g|}2}}{p_0q_0}\Big\{\int_0^1 y^{1-|\g|} e^{-l\sqrt{1+y^2}}I_0(jy) dy+\int_1^\infty y e^{-l\sqrt{1+y^2}}I_0(jy) (1+y^2)^{\f{a+|\g|}4} dy\Big\}\nonumber\\
&~~~~\mbox{Noting }~a\leq 2+\g, \eqref{A.5} ~\mbox{and}~\eqref{A.6}, \mbox{we have }\nonumber\\
&\lesssim \f{g^{\g-1}s^{\f32+\f{a+|\g|}2}}{p_0q_0} \left(\f{1}{(l^2-j^2)^{\f{2-|\g|}4-O(\v)}}+\f{l^{1+\f{a+|\g|}{2}}}{(l^2-j^2)^{1+\f{a+|\g|}{4}}}\right)e^{-\sqrt{l^2-j^2}}\label{2.25-1}\\
&\lesssim \f{g^{\g-1}s^{\f32+\f{a+|\g|}2}}{p_0q_0} l^{1+\f{a+|\g|}{2}}e^{-\f12|p-q|}
\lesssim\f{(p_0q_0)^{\f{1+|\g|}2}}{[|p\times q|+|p-q|]^{1+|\g|}}  (p_0q_0)^{\f{a+|\g|-2}4}e^{-\f14|p-q|},\label{2.25}
\end{align}
where $O(\v)\doteq \f{2(2-|\g|)^2\v}{4(4+4\v-2|\g|\v)}$ with $\v>0$ small enough so that $\f12-\f{|\g|}4-O(\v)\geq0$. Thus combining \eqref{2.20}-\eqref{2.25}, we have proved \eqref{2.13}.

\

\noindent{\underline{Estimation of $k_{2b}(p,q)$}:}
It follows from  \eqref{2.16} and \eqref{2.17} that
\begin{align}\label{2.26}
k_{2b}(p,q)\lesssim\f{g^{\g-1}s^{\f32}}{p_0q_0}\int_0^\infty y e^{-l\sqrt{1+y^2}}I_0(jy)\Big(\f{g}{\sin\f\psi2}\Big)^{-b-\g} \cos^{\g}\f\psi2dy.
\end{align}
As previous, we divide the proof into the following cases. \\[1.5mm]
\noindent{\it Case 1:}  For $\g\geq0$, noting $\z_2=\min\{2,b+\g\}$,  it follows from  \eqref{2.9}, \eqref{2.10}, \eqref{A.5} and \eqref{A.9}   that
\begin{align}
k_{2b}(p,q)&\lesssim\f{g^{\g-1}s^{\f32}}{p_0q_0}\int_0^\infty y e^{-l\sqrt{1+y^2}}I_0(jy)[g(1+y^2)^{\f14}]^{-b-\g} dy\nonumber\\
&\lesssim \f{g^{-b-1}s^{\f32}}{p_0q_0}\int_0^\infty y e^{-l\sqrt{1+y^2}}I_0(jy)(1+y^2)^{-\f{b+\g}4} dy\nonumber\\
&\lesssim \f{g^{-b-1}s^{\f32}}{p_0q_0}\f{l^{1-\f{\z_2}2}}{(l^2-j^2)^{1-\f{\z_2}{4}}}e^{-\sqrt{l^2-j^2}}
\lesssim (p_0q_0)^{-\f12-\f{\zeta_2}{4}}\f{(p_0q_0)^{\f{1+b}2}}{[|p\times q|+|p-q|]^{1+b}}e^{-\f18|p-q|}.\label{2.27}
\end{align}

\

\noindent{\it Case 2:}  For $-2<\g<0$. \\[1.5mm]
\noindent{1)} For  $|\g|<b$, it follows from \eqref{2.9}-\eqref{2.11}, \eqref{2.26} \eqref{A.5} and \eqref{A.6} that
\begin{align}
k_{2b}(p,q)&\lesssim\f{g^{\g-1}s^{\f32}}{p_0q_0}\int_0^\infty y e^{-l\sqrt{1+y^2}}I_0(jy)
[g(1+y^2)^{\f14}]^{-b+|\g|} \Big(\f{y}{\sqrt{1+y^2}}\Big)^{-|\g|}dy\nonumber\\
&\lesssim \f{g^{-b-1}s^{\f{3}2}}{p_0q_0}\Big\{\int_0^1y^{1-|\g|}e^{-l\sqrt{1+y^2}}I_0(jy)dy
+\int_1^{\infty}y e^{-l\sqrt{1+y^2}}I_0(jy)(1+y^2)^{\f{-b+|\g|}4}dy\Big\}\nonumber\\
&\lesssim \f{g^{-b-1}s^{\f{3}2}}{p_0q_0} \left(\f1{(l^2-j^2)^{\f{2-|\g|}{4}-O(\v)}}+\f{l^{1+\f{-b+|\g|}{2}}}{(l^2-j^2)^{1+\f{-b+|\g|}{4}}}\right)e^{-\sqrt{l^2-j^2}}\nonumber\\
&\lesssim (p_0q_0)^{\f14(|\g|-b-2)}\f{(p_0q_0)^{\f{1+b}2}}{[|p\times q|+|p-q|]^{1+b}}e^{-\f18|p-q|},~\mbox{for}~-2<\gamma< 0,~|\g|<b<2.\label{2.29}
\end{align}

\noindent{2)} For  $b\leq |\g|$, it follows from \eqref{2.9}-\eqref{2.11}, \eqref{2.26} \eqref{A.5} and \eqref{A.6} that
\begin{align}
k_{2b}(p,q)&\lesssim\f{g^{\g-1}s^{\f32}}{p_0q_0}\int_0^\infty y e^{-l\sqrt{1+y^2}}I_0(jy)
[\sqrt{s}(1+y^2)^{\f14}]^{-b+|\g|} \Big(\f{y}{\sqrt{1+y^2}}\Big)^{-|\g|}dy\nonumber\\
&\lesssim \f{g^{\g-1}s^{\f{3+|\g|-b}2}}{p_0q_0}\Big\{\int_0^1y^{1-|\g|}e^{-l\sqrt{1+y^2}}I_0(jy)dy
+\int_1^{\infty}y e^{-l\sqrt{1+y^2}}I_0(jy)(1+y^2)^{\f{-b+|\g|}4}dy\Big\}\nonumber\\
&\lesssim \f{g^{\g-1}s^{\f{3+|\g|-b}2}}{p_0q_0} \left(\f1{(l^2-j^2)^{\f{2-|\g|}{4}-O(\v)}}+\f{l^{1+\f{-b+|\g|}{2}}}{(l^2-j^2)^{1+\f{-b+|\g|}{4}}}\right)e^{-\sqrt{l^2-j^2}}\nonumber\\
&\lesssim (p_0q_0)^{\f14(|\g|-b-2)}\f{(p_0q_0)^{\f{1+|\g|}2}}{[|p\times q|+|p-q|]^{1+|\g|}}e^{-\f18|p-q|},~\mbox{for}~-2<\gamma< 0,~|\g|\geq b.\label{2.28}
\end{align}
 Combining \eqref{2.27}-\eqref{2.28}, we completed the proof of \eqref{2.14}.  $\hfill\Box$

\begin{lemma}\label{lem2.4}
For soft potentials and $b\in(0,2), \g>-\min\{\f43,4-2b\}$,  it holds that
\begin{align}\label{2.38}
\int_{\mathbb{R}^3}k_2(p,q)dq\lesssim p_0^{-\f{b}2-\xi_1}\lesssim 1,
\end{align}
with $\xi_1:=\f14\min\Big\{1, 2-b, 4-2b+\g, 4+3\g\Big\}>0$.\\[2mm]
For hard potentials and  $ \g>-\f43, a\in[0,\min\{2+\g,4+3\g\}), b\in[0,2)$, it holds that
\begin{align}\label{2.39}
\int_{\mathbb{R}^3}k_{2}(p,q)dq\lesssim p_0^{-\xi_2},
\end{align}
with $\xi_2:=\f14\min\Big\{1,2+\g-a,4+3\g-a,2-b\Big\}>0$.
\end{lemma}
\noindent{\bf Proof.}   From \eqref{2.12}, we need only to estimates $\int_{\mathbb{R}^3}k_{2b}(p,q)dq$ and  $\int_{\mathbb{R}^3}k_{2a}(p,q)dq$.\\[1.5mm]
\noindent{\underline{Estimation on $\int_{\mathbb{R}^3}k_{2b}(p,q)dq$}:}\\[1.5mm]
\noindent{1).} For $\g\geq0$, it follows from $\eqref{2.14}_1$, \eqref{A.10} and \eqref{A.17} that
\begin{align}\label{2.40}
\int_{\mathbb{R}^3}k_{2b}(p,q)dq&\lesssim p_0^{-1-\f{\z_2}2} p_0^{1+b}
\int_{\mathbb{R}^3}\f{e^{-\f1{16}|p-q|}}{[|p\times q|+|p-q|]^{1+b}} dq\nonumber\\
&\lesssim
\begin{cases}
p_0^{-1-\f{\z_2}2}\ln p_0,~\mbox{for}~0\leq b\leq 1,\\
p_0^{-2-\f{\z_2}2+b},~\mbox{for}~1< b <2,
\end{cases}
\lesssim p_0^{-\f{b}2} \left(p_0^{\f12(b-2-\z_2)}\ln p_0+p_0^{b-2+\f{b-\z_2}{2}}\right)\nonumber\\
&\lesssim p^{-\f{b}2} \left(p_0^{-1}\ln p_0+p_0^{b-2}\right),
\end{align}
where we have used the fact $\zeta_2\geq b$ since $0\leq b<2$.\\[2mm]
\noindent{2).} For $-b<\g<0,~b>|\g|$, it follows from $\eqref{2.14}_2$, \eqref{A.10} and \eqref{A.17} that
\begin{align}\label{2.42}
\int_{\mathbb{R}^3}k_{2b}(p,q)dq&\lesssim p_0^{\f12(|\g|+b)}
\int_{\mathbb{R}^3}\f{e^{-\f1{16}|p-q|}}{[|p\times q|+|p-q|]^{1+b}} dq\nonumber\\
&\lesssim
\begin{cases}
p_0^{\f12(|\g|-b-2)}\ln p_0,~\mbox{for}~ 0\leq b\leq 1,-b<\g<0,\\
p_0^{\f12(|\g|+b-4)} ,~\mbox{for}~1< b<2,~-b<\g<0.
\end{cases}
\end{align}

\noindent{3).} For $-2<\g<0,~b\leq|\g|$, it follows from $\eqref{2.14}_3$, \eqref{A.10} and \eqref{A.17} that
\begin{align}\label{2.41}
\int_{\mathbb{R}^3}k_{2b}(p,q)dq&\lesssim p_0^{\f12(3|\g|-b)}
\int_{\mathbb{R}^3}\f{e^{-\f1{16}|p-q|}}{[|p\times q|+|p-q|]^{1+|\g|}} dq\nonumber\\
&\lesssim
\begin{cases}
p_0^{\f12(|\g|-b-2)}\ln p_0,~\mbox{for}~ -1\leq\g< 0,\g\leq -b,\\
p_0^{\f12(3|\g|-b-4)},~\mbox{for}~ -2<\g<-1, \g\leq -b.
\end{cases}
\end{align}
Combining \eqref{2.40}-\eqref{2.41},  we have, for $b\in(0,2), \g>-\min\{\f43,4-2b\}$, that
\begin{align}\label{2.45}
\int_{\mathbb{R}^3}k_{2b}(p,q)dq\lesssim p_0^{-\f{b}2} p_0^{-\xi_1},
\end{align}
which yields  immediately \eqref{2.38}.   On the other hand, for $0\leq b<2, \g>-\f43$,  it follows from \eqref{2.40}-\eqref{2.41} that
\begin{align}\label{2.45-1}
\int_{\mathbb{R}^3}k_{2b}(p,q)dq\lesssim p_0^{-\xi_{21}},
\end{align}
where $\xi_{21}\doteq\f14\min\Big\{1,2-b,2+\g,4+3\g\Big\}>0$.

\

\noindent{\underline{Estimation on $\int_{\mathbb{R}^3}k_{2a}(p,q)dq$}:}\\[1.5mm]
\noindent{1).} For $\g\geq0$, it follows from $\eqref{2.13}_1$-$\eqref{2.13}_4$, \eqref{A.10} and \eqref{A.17} that
\begin{align}\label{2.43}
\int_{\mathbb{R}^3}k_{2a}(p,q)dq\lesssim  p_0^{\f12(a-\g-2)}+p_0^{-1+\f{\z_1}2} \lesssim  p_0^{\f12(a-\g-2)}+p_0^{-1},
\end{align}
where we have used the fact $\z_1\leq 0$.\\

\noindent{2).} For $-2<\g<0$, it follows from $\eqref{2.13}_5$, \eqref{A.10} and \eqref{A.17} that
\begin{align}\label{2.44}
\int_{\mathbb{R}^3}k_{2a}(p,q)dq&\lesssim  p_0^{\f12(a+3|\g|)}\int_{\mathbb{R}^3}\f{e^{-\f1{16}|p-q|}}{[|p\times q|+|p-q|]^{1+|\g|}} dq\nonumber\\
&\lesssim
\begin{cases}
p_0^{\f12(a+|\g|-2)}\ln p_0,~\mbox{for}~-1\leq\g<0,\\
p_0^{\f12(a+3|\g|)-2},~\mbox{for}~-2<\g<-1.
\end{cases}
\end{align}
Finally, noting $ \g>-\f43, a\in[0,\min\{2+\g,4+3\g\}), b\in[0,2)$, if follows from \eqref{2.40}-\eqref{2.41} and \eqref{2.43}-\eqref{2.44} that
\begin{align}\label{2.46}
\int_{\mathbb{R}^3}k_{2a}(p,q)dq\lesssim p^{-\xi_{22}},
\end{align}
where $\xi_{22}\doteq \f14\Big\{1,2+\g-a,4+3\g-a\Big\}>0$. Thus \eqref{2.39} follows immediately from \eqref{2.45-1} and \eqref{2.46}.
 Therefore  the proof of Lemma \ref{lem2.4} is completed. $\hfill\Box$

\

\begin{lemma}\label{lem2.5}
Under the assumptions of Lemma \ref{lem2.4}, for any given $\a\geq0$, it holds that
\begin{align}\label{2.47}
\int_{\mathbb{R}^3} (1+|q|)^{-\a} |k(p,q)|dq\lesssim  p_0^{-\a-\f{b}2-\xi_1}\lesssim 1,~~\mbox{for soft potentials},
\end{align}
and for hard potentials,
\begin{align}\label{2.48}
\int_{\mathbb{R}^3}(1+|q|)^{-\a}  |k(p,q)|dq\lesssim p_0^{-\a-\xi_2}.
\end{align}
\end{lemma}
\noindent{\bf Proof.}  Using  \eqref{2.11-3}, a direct calculation shows that
\begin{align}\label{2.49}
\int_{\mathbb{R}^3} (1+|q|)^{-\a} k_1(p,q)dq\lesssim e^{-\f18p_0}.
\end{align}
On the other hand,  using  the same arguments as in Lemma \ref{lem2.4} and the following facts
\begin{align} \label{2.50}
(1+|q|)^{-\a} e^{-\f18|p-q|}\leq C(1+|p|)^{-\a} e^{-\f1{10}|p-q|},
\end{align}
one can obtain, for soft potentials,
\begin{equation*}\label{2.51}
\int_{\mathbb{R}^3} (1+|q|)^{-\a} |k_2(p,q)|dq\lesssim  p_0^{-\a-\f{b}2-\xi_1}\lesssim 1,
\end{equation*}
and for hard potentials,
\begin{equation*}\label{2.55}
\int_{\mathbb{R}^3}(1+|q|)^{-\a}  |k_2(p,q)|dq\lesssim p_0^{-\a-\xi_2}.
\end{equation*}
Therefore the proof of Lemma \ref{lem2.5} is completed.  $\hfill\Box$

\

Motivated by  Guo  \cite{Guo}, we have the following lemma, which will play an important role in the following a priori estimates later.
\begin{lemma}\label{lem2.6}
Let $F(t)$ satisfy \eqref{1.17}, \eqref{1.19} and the additional entropy inequality \eqref{1.20}
then it holds that
\begin{align}\label{2.52}
&\int_{\Omega}\int_{\mathbb{R}^3}\f{|F(t,x,p)-J(p)|^2}{4J(p)}I_{\{|F(t,x,p)-J(p)|< J(p)\}}dpdx\nonumber\\
&~~~~~~~~~~~~~~~~~~~~~~~~~~~~~~~+\int_{\Omega}\int_{\mathbb{R}^3}\f{1}{4}|F(t,x,p)-J(p)|I_{\{|F(t,x,p)-J(p)|\geq J(p)\}}dpdx\nonumber\\
&\leq \int_{\Omega}\int_{\mathbb{R}^3}[F_0\ln F_0-J\ln J] dpdx+[\ln(4\pi)-1] M_0 + E_0\equiv \mathcal{E}(F_0).
\end{align}
\end{lemma}
\noindent{\bf Proof.} The Taylor expansion implies that
\begin{equation*}
F(t)\ln F(t)-J\ln J=(1+\ln J)[F(t)-J]+\f{|F(t)-J|^2}{2\tilde{F}},
\end{equation*}
where $\tilde{F}$ is between $F(t)$ and $J$. Noting $\ln J=-\ln(4\pi)-p_0$, we have
\begin{align}\label{2.53}
0&\leq\int_{\Omega}\int_{\mathbb{R}^3} \f{|F(t)-J|^2}{2\tilde{F}}dqdx\nonumber\\
&= \int_{\Omega}\int_{\mathbb{R}^3} [F(t)\ln F(t)-J\ln{J}]dqdx+\int_{\Omega}\int_{\mathbb{R}^3} \{[\ln(4\pi)-1]+p_0\}[F(t)-J]dqdx\nonumber\\
&=  \int_{\Omega}\int_{\mathbb{R}^3} [F(t)\ln F(t)-J\ln{J}]dqdx+ [\ln(4\pi)-1]M_0+E_0\nonumber\\
&\leq \int_{\Omega}\int_{\mathbb{R}^3} [F_0\ln F_0-J\ln{J}]dqdx+ [\ln(4\pi)-1]M_0+E_0,
\end{align}
where we have used \eqref{1.17}, \eqref{1.19} and \eqref{1.20} above.
Noting that $|F-J|\geq J$ yields that $F\geq 2J$ or $F=0$, thus we have  that
\begin{equation*}
\f{|F-J|}{\tilde{F}}\geq \f12,
\end{equation*}
which, together with \eqref{2.53},  yields that
\begin{align}
&\int_{\Omega}\int_{\mathbb{R}^3}\f{|F(t,x,p)-J(p)|^2}{4J(p)}I_{\{|F(t,x,p)-J(p)|< J(p)\}}dpdx\nonumber\\
&~~~~~~~~~~~~~~~~~~~~~~~~~+\int_{\Omega}\int_{\mathbb{R}^3}\f{1}{4}|F(t,x,p)-J(p)|I_{\{|F(t,x,p)-J(p)|\geq J(p)\}}dpdx\nonumber\\
&\leq \int_{\Omega}\int_{\mathbb{R}^3} [F_0\ln F_0-J\ln{J}]dqdx+ [\ln(4\pi)-1] M_0 + E_0.\nonumber
\end{align}
Therefore the proof of Lemma \ref{lem2.6} is completed. $\hfill\Box$


\section{Local Existence Result}
As mentioned previously, Bichteler \cite{Bichteler} proved the local existence and uniqueness  in the $L^\infty$ framework under smallness conditions on the initial data. To prove Theorem \ref{thm1.1}, firstly,  we need to  establish the local existence of unique solutions to the relativistic Boltzmann equation \eqref{1.1} with general bounded initial data in $L^\infty$ space. We point out that the Lorentz transformation is essentially used in this section.

\begin{theorem}[Local Existence]\label{thm7.1} We assume that $-2<\g,~a\in[0,2]\cap [0,2+\g),b\in[0,\min\{4,4+\g\})$ for hard potentials,  and $-2<\g,b\in(0,\min\{4,4+\g\})$ for soft potentials.  Let $\Omega=\mathbb{T}^3$ or $\mathbb{R}^3$, $\b>14$,  $F_0(x,p)=J(p)+\sqrt{J(p)}f_0(x,p)\geq 0$ and $\|w_\b f_0\|_{L^\infty}<\infty$, then there exists a positive time
	\begin{align}\label{LT}
	t_1:=(8\tilde{C}_4[1+\|w_\b f_0\|_{L^\infty}])^{-1}>0,
	\end{align}
	such that the relativistic  Boltzmann equation \eqref{1.5}, \eqref{1.5-1} has a unique mild solution $F(t,x,p)=J(p)+\sqrt{J(p)}f(t,x,p)\geq 0$ on the time interval $t\in[0,t_1]$ and  satisfies
	\begin{align}\label{7.1}
	\|w_\b f(t)\|_{L^\infty}\leq 2\|w_\b f_0\|_{L^\infty},~\mbox{for}~0\leq t\leq t_1,
	\end{align}
	where the positive constant $\tilde{C}_4\geq1$ depends only on $a,b,\g,\b$. In addition,  the conservations of mass, momentum, energy \eqref{1.17}-\eqref{1.19} as well as the additional entropy inequality \eqref{1.20} hold.
	Furthermore, if the initial data $f_0$ is continuous, then the solution $f(t,x,p)$ is continuous in $[0,t_1]\times\Omega\times\mathbb{R}^3$.
\end{theorem}
\noindent{\bf Proof.} To prove the local existence of the relativistic Boltzmann equation \eqref{1.5}, \eqref{1.5-1}, we consider the iteration, for $n=0,1,2,\cdots$,
\begin{align}\label{7.2}
& F^{n+1}_t+\hat{p}\cdot\nabla_x F^{n+1}+F^{n+1}\cdot\int_{\mathbb{R}^3}\int_{\mathbb{S}^2}v_{\phi}\s(g,\t)F^{n}(t,x,q)d\omega dq=Q_+(F^n,F^n),
\end{align}
with
\begin{align}\label{7.3}
F^{n+1}(t,x,p)\Big|_{t=0}=F_0(x,p)\geq0~~~\mbox{and}~~F^0(t,x,p)=J(p).
\end{align}
Denote
\begin{equation*}\label{7.4}
f^{n+1}(t,x,p):=\f{F^{n+1}(t,x,p)-J(p)}{\sqrt{J(p)}},
\end{equation*}
then \eqref{7.2} can be written equivalently as
\begin{align}\label{7.5}
&f^{n+1}_t+\hat{p}\cdot\nabla_x f^{n+1}+f^{n+1}\cdot \int_{\mathbb{R}^3}\int_{\mathbb{S}^2}v_{\phi}\s(g,\t)\Big\{J(q)+\sqrt{J(q)}f^n(t,x,q)\Big\}d\omega dq\nonumber\\
&~~~~~~~~~~~~~~~~~~=Kf^n+\f{1}{\sqrt{J(p)}}Q_+(\sqrt{J}f^n,\sqrt{J}f^n)
\end{align}
with $n=0,1,2,\cdots$ and
\begin{equation*}\label{7.6}
f^{n+1}(t,x,p)\Big|_{t=0}=f_0(x,p)~~~\mbox{and}~~~f^0(t,x,p)=0.
\end{equation*}
Hence we get an approximation sequence $F^{n+1},n=0,1,\cdots$ by solving the linear equation \eqref{7.2}, \eqref{7.3}.

Firstly, we consider the positivity of $F^{n+1}$.  It is noted that
\begin{align}\label{7.7}
F^{n+1}(t,x,p)&=e^{-\int_0^tB^n(\tau,x-\hat{p}(t-\tau),p)d\tau}\cdot F_0(x-\hat{p}t,p)\nonumber\\
&~~+\int_0^te^{-\int_s^tB^n(\tau,x-\hat{p}(t-\tau),p)d\tau}\cdot Q_+(F^n,F^n)(s,x-\hat{p}(t-s),p)ds,
\end{align}
where
\begin{equation}\label{7.8}
B^n(\tau,y,p):=\int_{\mathbb{R}^3}\int_{\mathbb{S}^2}v_{\phi}\s(g,\t)F^n(\tau,y,q)d\omega dq=\int_{\mathbb{R}^3}\int_{\mathbb{S}^2}v_{\phi}\s(g,\t)\Big\{J(q)+\sqrt{J}f^n(\tau,y,q)\Big\}d\omega dq.
\end{equation}
By the induction argument, if $F^n\geq0$, then it holds that
$$B^n(\tau,x-\hat{p}(t-\tau),p)\geq0~~\mbox{and}~~Q_+(F^n,F^n)(\tau,x-\hat{p}(t-\tau),p)\geq0,$$
which, together with \eqref{7.7},  yields immediately  that
\begin{equation*}\label{7.9}
F^{n+1}(t,x,p)\geq e^{-\int_0^tB^n(\tau,x-\hat{p}(t-\tau),p)d\tau}\cdot F_0(x-\hat{p}t,p)\geq 0.
\end{equation*}
Hence we have shown $F^{n+1}\geq0$ for any $n=0,1,\cdots$.

Next we consider the  uniform $L^\infty$-estimate for the above approximation sequence. For this, it is more convenient to use the equivalent form $f^{n+1}$. Indeed, it follows from \eqref{7.5} that
\begin{align}\label{7.10}
f^{n+1}(t,x,p)
&=e^{-\int_0^tB^n(\tau,x-\hat{p}(t-\tau),p)d\tau}\cdot f_0(x-\hat{p}t,p) \nonumber\\
&~~+\int_0^te^{-\int_s^tB^n(\tau,x-\hat{p}(t-\tau),p)d\tau}\cdot (Kf^n)(s,x-\hat{p}(t-s),p)ds\\
&~~+\int_0^te^{-\int_s^tB^n(\tau,x-\hat{p}(t-\tau),p)d\tau}\cdot\f{1}{\sqrt{J(p)}}Q_+(\sqrt{J}f^n,\sqrt{J}f^n)(s,x-\hat{p}(t-s),p)ds,\nonumber
\end{align}
which yields immediately that
\begin{align}\label{7.11}
|w_\b(p) f^{n+1}(t,x,p)|&\leq \|w_\b f_0\|_{L^\infty} +\int_0^t\f{w_\b(p)}{\sqrt{J(p)}}\Big|Q_+(\sqrt{J}f^n,\sqrt{J}f^n)(s,\hat{x},p)\Big|ds
\nonumber\\
&~~~~+\int_0^t\Big|w_\b(p) (Kf^n)(s,\hat{x},p)\Big| ds
\end{align}
where have denoted $\hat{x}:= x-\hat{p}(t-s)$.\\

To estimate the second term on the RHS of \eqref{7.11},
we  define
\begin{align}\label{7.15-2}
&\f{w_\b(p)}{\sqrt{J(p)}}Q_+^b(\sqrt{J}|f^n|,\sqrt{J}|f^n|)(s,\hat{x},p):=w_\b(p) \int_{\mathbb{R}^3\times\mathbb{S}^2}v_\phi \s_b(g,\t)\sqrt{J(q)}|f^n(s,\hat{x},p')f^n(s,\hat{x},q')|d\omega dq,
\end{align}
and
\begin{align}\label{7.15-3}
&\f{w_\b(p)}{\sqrt{J(p)}}Q_+^a(\sqrt{J}|f^n|,\sqrt{J}|f^n|)(s,\hat{x},p):=w_\b(p) \int_{\mathbb{R}^3\times\mathbb{S}^2}v_\phi \s_a(g,\t)\sqrt{J(q)}|f^n(s,\hat{x},p')f^n(s,\hat{x},q')|d\omega dq,
\end{align}
where $\s_b(\cdot,\cdot)$ and $\s_a(\cdot,\cdot)$ are defined in \eqref{2.15-1} above.
Therefore, to estimate  the second term on the RHS of \eqref{7.11}, we need only to estimate \eqref{7.15-2} and \eqref{7.15-3}.
Using \eqref{A.25} with $d=1$, one has
\begin{align}\label{7.15-4}
\f{w_\b(p)}{\sqrt{J(p)}}Q_+^b(\sqrt{J}|f^n|,\sqrt{J}|f^n|)(s,\hat{x},p)
&\lesssim \|w_\b f^n\|^2_{L^\infty} \int_{\mathbb{R}^3\times\mathbb{S}^2}v_\phi \s_b(g,\t)\sqrt{J(q)}d\omega dq \nonumber\\
&\lesssim p_0^{-\f{b}2}\|w_\b f^n\|^2_{L^\infty}
\end{align}
On the other hand, using
$ p_0\leq p_0'+q_0'$, one has
\begin{equation*}\label{3.10}
\mbox{either}~\f12 p_0\leq p_0'~~\mbox{or}~~\f12p_0\leq q_0',
\end{equation*}
which yields
\begin{align}\label{7.15}
I_a&:= \f{w_\b(p)}{\sqrt{J(p)}}Q_+^a(\sqrt{J}|f^n|,\sqrt{J}|f^n|)(s,\hat{x},p)\nonumber\\
&\equiv\f{w_\b(p)}{2p_0} \int_{\mathbb{R}^3\times\mathbb{R}^3\times\mathbb{R}^3} s\s_a(g,\t)\d^{(4)}(p^\mu+q^\mu-{p^\mu}'-{q^\mu}')\sqrt{J(q)}\Big|f^n(s,\hat{x},p')f^n(s,\hat{x},q')\Big|\f{dp'dq'dq}{p_0'q_0'q_0}\nonumber\\
&\leq   \f{C}{p_0}\int_{\mathbb{R}^3\times\mathbb{R}^3\times\mathbb{R}^3}s\s_a(g,\t)\d^{(4)}(p^\mu+q^\mu-{p^\mu}'-{q^\mu}')\sqrt{J(q)}\Big|w_\b(q')f^n(s,\hat{x},q')f^n(s,\hat{x},p')
\Big| \f{dp'dq'dq}{p_0'q_0'q_0}\nonumber\\
&~+\f{C}{p_0}\int_{\mathbb{R}^3\times\mathbb{R}^3\times\mathbb{R}^3}s\s_a(g,\t)\d^{(4)}(p^\mu+q^\mu-{p^\mu}'-{q^\mu}')\sqrt{J(q)}\Big|f^n(s,\hat{x},q')w_\b(p')f^n(s,\hat{x},p')
\Big|\f{dp'dq'dq}{p_0'q_0'q_0}\nonumber\\
&:= I_{a1}+I_{a2}.
\end{align}
For $I_{a1}$, we first exchange  $p'$ and $q$ to get that
\begin{align}\label{7.16}
I_{a1}&= \f{C}{p_0}\int_{\mathbb{R}^3\times\mathbb{R}^3\times\mathbb{R}^3}
\bar{s}\s_a(\bar{g},\bar{\t})\d^{(4)}(p^\mu+{p^\mu}'-q^\mu-{q^\mu}')\nonumber\\
&~~~~~~~~~~~~~~~~~~~~~~~~~~\times\sqrt{J(p')}\Big|w_\b(q')f^n(s,\hat{x},q')f^n(s,\hat{x},q)
\Big| \f{dp'dq'dq}{p_0'q_0'q_0},
\end{align}
where
\begin{align}\label{7.17}
\begin{cases}
\bar{g}^2:= |g(p^\mu,{p^\mu}')|^2\equiv g^2+\f12(p^\mu+q^\mu)(p_\mu+q_\mu-p_\mu'-q_\mu'),\\
\bar{s}=4+\bar{g}^2,~~\cos\bar{\t}=1-2\Big(\f{g}{\bar{g}}\Big)^2.
\end{cases}
\end{align}
For $I_{a2}$, we first exchange $q'$ and $p'$, then $p'$ and  $q$ to obtain
\begin{align}\label{7.18}
I_{a2}&= \f{C}{p_0}\int_{\mathbb{R}^3\times\mathbb{R}^3\times\mathbb{R}^3}
\bar{s}\s_a(\bar{g},\bar{\t})\d^{(4)}(p^\mu+{p^\mu}'-q^\mu-{q^\mu}')\nonumber\\
&~~~~~~~~~~~~~~~~~~~~~~~~~~~~~~~~~~~\times\sqrt{J(p')}\Big|w_\b(q')f^n(s,\hat{x},q')f^n(s,\hat{x},q)
\Big| \f{dp'dq'dq}{p_0'q_0'q_0}.
\end{align}
Combining \eqref{7.15}-\eqref{7.18}, one gets that
\begin{align}\label{7.19}
I_a\leq C\|w_\b f^n(s)\|_{L^\infty}\int_{\mathbb{R}^3} A_a(p,q) |f(s,\hat{x},q)| dq
\leq C\|w_\b f^n(s)\|^2_{L^\infty}\int_{\mathbb{R}^3} A_a(p,q) q_0^{-\b}dq,
\end{align}
where
\begin{equation}\label{7.19-1}
A_a(p,q):=\f{1}{p_0q_0}\int_{\mathbb{R}^3\times\mathbb{R}^3}
\bar{s}\s_a(\bar{g},\bar{\t})\d^{(4)}(p^\mu+{p^\mu}'-q^\mu-{q^\mu}')\sqrt{J(p')}\f{dp'dq'}{p_0'q_0'}.
\end{equation}
Now we estimate $A_a(p,q)$. Define
$$u=u(r)=\begin{cases}
0,~~\mbox{if}~~r<0,\\
1,~~\mbox{if}~~r\geq0.
\end{cases}$$
Let $\underline{g}:=g({p^\mu}',{q^\mu}')$ and $\underline{s}:=s({p^\mu}',{q^\mu}')$, then
the  following  identity holds
\begin{align}\label{7.20}
\int_{\mathbb{R}^3\times\mathbb{R}^3} G(p^\mu,q^\mu|{p^\mu}',{q^\mu}')\f{dp'dq'}{p_0'q_0'}
=16\int_{\mathbb{R}^4\times\mathbb{R}^4}G(p^\mu,q^\mu|{p^\mu}',{q^\mu}')d\Theta({p^\mu}',{q^\mu}'),
\end{align}
with
\begin{equation*}\label{7.21}
d\Theta({p^\mu}',{q^\mu}')\doteq u(p_0'+q_0')u(\underline{s}-4)\d(\underline{s}-\underline{g}^2-4)\d(({p^\mu}'+{q^\mu}')({p_\mu}'-{q_\mu}'))d{p^\mu}'d {q^\mu}',
\end{equation*}
whose proof can be found in \cite{Strain3,Strain,Groot}. Combining \eqref{7.20} and \eqref{7.19-1}, we obtain
\begin{align}\label{7.22}
A_a(p,q)=\f{16}{p_0q_0}\int_{\mathbb{R}^4\times\mathbb{R}^4}
\bar{s}\s_a(\bar{g},\bar{\t})\d^{(4)}(p^\mu+{p^\mu}'-q^\mu-{q^\mu}')\sqrt{J(p')}d\Theta({p^\mu}',{q^\mu}').
\end{align}
We consider the change of variables
\begin{equation*}\label{7.23}
\bar{p}^{\mu}={p^\mu}'+{q^\mu}',~~~\bar{q}^{\mu}={p^\mu}'-{q^\mu}',
\end{equation*}
which yields immediately that
\begin{equation*}\label{7.24}
{p^\mu}'=\f12(\bar{p}^{\mu}+\bar{q}^{\mu}),~~{q^\mu}'=\f12(\bar{p}^{\mu}-\bar{q}^{\mu}).
\end{equation*}
Applying the above change of variables, we have that
\begin{equation*}\label{7.28}
A_a(p,q)
=\f{1}{p_0q_0}\int_{\mathbb{R}^4\times\mathbb{R}^4}
\bar{s}\s_a(\bar{g},\bar{\t})\d^{(4)}(p^\mu-q^\mu+\bar{q}^\mu) e^{-\f14(\bar{p}_0+\bar{q}_0)} d\Theta(\bar{p}^\mu,
\bar{q}^\mu),
\end{equation*}
with $ d\Theta(\bar{p}^\mu,
\bar{q}^\mu):= u(\bar{p}_0)u(-\bar{p}^{\mu}\bar{p}_{\mu}-4)\d(-\bar{p}^{\mu}\bar{p}_{\mu}-\bar{q}^{\mu}\bar{q}_{\mu}-4)\d(\bar{p}^{\mu}\bar{q}_{\mu})d\bar{p}^\mu d \bar{q}^\mu$.
Here we have used $\sqrt{J(p')}=e^{-\f12{p_0}'}=e^{-\f14({\bar{p}_0}+\bar{q}_0)}$. Also, $\bar{g}\geq0$ is now given by
\begin{equation*}\label{7.29}
\bar{g}^2=g^2+\f12(p^\mu+q^\mu)(p_\mu+q_\mu-\bar{p}_\mu),~~\bar{s}=4+\bar{g}^2,~~\cos\bar{\t}=1-2\Big(\f{g}{\bar{g}}\Big)^2.
\end{equation*}
Now we calculate the delta function for $\d^{(4)}(p^\mu-q^\mu+\bar{q}^\mu)$ to get that
\begin{align}\label{7.30}
A_a(p,q)&=\f{e^{-\f14(q_0-p_0)}}{p_0q_0}\int_{\mathbb{R}^4}
\bar{s}\s_a(\bar{g},\bar{\t}) e^{-\f14\bar{p}_0} u(\bar{p}_0)u(-\bar{p}^{\mu}\bar{p}_{\mu}-4)\d(-\bar{p}^{\mu}\bar{p}_{\mu}-g^2-4)\d(\bar{p}^{\mu}(q_{\mu}-p_{\mu}))d\bar{p}^\mu \nonumber\\
&~~~\mbox{Noting} ~~u(\bar{p}_0)\d(-\bar{p}^{\mu}\bar{p}_{\mu}-g^2-4)=
\f{\d(\bar{p}_0-\sqrt{s+|\bar{p}|^2})}{2\sqrt{s+|\bar{p}|^2}},~\mbox{we have}\nonumber\\
&=\f{e^{-\f14(q_0-p_0)}}{p_0q_0}\int_{\mathbb{R}^4}
\bar{s}\s_a(\bar{g},\bar{\t}) e^{-\f14\bar{p}_0} u(-\bar{p}^{\mu}\bar{p}_{\mu}-4)\f{\d(\bar{p}_0-\sqrt{s+|\bar{p}|^2})}{2\sqrt{s+|\bar{p}|^2}}\d(\bar{p}^{\mu}(q_{\mu}-p_{\mu}))d\bar{p}^\mu \nonumber\\
&=\f{e^{-\f14(q_0-p_0)}}{2p_0q_0}\int_{\mathbb{R}^3}
\bar{s}\s_a(\bar{g},\bar{\t})  \d(\bar{p}^{\mu}(q_{\mu}-p_{\mu}))e^{-\f14\bar{p}^\mu\bar{U}_\mu}\f{d\bar{p}}{\bar{p}_0},
\end{align}
where we have used the notations $\bar{p}_0=\sqrt{s+|\bar{p}|^2}$, $\bar{U}^\mu=(-1,0,0,0)^t$ and the fact $-\bar{p}^{\mu}\bar{p}_{\mu}-4=s-4\geq0$. We introduce a Lorentz transformation $\Lambda=\Big(\Lambda^{\mu\nu}\Big)$ such that
\begin{align}
A_\nu=\Lambda^{\mu\nu}(p_\mu+q_{\mu})=(\sqrt{s},0,0,0),~~B_\nu=-\Lambda^{\mu\nu}(p_\mu-q_{\mu})=(0,0,0,g).\nonumber
\end{align}
Indeed, Strain  \cite{Strain3} gives  details of the Lorentz transformation
\begin{align}
\Lambda=\Big(\Lambda^{\mu\nu}\Big)
=
\left(\begin{array}{cccc}
\f{p_0+q_0}{\sqrt{s}},&-\f{p_1+q_1}{\sqrt{s}},&-\f{p_2+q_2}{\sqrt{s}},&-\f{p_3+q_3}{\sqrt{s}}\\[2mm]
\f{2|p\times q|}{g\sqrt{s}},&\Lambda^{11},&\Lambda^{21},&\Lambda^{31}\\[2mm]
0,&\f{(p\times q)_1}{|p\times q|},&\f{(p\times q)_2}{|p\times q|},&\f{(p\times q)_3}{|p\times q|}\\[2mm]
\f{p_0-q_0}{g},&-\f{p_1-q_1}{g},&-\f{p_2-q_2}{g},&-\f{p_3-q_3}{g}
\end{array} \right),\nonumber
\end{align}
where $\Lambda^{i1},~i=1,2,3$ can be found in \cite{Strain3}, we omit the details here.

Define $U^\mu:=\Lambda^{\nu\mu}\bar{U}_{\nu}$,  we have
\begin{equation*}\label{7.33}
U^\mu=\Big(-\f{p_0+q_0}{\sqrt{s}},-\f{2|p\times q|}{g\sqrt{s}},~0,~-\f{p_0-q_0}{g}\Big).
\end{equation*}
Using the above Lorentz transformation, one can get that
\begin{align}\label{7.34}
\int_{\mathbb{R}^3}
\bar{s}\s_a(\bar{g},\bar{\t})  \d(\bar{p}^{\mu}(q_{\mu}-p_{\mu}))e^{-\f14\bar{p}^\mu\bar{U}_\mu}\f{d\bar{p}}{\bar{p}_0}
=\int_{\mathbb{R}^3}
\bar{s}_{\Lambda}\s_a(\bar{g}_{\Lambda},\bar{\t}_{\Lambda})  \d(\bar{p}^{\mu}B_\mu)e^{-\f14\bar{p}^\mu U_\mu}\f{d\bar{p}}{\bar{p}_0},
\end{align}
where we have used $\bar{p}^{\mu}$ and $\f{d \bar{p}}{\bar{p}_0}$ are Lorentz invariants. Here $\bar{g}_{\Lambda},~\bar{s}_{\Lambda}\geq0,~\bar{\theta}_{\Lambda}$ are given by
\begin{align}\label{7.35}
\bar{g}_{\Lambda}^2=g^2+\f12A^\mu(A_\mu-\bar{p}_\mu)=g^2+\f12\sqrt{s}(\bar{p}_0-\sqrt{s}),~~\bar{s}_{\Lambda}=4+\bar{g}_{\Lambda}^2,~~\cos\bar{\t}_{\Lambda}=1-2\Big(\f{g}{\bar{g}_{\Lambda}}\Big)^2.
\end{align}
To calculate \eqref{7.34}, we use the polar coordinate
\begin{equation*}
d\bar{p}=|\bar{p}|^2 \sin\psi d|\bar p|d\psi d\varphi,~~\bar{p}=|\bar p|(\sin\psi\cos\varphi,\sin\psi\sin\varphi,\cos\psi),
\end{equation*}
which, together with the fact $\bar{p}^\mu B_{\mu}=\bar{p}_3 g$, yields that
\begin{align}
A_a(p,q)&=\f{e^{-\f14(q_0-p_0)}}{2p_0q_0} \int_{0}^{2\pi}d\varphi\int_0^{\pi} \sin\psi d\psi
\int_0^{\infty}\bar{s}_{\Lambda}\s_a(\bar{g}_{\Lambda},\bar{\t}_{\Lambda})  \d(|\bar{p}|g\cos\psi)e^{-\f14\bar{p}^\mu U_\mu}\f{|\bar{p}|^2d|\bar{p}|}{\bar{p}_0}\nonumber\\
&=\f{e^{-\f14(q_0-p_0)}}{2gp_0q_0} \int_{0}^{2\pi}d\varphi
\int_0^{\infty}\bar{s}_{\Lambda}\s_a(\bar{g}_{\Lambda},\bar{\t}_{\Lambda})  e^{-\bar{p}_0\f{p_0+q_0}{4\sqrt{s}}} e^{\f{|p\times q|}{2g\sqrt{s}}|\bar{p}|\cos\varphi} \f{|\bar{p}|d|\bar{p}|}{\bar{p}_0}\nonumber\\
&=\f{\pi e^{-\f14(q_0-p_0)}}{gp_0q_0}
\int_0^{\infty}\bar{s}_{\Lambda}\s_a(\bar{g}_{\Lambda},\bar{\t}_{\Lambda})  e^{-\bar{p}_0\f{p_0+q_0}{4\sqrt{s}}} I_0\Big(\f{|p\times q|}{2g\sqrt{s}}|\bar{p}|\Big) \f{|\bar{p}|d|\bar{p}|}{\bar{p}_0}.\nonumber
\end{align}
Denoting $z\doteq|\bar{p}|$, it follows from \eqref{7.35} that
\begin{equation*}
g=\bar{g}_{\Lambda}\sqrt{\f{1-\cos\bar{\t}_{\Lambda}}{2}}=\bar{g}_{\Lambda}\sin\f{\bar{\t}_{\Lambda}}{2},~~\bar{s}_{\Lambda}=\f12s+\f12s\sqrt{z^2/s+1},
\end{equation*}
with
\begin{equation*}
\sin\f{\bar{\t}_{\Lambda}}{2}=\f{g}{\bar{g}_{\Lambda}}=\f{\sqrt2 g}{\sqrt{g^2-4+s\sqrt{z^2/s+1}}}.
\end{equation*}
Noting $\bar{p}_0=\sqrt{s+z^2}$, one gets
\begin{align}\label{7.38}
A_a(p,q)&=\f{\pi  e^{-\f14(q_0-p_0)}s}{2gp_0q_0}
\int_0^{\infty}[1+\sqrt{z^2/s+1}]\s_a\Big(\f{g}{\sin\f{\bar{\t}_{\Lambda}}{2}},\bar{\t}_{\Lambda}\Big)  e^{-\f{p_0+q_0}{4\sqrt{s}}\sqrt{s+z^2}} I_0\Big(\f{|p\times q|}{2g\sqrt{s}}z\Big) \f{zdz}{\sqrt{s+z^2}}\nonumber\\
&=\f{\pi  e^{-\f14(q_0-p_0)}s^{\f32}}{2gp_0q_0} \int_{0}^{\infty} \f{y[1+\sqrt{y^2+1}]}{\sqrt{1+y^2}}\s_a\Big(\f{g}{\sin\f{\psi}{2}},\psi\Big)e^{-\f{l}{2}\sqrt{1+y^2}} I_0\Big(\f{j}{2}y\Big)dy,
\end{align}
where we have made the change of variable $y=\f{z}{\sqrt{s}}$ and defined
\begin{equation*}\label{7.39}
\sin\f{\psi}{2}:=\f{\sqrt2 g}{\sqrt{g^2-4+(g^2+4)\sqrt{y^2+1}}}.
\end{equation*}

It follows from $\eqref{2.7}$ that $s\lesssim p_0q_0$, then
by the same arguments as in \eqref{2.19-1}, \eqref{2.22-1}, \eqref{2.25-1} and using \eqref{2.7}, \eqref{2.9}, one has
\begin{align}\label{7.55}
A_a(p,q)\lesssim
\begin{cases}
\f{(p_0q_0)^{1+a-\f\g2}}{(1+|p-q|)^{2+\f{a-\g}{2}}}e^{-\f14|p-q|}e^{-\f14(q_0-p_0)},
~~a\geq\g\geq1,\\[1.5mm]
\f{(p_0q_0)^{2+a-\f{3\g}2}}{|p-q|^{1-\g}(1+|p-q|)^{2+\f{a-\g}{2}}}e^{-\f14|p-q|}e^{-\f14(q_0-p_0)},~~a\geq\g,~0\leq \g<1,\\[3mm]
\f{(p_0q_0)^{1+\f{a}{2}+\f{\z_1}2}}{(1+|p-q|)^{2+\f{\z_1}{2}}}e^{-\f14|p-q|}e^{-\f14(q_0-p_0)},
~~1\leq a<\g, \\[3mm]
\f{(p_0q_0)^{2-\f{a}{2}+\f{\z_1}2}}{|p-q|^{1-a}(1+|p-q|)^{2+\f{\z_1}{2}}}e^{-\f14|p-q|}e^{-\f14(q_0-p_0)},~~a<\g,~0\leq a<1,\\[3.5mm]
\Big[\f{(p_0q_0)^{1+|\g|+\f{a}{2}}e^{-\f14|p-q|}}{|p-q|^{1+|\g|}(1+|p-q|)^{1-\f{|\g|}{2}-O(\v)}}+
\f{(p_0q_0)^{2+a+\f{3|\g|}{2}}e^{-\f14|p-q|}}{|p-q|^{1+|\g|}(1+|p-q|)^{2+\f{a+|\g|}{2}}}\Big]e^{-\f14(q_0-p_0)},-2<\g<0,
\end{cases}
\end{align}
where we have used the fact $p_0+q_0\leq p_0q_0$.
Noting that
\begin{equation*}\label{7.56}
e^{-\f14|p-q|}e^{-\f14(q_0-p_0)}\leq 1,
\end{equation*}
which, together with \eqref{7.55} and \eqref{A.18}, yields, for $\b>14$, that
\begin{equation}\label{7.57}
\int_{\mathbb{R}^3}A_a(p,q)  q_0^{-\b}dq \lesssim
\begin{cases}
p_0^{-1+\f{a}2}\lesssim 1,~\g\geq0,~0\leq a \leq 2,\\
p_0^{-1+\f{a}2}+p_0^{\f12(a+|\g|-2)+O(\v)}\lesssim 1, ~-2<\g<0,~0\leq  a < \min\{2+\g\}.
\end{cases}
\end{equation}
where we have chosen $\v>0$ sufficiently small so that $\f12(a-2-\g)+O(\v)\leq 0$.
It is here that we need the condition $a\leq 2$ and $a<2+\g$.
Substituting \eqref{7.57} into \eqref{7.19},  one obtains
\begin{align}
I_a &=\f{w_\b(p)}{\sqrt{J(p)}}Q_+^a(\sqrt{J}|f^n|,\sqrt{J}|f^n|)(s,\hat{x},p)
\leq C\|w_\b f^n(s)\|^2_{L^\infty},\nonumber
\end{align}
which, together with \eqref{7.15-4}, yields that
\begin{align}\label{7.62-1}
\f{w_\b(p)}{\sqrt{J(p)}}\Big|Q_+(\sqrt{J}f^n,\sqrt{J}f^n)(s,\hat{x},p)\Big|
\leq C\|w_\b f^n(s)\|^2_{L^\infty}.
\end{align}
We remark that Lemmas \ref{lemA.2} and \ref{lemA.3} are very important for us to bound $A_a(p,q)$. Indeed, if one use the corresponding lemmas of \cite{Glassey1}, it will be hard to control the part $e^{-\f14(q_0-p_0)}$.

\

For the last term in the RHS of \eqref{7.11}, we note
\begin{align}\label{7.13-1}
I_K&:= \Big|w_\b(p) (Kf^n)(s,\hat{x},p)\Big|\leq  \Big|w_\b(p) (K_1f^n)(s,\hat{x},p)\Big|+ \Big|w_\b(p) (K_2f^n)(s,\hat{x},p)\Big|\nonumber\\
&=I_{K1}+I_{K2}.
\end{align}
It follows from  \eqref{2.49}, for $\b\geq0$, that
\begin{align}\label{7.13-2}
I_{K1}\leq \|w_\b f^n(s)\|_{L^\infty}\cdot w_{\b}(p)\int_{\mathbb{R}^3}(1+|q|)^{-\b} k_1(p,q)dq
\lesssim \|w_\b f^n(s)\|_{L^\infty} .
\end{align}
Noting the definition of $K_2f$ in \eqref{1.25}, by the similar arguments as in  \eqref{7.15-2}-\eqref{7.62-1}, we can obtain that
\begin{equation}\label{7.13-7}
I_{K2}\leq C \|w_\b f(s)\|_{L^\infty},~\mbox{for}~\b\geq0,
\end{equation}
which, together with \eqref{7.13-2}, yields that
\begin{align}\label{7.13}
I_K\leq C\|w_\b f(s)\|_{L^\infty},~\mbox{for}~\b\geq0.
\end{align}
Substituting \eqref{7.13} and \eqref{7.62-1} into \eqref{7.11}, we have, for $\b>14$, that
\begin{equation}\label{7.63}
\|w_{\b}(p)f^{n+1}(t)\|_{L^\infty}\leq \|w_{\b}f_0\|_{L^\infty}+C_1t\Big\{\sup_{0\leq s\leq t}\|w_{\b}f^n(s)\|_{L^\infty}+\sup_{0\leq s\leq t}\|w_{\b}f^n(s)\|^2_{L^\infty}\Big\},
\end{equation}
where $C_1\geq1$ depends only on $a,b,\g,\b$.  By the induction argument, we can prove that if
\begin{align}\label{7.64}
\sup_{0\leq s\leq t_1}\|w_{\b}f^n(s)\|_{L^\infty}\leq 2\|w_{\b}f_0\|_{L^\infty}~~\mbox{with}~~t_1=\Big(8C_1[1+\|w_{\b}f_0\|_{L^\infty}]\Big)^{-1},
\end{align}
then it follows from \eqref{7.63}  and \eqref{7.64}, for $\b>14$, that
\begin{equation}\label{7.65}
\sup_{0\leq s\leq t_1}\|w_{\b}f^{n+1}(s)\|_{L^\infty}\leq 2\|w_{\b}f_0\|_{L^\infty}~\mbox{with}~t_1=\Big(8C_1[1+\|w_{\b}f_0\|_{L^\infty}]\Big)^{-1},~\mbox{for}~n\geq0.
\end{equation}

\

Now we show that $f^{n+1},~n=0,1,2,\cdots$ is a Cauchy sequence.  It follows from \eqref{7.10} that
\begin{align}\label{7.66}
&\Big|\sqrt{w_{\b}(p)}(f^{n+2}-f^{n+1})(t,x,p)]\Big|\nonumber\\
&\leq |\sqrt{w_{\b}(p)}f_0(x-\hat{p}t,p)|\cdot\int_0^t\Big|(B^{n+1}-B^n)(\tau,x-\hat{p}(t-\tau),p)\Big|d\tau\nonumber\\
&~~~~+\int_0^t|\sqrt{w_{\b}(p)}Kf^{n+1}(s,x-\hat{p}(t-s),p)|\cdot\int_s^t\Big|(B^{n+1}-B^n)(\tau,x-\hat{p}(t-\tau),p)\Big|d\tau ds\nonumber\\
&~~~~+\int_0^t\f{\sqrt{w_{\b}(p)}}{\sqrt{J(p)}}|Q_+(\sqrt{J}f^{n+1},\sqrt{J}f^{n+1})(s,x-\hat{p}(t-s),p)|\nonumber\\
&~~~~~~~~~~~~~~~~~~~~~~~~~~~~~~~~~~~~~~~~~~~\times\int_s^t\Big|(B^{n+1}-B^n)(\tau,x-\hat{p}(t-\tau),p)\Big|d\tau ds\nonumber\\
&~~~~+\int_0^t|\sqrt{w_{\b}(p)}(Kf^{n+1}-Kf^n)(s,x-\hat{p}(t-s),p)|ds\nonumber\\
&~~~~+\int_0^t\f{\sqrt{w_{\b}(p)}}{\sqrt{J(p)}}\Big|\Big(Q_+(\sqrt{\mu}f^{n+1},\sqrt{\mu}f^{n+1})-Q_+(\sqrt{J}f^{n},\sqrt{J}f^{n})\Big)(s,x-\hat{p}(t-s),p)\Big|ds\nonumber\\
&\triangleq I_1+I_2+I_3+I_4+I_5.
\end{align}
A direct calculation shows that
\begin{align}\label{7.66-1}
\Big|(B^{n+1}-B^n)(\tau,z,p)\Big|
\leq C\nu(p)\|(f^{n+1}-f^{n})(\tau)\|_{L^\infty},
\end{align}
which, together with \eqref{7.13}, \eqref{7.62-1} an \eqref{n1.20-1}, yields  
\begin{align}\label{7.67}
&I_1+I_2+I_3\leq C\Big\{ \|\nu\sqrt{w_{\b}} f_0\|_{L^\infty}+\int_0^t\Big|\sqrt{w_{\b}(p)}\nu(p)Kf^{n+1}(s,x-\hat{p}(t-s),p)\Big|ds\nonumber\\
&~~~+\int_0^t\f{\sqrt{w_{\b}(p)}\nu(p)}{\sqrt{J(p)}}\Big|Q_+(\sqrt{J}f^{n+1},\sqrt{J}f^{n+1})(s,x-\hat{p}(t-s),p)\Big|ds\Big\}\int_0^t\|(f^{n+1}-f^{n})(\tau)\|_{L^\infty}d\tau\nonumber\\
&\leq Ct\Big\{  \|w_{\b}f_0\|_{L^\infty}+t\sup_{0\leq s\leq t}\|w_{\b}f^{n+1}(s)\|_{L^\infty}+t\sup_{0\leq s\leq t}\|w_{\b}f^{n+1}(s)\|^2_{L^\infty}\Big\}\sup_{0\leq \tau\leq t}\|(f^{n+1}-f^{n})(\tau)\|_{L^\infty}\nonumber\\
&\leq Ct\|w_{\b}f_0\|_{L^\infty}\cdot\sup_{0\leq \tau\leq t}\|(f^{n+1}-f^{n})(\tau)\|_{L^\infty},~\mbox{for}~0\leq t\leq t_1, ~\b>14,
\end{align}
where we have used  \eqref{7.65} in the last inequality.
By the same argument as in \eqref{7.13}, one can obtain that
\begin{align}\label{7.68}
I_4\leq Ct\sup_{0\leq \tau\leq t}\|\sqrt{w_{\b}}(f^{n+1}-f^{n})(\tau)\|_{L^\infty},~\mbox{for}~\b\geq0.
\end{align}
For $I_5$, we note
\begin{align}\label{7.69}
&\Big|\Big(Q_+(\sqrt{J}f^{n+1},\sqrt{J}f^{n+1})-Q_+(\sqrt{J}f^{n},\sqrt{J}f^{n})\Big)(s,\hat{x},p)\Big|\nonumber\\
&\leq \Big|Q_+(\sqrt{J}f^{n+1},\sqrt{J}(f^{n+1}-f^n))(s,\hat{x},p)\Big|+\Big|Q_+(\sqrt{J}(f^{n+1}-f^n),\sqrt{J}f^{n})(s,\hat{x},p)\Big|.
\end{align}
By similar arguments as in \eqref{7.15}, one can obtain
\begin{align}\label{7.70}
&\f{\sqrt{w_{\b}(p)}}{\sqrt{J(p)}}\Big|Q_+(\sqrt{J}f^{n+1},\sqrt{J}(f^{n+1}-f^n))(s,\hat{x},p)\Big|\nonumber\\
&\leq C\int_{\mathbb{R}^3}\int_{\mathbb{S}^2}v_{\phi}\s(g,\t)\sqrt{J(q)}|f^{n+1}(s,\hat{x},p')|\cdot\Big|\sqrt{w_{\b}(q')}(f^{n+1}-f^n)(s,\hat{x},q')\Big|d\omega dq\nonumber\\
&~~~~~+C\f{\nu(p)}{\sqrt{w_\b(p)}}\|w_{\b}f^{n+1}(s)\|_{L^\infty}\cdot\|\sqrt{w_{\b}}(f^{n+1}-f^n)(s)\|_{L^\infty}.
\end{align}
Using similar arguments as in \eqref{7.15}-\eqref{7.62-1},  one can get
\begin{align}
& \int_{\mathbb{R}^3}\int_{\mathbb{S}^2}v_{\phi}\s(g,\t)\sqrt{J(q)}|f^{n+1}(s,\hat{x},p')|\cdot\Big|\sqrt{w_{\b}(q')}(f^{n+1}-f^n)(s,\hat{x},q')\Big|dqd\omega\nonumber\\
&\leq C\|w_{\b}f^{n+1}(s)\|_{L^\infty} \cdot\|\sqrt{w_{\b}}(f^{n+1}-f^n)(s)\|_{L^\infty},~\mbox{for}~\b>14,\nonumber
\end{align}
which, together with \eqref{7.70}, yields
\begin{align}\label{7.72}
&\f{\sqrt{w_{\b}(p)}}{\sqrt{J(p)}}\Big|Q_+(\sqrt{J}f^{n+1},\sqrt{J}(f^{n+1}-f^n))(s,\hat{x},p)\Big|\nonumber\\
&\leq C\|w_{\b}f^{n+1}(s)\|_{L^\infty}\cdot\|\sqrt{w_{\b}}(f^{n+1}-f^n)(s)\|_{L^\infty},~\mbox{for}~\b>14.
\end{align}
Similarly, it holds
\begin{align}\label{7.73}
&\f{\sqrt{w_{\b}(p)}}{\sqrt{J(p)}}\Big|Q_+(\sqrt{J}(f^{n+1}-f^n),\sqrt{J}f^{n})(s,\hat{x},p)\Big|\nonumber\\
&\leq C\|w_{\b}f^{n}(s)\|_{L^\infty}\cdot\|\sqrt{w_{\b}}(f^{n+1}-f^n)(s)\|_{L^\infty},~\mbox{for}~\b>14.
\end{align}
Then it follows from \eqref{7.69}, \eqref{7.72} and \eqref{7.73} for, $0\leq t\leq t_1$, that
\begin{align}\label{7.74}
I_5&\leq Ct\sup_{0\leq s\leq t}\Big\{\|w_{\b}f^{n+1}(s)\|_{L^\infty}+\|w_{\b}f^{n}(s)\|_{L^\infty}\Big\}\cdot\sup_{0\leq s\leq t}\|\sqrt{w_{\b}}(f^{n+1}-f^n)(s)\|_{L^\infty}\nonumber\\
&\leq Ct\|w_{\b}f_0\|_{L^\infty}\sup_{0\leq s\leq t}\|\sqrt{w_{\b}}(f^{n+1}-f^n)(s)\|_{L^\infty}~~\mbox{for}~~\b>14,
\end{align}
where we have used  \eqref{7.65} in the last inequality.

Substituting  \eqref{7.67}, \eqref{7.68} and  \eqref{7.74} into \eqref{7.66}, one obtains, for $0\leq t\leq t_1$, that
\begin{align}
&\sup_{0\leq s\leq t_1}\|\sqrt{w_{\b}}(f^{n+2}-f^{n+1})(s)\|\leq Ct_1(1+\|w_{\b}f_0\|_{L^\infty})\cdot\sup_{0\leq s\leq t_1}\|\sqrt{w_{\b}}(f^{n+1}-f^n)(s)\|_{L^\infty}\nonumber\\
&\leq \f{C}{8C_1}\sup_{0\leq s\leq t_1}\|\sqrt{w_{\b}}(f^{n+1}-f^n)(s)\|_{L^\infty}\leq \f12\sup_{0\leq s\leq t_1}\|\sqrt{w_{\b}}(f^{n+1}-f^n)(s)\|_{L^\infty},\nonumber
\end{align}
where we have chosen $C_1$ suitably large so that $\f{C}{8C_1}\leq \f12$. Thus, by the induction argument, it is straightforward to obtain that
\begin{align}
\sup_{0\leq s\leq t_1}\|\sqrt{w_{\b}}(f^{n+2}-f^{n+1})(s)\|\leq 2^{-n-1}\|\sqrt{w_{\b}}(f^1-f^0)\|_{L^\infty}
\leq 2^{-n}\|w_{\b}f_0\|_{L^\infty},\nonumber
\end{align}
which yields immediately that $f^{n+1},~n=0,1,2,\cdots$ is a Cauchy sequence.  Then there exists a limit $f$ such that
\begin{equation}
\sup_{0\leq s\leq t_1}\|\sqrt{w_{\b}}(f^n-f)(s)\|_{L^\infty}\rightarrow0 ~~\mbox{as}~~n\rightarrow +\infty,\nonumber
\end{equation}
and the limit function $f$ is indeed a mild solution to the relativistic Boltzmann equation \eqref{1.22}, \eqref{1.27-1}. Moreover, it follows from \eqref{7.65} that
\begin{equation*}
\sup_{0\leq t\leq t_1}\|w_{\b}f(t)\|_{L^\infty}\leq 2\|w_{\b}f_0\|_{L^\infty}.
\end{equation*}

Now we consider the uniqueness. Let $\tilde{f}$ be another solution of the relativistic Boltzmann equation \eqref{1.22}, \eqref{1.27-1} with $\sup_{0\leq t\leq t_1}\|w_{\b}\tilde{f}(t)\|_{L^\infty}<+\infty$, by the same arguments as in \eqref{7.66}-\eqref{7.74},  it is straightforward to obtain that
\begin{equation*}\label{7.75}
\|\sqrt{w_{\b}}(f-\tilde{f})(t)\|_{L^\infty}\leq C(1+\|w_{\b}f\|_{L^\infty}+\|w_{\b}\tilde{f}\|_{L^\infty})\cdot\int_0^t\|\sqrt{w_{\b}}(f-\tilde{f})(s)\|_{L^\infty}ds,
\end{equation*}
which, together with the Gronwall inequality, yields immediately the uniqueness, i.e., $f=\tilde{f}$.

Multiplying  \eqref{7.2} (with $F^{n+1(t)}=F^n(t)=F(t)$) by  $1,p,p_0$ and $F(t)$  integrating by parts, one can obtain \eqref{1.17}-\eqref{1.20}.

Finally, if $F_0$ (or equivalent $f_0$) is continuous, it is obvious that $F^{n+1}(t,x,p)$ (or equivalent $f^{n+1}(t,x,p)$) is continuous in $[0,\infty)\times\Omega\times\mathbb{R}^3$ since \eqref{7.2} is a linear equation. The continuity of $f(t,x,p)$  is an immediate consequence of  $\sup_{0\leq s\leq t_1}\|(f^{n+1}-f)(s)\|_{L^\infty}\rightarrow 0$ as $n\rightarrow +\infty$. Therefore the proof of Theorem \ref{thm7.1} is completed. $\hfill\Box$


\section{Global Estimates}

This section is devoted to the proof of Theorem \ref{thm1.1}. Indeed, we need only to obtain the uniform a priori  estimates to the solutions of relativistic Boltzmann equation \eqref{1.22} since we have already proved the local existence of unique solution to the relativistic Boltzmann equation for general bounded $L^\infty$ initial data in Theorem \ref{thm7.1}.

\subsection{Weighted $L^\infty$-Estimates}

Define
\begin{equation}\label{4.1-0}
h(t,x,p):= w_{\b}(p)f(t,x,p).
\end{equation}
Multiplying \eqref{1.22} by $w_\b(p)$, one gets
\begin{align}\label{4.1}
h_t+\hat{p}\cdot\nabla_xh+\nu(p)h-K_{w_{\b}}h=\Gamma_{w_{\b}}(h,h),
\end{align}
where
\begin{align}
(K_{w_{\b}}h)(p):= w_{\b}(p)\Big(K\f{h}{w_{\b}}\Big)(p) \nonumber
\end{align}
and
\begin{align}\label{4.3}
\Gamma_{w_{\b}}(h,h)&:= w_{\b}(p)\Gamma(\f{h}{w_{\b}},\f{h}{w_{\b}})= w_{\b}(p)\Gamma_+(\f{h}{w_{\b}},\f{h}{w_{\b}})- w_{\b}(p)\Gamma_-(\f{h}{w_{\b}},\f{h}{w_{\b}})\nonumber\\
&\equiv w_{\b}(p)\Gamma(f,f)= w_{\b}(p)\Gamma_+(f,f)- w_{\b}(p)\Gamma_-(f,f).
\end{align}
Then the mild form of \eqref{4.1} can be written as
\begin{align}\label{4.4}
h(t,x,p)&=e^{-\nu(p)t}h_0(x-\hat{p}t,p)
+\int_0^te^{-\nu(p)(t-s)} (K_{w_{\b}} h)(s,x-\hat{p}(t-s),p)ds\nonumber\\
&~~~~+\int_0^te^{-\nu(p)(t-s)} (\Gamma_{w_{\b}}(h,h))(s,x-\hat{p}(t-s),p)ds.
\end{align}

Firstly, we give a useful estimation on the nonlinear term $\Gamma(f,f)$. It is noted that the following key  lemma holds for the full ranges of $\g, b, a$ as in \cite{Strain,Dudy-E4}.
\begin{lemma}\label{lem4.1}
We assume that $-2<\g, a\in[0, 2+\g], b\in[0,\min\{4,4+\g\})$ for hard potentials,  and $-2<\g,b\in(0,\min\{4,4+\g\})$ for soft potentials.
Let $1< d<\min\left\{\f98, \f2{\max\{-\g,1\}},\f3{\max\{b-1,1\}}\right\}$,  then  for any $\a\geq0$, it holds that
	\begin{align}\label{4.6}
	\begin{cases}
	\Big|w_{\a}(p)\Gamma_-(f,f)(s,x,p)\Big|\leq C\nu(p)\|w_{\a}f(s)\|_{L^\infty}\cdot\|f(s)\|^{\f{4d+1}{5d}}_{L^\infty}\cdot\Big(\int_{\mathbb{R}^3}|f(s,x,q)|dq\Big)^{\f{d-1}{5d}},\\[5mm]
	\Big|w_{\a}(p)\Gamma_+(f,f)(s,x,p)\Big|\leq C\nu(p)\|w_{\a}f(s)\|_{L^\infty}\cdot\|w_{1}f(s)\|^{\f{4d+1}{5d}}_{L^\infty}\cdot\Big(\int_{\mathbb{R}^3}|f(s,x,q)|dq\Big)^{\f{d-1}{5d}},
	\end{cases}
	\end{align}
where the positive constant $C>0$ depends only on $a,b,\g, d$.
\end{lemma}

\noindent{\bf Proof.} It is noted that
\begin{align}\label{3.7}
&\Big|w_\a(p) \Gamma_-(f,f)(s,x,p)\Big|
=\Big|w_\a(p)\int_{\mathbb{R}^3\times \mathbb{S}^2} v_{\phi}\s(g,\t)\sqrt{J(q)}f(s,x,p)f(s,x,q)d\omega dq\Big|\nonumber\\
&\leq  \|w_\a f(s)\|_{L^\infty}\int_{\mathbb{R}^3\times \mathbb{S}^2} v_{\phi}\s(g,\t)\sqrt{J(q)}|f(s,x,q)|d\omega dq\nonumber\\
&\leq \|w_\a f(s)\|_{L^\infty}\left(\int_{\mathbb{R}^3\times \mathbb{S}^2} |v_{\phi}\s(g,\t)|^{d}\sqrt{J(q)}d\omega dq\right)^{\f1d}
\left(\int_{\mathbb{R}^3\times \mathbb{S}^2} \sqrt{J(q)}|f(s,x,q)|^{\f{d}{d-1}}d\omega dq\right)^{1-\f1d}\nonumber\\
&\lesssim \nu(p)\|w_\a f(s)\|_{L^\infty}
\left(\int_{\mathbb{R}^3} |f(s,x,q)|^{\f{5d}{d-1}}dq\right)^{\f{d-1}{5d}}\nonumber\\
&\lesssim \nu(p)\|w_\a f(s)\|_{L^\infty}\cdot \|f(s)\|_{L^\infty}^{\f{4d+1}{5d}}
\left(\int_{\mathbb{R}^3} |f(s,x,q)|dq\right)^{\f{d-1}{5d}},
\end{align}
where we have used Lemma \ref{lemA.7} below.\\

Next, we consider the gain term $\Gamma_+(f,f)$ which is much more complicated. Using Lemma \ref{lemA.7}, one
can get that
\begin{align}\label{4.8}
&\Big|w_{\a}(p)\Gamma_+(f,f)(s,x,p)\Big|=\Big|w_\a(p)\int_{\mathbb{R}^3\times \mathbb{S}^2} v_{\phi}\s(g,\t)\sqrt{J(q)}f(s,x,p')f(s,x,q')d\omega dq\Big|\nonumber\\
&\leq \left(\int_{\mathbb{R}^3\times \mathbb{S}^2} |v_{\phi}\s(g,\t)|^{d}\sqrt{J(q)}d\omega dq\right)^{\f1d}\cdot\Big(\int_{\mathbb{R}^3\times \mathbb{S}^2} \sqrt{J(q)}|w_\a(p) f(s,x,p')f(s,x,q')|^{\f{d}{d-1}}d\omega dq\Big)^{1-\f1d}\nonumber\\
&\lesssim \nu(p)\left\{ \int_{\mathbb{R}^3\times \mathbb{S}^2} \sqrt{J(q)}|w_\a(p)  f(s,x,p')f(s,x,q')|^{\f{d}{d-1}}d\omega dq \right\}^{1-\f1d}\nonumber\\
&\cong \nu(p)\Big\{ \int_{\mathbb{R}^3\times \mathbb{R}^3\times \mathbb{R}^3} \f{\sqrt{s}}{2g}\d^{(4)}(p^\mu+q^\mu-{p^\mu}'-{q^\mu}')\sqrt{J(q)}\nonumber\\
&~~~~~~~~~~~~~~~~~~~~~~~~~~~~~~~~~~~~~~~~~~~~~~~~~~~~\times| w_\a(p)  f(s,x,p')f(s,x,q')|^{\f{d}{d-1}}\f{dp'dq'}{p_0'q_0'} dq \Big\}^{1-\f1d},
\end{align}
where we have used the fact $d\omega=\f{\sqrt{s}}{2g}\d^{(4)}(p^\mu+q^\mu-{p^\mu}'-{q^\mu}')\f{dp'dq'}{p_0'q_0'}$ in the last equality. By the same arguments as in \eqref{7.15}-\eqref{7.18}, one obtains
\begin{align}\label{4.9}
I&:=\int_{\mathbb{R}^3\times \mathbb{R}^3\times \mathbb{R}^3} \f{\sqrt{s}}{2g}\d^{(4)}(p^\mu+q^\mu-{p^\mu}'-{q^\mu}')\sqrt{J(q)}| w_\a(p)  f(s,x,p')f(s,x,q')|^{\f{d}{d-1}}\f{dp'dq'}{p_0'q_0'} dq  \nonumber\\
&\lesssim \int_{\mathbb{R}^3\times \mathbb{R}^3\times \mathbb{R}^3} \f{\sqrt{s}}{g}\d^{(4)}(p^\mu+q^\mu-{p^\mu}'-{q^\mu}')J(q)^{\f14}| w_\a(p')  f(s,x,p')f(s,x,q')|^{\f{d}{d-1}}\f{dp'dq'dq}{p_0'q_0'q_0}  \nonumber\\
&~~+\int_{\mathbb{R}^3\times \mathbb{R}^3\times \mathbb{R}^3} \f{\sqrt{s}}{g}\d^{(4)}(p^\mu+q^\mu-{p^\mu}'-{q^\mu}')J(q)^{\f14}| f(s,x,p') w_\a(q') f(s,x,q')|^{\f{d}{d-1}}\f{dp'dq'dq}{p_0'q_0'q_0} \nonumber\\
&\lesssim \int_{\mathbb{R}^3\times \mathbb{R}^3\times \mathbb{R}^3} \f{\sqrt{\bar s}}{\bar g}\d^{(4)}(p^\mu+{p^\mu}'-q^\mu-{q^\mu}')J(p')^{\f14}| w_\a(q') f(s,x,q')f(s,x,q) |^{\f{d}{d-1}}\f{dp'dq'dq}{p_0'q_0'q_0} \nonumber\\
&\lesssim \|w_\a f(s)\|_{L^\infty}^{\f{d}{d-1}}\int_{\mathbb{R}^3\times \mathbb{R}^3\times \mathbb{R}^3} \f{\sqrt{\bar s}}{\bar g}\d^{(4)}(p^\mu+{p^\mu}'-q^\mu-{q^\mu}')J(p')^{\f14}
| f(s,x,q) |^{\f{d}{d-1}}\f{dp'dq'dq}{p_0'q_0'q_0},
\end{align}
where $\bar{g}, \bar{s}$ are defined in \eqref{7.17}.
Then, using the same procedures as in \eqref{7.22}-\eqref{7.38}, we can get that
\begin{align}\label{4.10}
\tilde{A}(p,q)&:= \f{1}{p_0q_0}\int_{\mathbb{R}^3\times \mathbb{R}^3} \f{\sqrt{\bar s}}{\bar g}\d^{(4)}(p^\mu+{p^\mu}'-q^\mu-{q^\mu}')J(p')^{\f14}
\f{dp'dq'}{p_0'q_0'}\nonumber\\
&\lesssim e^{-\f18(q_0-p_0)}\f{s}{gp_0q_0} \int_{0}^{\infty} \f{y\sqrt{1+\sqrt{y^2+1}}}{\sqrt{1+y^2}} \Big(\f{g}{\sin\f{\psi}{2}}\Big)^{-1}e^{-\f{l}{4}\sqrt{1+y^2}} I_0\Big(\f{j}{4}y\Big)dy\nonumber\\
&\lesssim e^{-\f18(q_0-p_0)}\f{s}{g^2p_0q_0} \int_{0}^{\infty} y(1+y^2)^{-\f12} e^{-\f{l}{4}\sqrt{1+y^2}} I_0\Big(\f{j}{4}y\Big)dy\nonumber\\
&\lesssim e^{-\f18(q_0-p_0)} e^{-\f18|p-q|}\f{s}{g^2p_0q_0} \f{1}{\sqrt{l^2-j^2}}
\lesssim \f{p_0q_0}{|p-q|^2(1+|p-q|)}.
\end{align}
Combining  \eqref{4.9} and \eqref{4.10}, one has that
\begin{align}\label{4.11}
I&\lesssim \|w_\a f(s)\|_{L^\infty}^{\f{d}{d-1}}
\int_{\mathbb{R}^3}\f{p_0^2q_0}{|p-q|^2(1+|p-q|)}| f(s,x,q) |^{\f{d}{d-1}}dq\nonumber\\
&\lesssim \|w_\a f(s)\|_{L^\infty}^{\f{d}{d-1}}
\left(\int_{\mathbb{R}^3}\f{p_0^{\f52}q_0^{\f54}(1+|q|)^{-10}}{|p-q|^{\f52}(1+|p-q|)^{\f54}}dq\right)^{\f45}
\left(\int_{\mathbb{R}^3}(1+|q|)^{40}| f(s,x,q) |^{\f{5d}{d-1}}dq\right)^{\f15}\nonumber\\
&\lesssim p_0^{-1}\|w_\a f(s)\|_{L^\infty}^{\f{d}{d-1}}\cdot \Big\|q_0^{\f{40(d-1)}{4d+1}}f(s)\Big\|_{L^\infty}^{\f{4d+1}{5(d-1)}}
\left(\int_{\mathbb{R}^3}| f(s,x,q) |dq\right)^{\f15}\nonumber\\
&\lesssim p_0^{-1}\|w_\a f(s)\|_{L^\infty}^{\f{d}{d-1}}\cdot \Big\|w_1f(s)\Big\|_{L^\infty}^{\f{4d+1}{5(d-1)}}
\left(\int_{\mathbb{R}^3}| f(s,x,q) |dq\right)^{\f15},
\end{align}
where we have used the fact $\f{40(d-1)}{4d+1}\leq1$ due to $d<\f98$. Then, substituting  \eqref{4.11} into \eqref{4.8}, one has
\begin{align}\label{4.12}
\Big|w_{\a}(p)\Gamma_+(f,f)(s,x,p)\Big|\leq C\nu(p) \|w_\a f(s)\|_{L^\infty}\cdot \|w_1f(s)\|_{L^\infty}^{\f{4d+1}{5d}}
\left(\int_{\mathbb{R}^3}| f(s,x,q) |dq\right)^{\f{d-1}{5d}}.
\end{align}
Therefore the proof of Lemma \ref{lem4.1} is completed.  $\hfill\Box$

\

\begin{lemma}\label{lem4.2}
Under the hypothesis {\bf H)},  it holds, for $\b>14$, that
	\begin{align}\label{4.13}
	\sup_{0\leq s\leq t}\|h(s)\|_{L^\infty}
	&\leq C_2\Big\{\|h_0\|_{L^\infty}+\|h_0\|^2_{L^\infty}+\sqrt{\mathcal{E}(F_0)}+\mathcal{E}(F_0)\Big\}\nonumber\\
	&~~~+C_2\sup_{t_1\leq s\leq t, z\in\Omega}\left\{ \|h(s)\|^{\f{9d+1}{5d}}_{L^\infty}\Big(\int_{\mathbb{R}^3}|f(s,z,\eta)|d\eta\Big)^{\f{d-1}{5d}}\right\},
	\end{align}
	where $t_1>0$  is defined in \eqref{LT}, and the  constant $C_2\geq1$ depends only on $a,b,\g,\b$.
\end{lemma}

\noindent{\bf Proof.} It follows from \eqref{4.4} that
\begin{align}\label{4.14}
|h(t,x,p)|&\leq e^{-\nu(p)t}|h_0(x-\hat{p}t,p)|+\int_0^te^{-\nu(p)(t-s)} \Big|(K_{w_{\b}} h)(s,x-\hat{p}(t-s),p)\Big|ds\nonumber\\
&~~~~+\int_0^te^{-\nu(p)(t-s)} \Big|(\Gamma_{w_{\b}}(h,h))(s,x-\hat{p}(t-s),p)\Big|ds\nonumber\\
&:=e^{-\nu(p)t}\|h_0\|_{L^\infty}+D_1+D_2.
\end{align}
For $D_2$,   it follows from \eqref{4.3} and \eqref{4.6} that
\begin{align}\label{4.16}
D_2&\leq C\sup_{0\leq s\leq t,~z\in\Omega}\left\{ \|h(s)\|^{\f{9d+1}{5d}}_{L^\infty}\Big(\int_{\mathbb{R}^3}|f(s,z,\eta)|d\eta\Big)^{\f{d-1}{5d}}\right\},~~\mbox{for}~\b\geq1.
\end{align}

It remains to consider $D_1$ whose estimation is rather complicated. Let $k_{w_\b}(p,q)$ be the corresponding kernel associated with $K_{w_\b}$, then it holds that
\begin{align}\label{4.20}
k_{w_\b}(p,q)=k(p,q)\f{w_\b(p)}{w_\b(q)},
\end{align}
which, together with \eqref{2.47} and \eqref{2.48}, yields  that
\begin{equation}\label{4.20-1}
\int_{\mathbb{R}^3}\Big|k_{w_\b}(p,q)\Big| dq\leq C
~~\mbox{and}~~\int_{\mathbb{R}^3}\Big|k_{w_\b}(p,q)\Big| dq\leq C \f{\nu(p)}{p_0^{\xi}},
\end{equation}
where $\xi:=\min\{\xi_1,\xi_2\}>0$, and $\xi_1, \xi_2$ are defined in Lemma \ref{lem2.4}.
Denoting $\hat{x}:= x-\hat{p}(t-s)$, similar to  \cite{Vidav,Strain}, we use \eqref{4.4} again to get
\begin{align}\label{4.23}
D_1&\leq \|h_0\|_{L^\infty}\int_0^te^{-\nu(p)(t-s)} \int_{\mathbb{R}^3}e^{-\nu(q)s}|k_{w_\b}(p,q)|dq ds\nonumber\\
&~~+\int_0^te^{-\nu(p)(t-s)} \int_{\mathbb{R}^3}|k_{w_\b}(p,q)| \int_0^se^{-\nu(q)(s-\tau)} |(\Gamma_{w_\b}(h,h))(\tau,\hat{x}-\hat{q}(s-\tau),q)|d\tau dqds\nonumber\\
&~~+\int_0^te^{-\nu(p)(t-s)} \int_{\mathbb{R}^3}\int_{\mathbb{R}^3}|k_{w_\b}(p,q)k_{w_\b}(q,\eta)|\nonumber\\
&~~~~~~~~~~~~~~~~~~~~~~~~~~~~~~~~~~~~~~\times \int_0^se^{-\nu(q)(s-\tau)} |h(\tau,\hat{x}-\hat{q}(s-\tau),\eta)|d\tau dqd\eta ds \nonumber\\
&:=D_{11}+D_{12}+D_{13}.
\end{align}
It follows from  \eqref{4.16}  and \eqref{4.20-1}, for $\b\geq1$, that
\begin{align}\label{4.24}
D_{11}+D_{12}
&\leq C\left\{\|h_0\|_{L^\infty}+\sup_{0\leq s\leq t,z\in\Omega}\Big[ \|h(s)\|^{\f{9d+1}{5d}}_{L^\infty}\Big(\int_{\mathbb{R}^3}|f(s,z,\eta)|d\eta\Big)^{\f{d-1}{5d}}\Big]\right\} \nonumber\\
&~~~~~~~~~~~~~~~~~~~~~~~~~~~~~~~~~~~~~~~~~~~~~~\times \int_0^te^{-\nu(p)(t-s)} \int_{\mathbb{R}^3}|k_{w_\b}(p,q)| dqds\nonumber\\
&\leq C\left\{\|h_0\|_{L^\infty}+\sup_{0\leq s\leq t,z\in\Omega}\Big[ \|h(s)\|^{\f{9d+1}{5d}}_{L^\infty}\Big(\int_{\mathbb{R}^3}|f(s,z,\eta)|d\eta\Big)^{\f{d-1}{5d}}\Big]\right\}.
\end{align}

\

We now concentrate on the last term $D_{13}$ on the RHS of \eqref{4.23}.  Motivated by  \cite{DHWY,Guo}, we divide the proof  into the following three cases.\\[2mm]
\noindent{\it Case 1}. For $|p|\geq N$, it follows from \eqref{4.20-1} that
\begin{equation*}\label{4.27}
\int_{\mathbb{R}^3}\Big|k_{w_\b}(p,q)\Big| dq\leq C\f{\nu(p)}{N^\xi},
\end{equation*}
which yields immediately that
\begin{align}\label{4.28}
D_{13}&\leq C\sup_{0\leq s\leq t}\|h(s)\|_{L^\infty}\int_0^te^{-\nu(p)(t-s)} \int_{\mathbb{R}^3}|k_{w_\b}(p,q)|\int_0^se^{-\nu(q)(s-\tau)} \f{\nu(q)}{(1+|q|)^\xi}d\tau dqds\nonumber\\
&\leq \f{C}{N^\xi}\sup_{0\leq s\leq t}\|h(s)\|_{L^\infty}.
\end{align}

\noindent{\it Case 2}. For  $|p|\leq N,~|q|\geq2N$ (or $|p|\leq N,|q|\leq2N,~|\eta|\geq 3N$),  we have  $|p-q|\geq N$ (or $|q-\eta|\geq N$), and
\begin{align}
|k_{w_\b}(p,q)|\leq Ce^{-\f{N}{20}}\Big|k_{w_\b}(p,q)e^{\f{|p-q|}{20}}\Big|~
\Big(\mbox{or}~ |k_{w_\b}(q,\eta)|\leq Ce^{-\f{N}{20}}\Big|k_{w_\b}(q,\eta)e^{\f{|q-\eta|}{20}}\Big|\Big).\label{4.29}
\end{align}
By similar arguments as in Lemma \ref{lem2.5}, it holds that
\begin{align}
\small{\int_{\mathbb{R}^3}\Big|k_{w_\b}(p,q)e^{\f{|p-q|}{20}}\Big|dq\leq \f{C\nu(p)}{(1+|p|)^\xi}~~\mbox{and}~\int_{\mathbb{R}^3}\Big|k_{w_\b}(q,\eta)e^{\f{|q-\eta|}{20}}\Big|d\eta\leq \f{C\nu(q)}{(1+|q|)^\xi}.}\label{4.31}
\end{align}
Then it follows from \eqref{4.29} and \eqref{4.31} that
\begin{align}\label{4.32}
&\int_0^te^{-\nu(p)(t-s)} \Big\{\int_{|p|\leq N,|q|\geq2N}+\int_{|q|\leq 2N,|\eta|\geq3N}\Big\}|k_{w_\b}(p,q)k_{w_\b}(q,\eta)|\nonumber\\
&~~~~~~~~~~~~~~~~~~~~~~~~~~~~~~~~~~~\times \int_0^se^{-\nu(q)(s-\tau)} |h(\tau,\hat{x}-\hat{q}(s-\tau),\eta)|d\tau d\eta dqds\nonumber\\
&\leq Ce^{-\f{N}{20}}\sup_{0\leq s\leq t}\|h(s)\|_{L^\infty}.
\end{align}

\noindent{\it Case 3}. For $|p|\leq N,~|q|\leq 2N,~|\eta|\leq 3N$,  we get
\begin{align}\label{4.33}
&\int_0^te^{-\nu(p)(t-s)} \int_{|q|\leq 2N,|\eta|\leq3N}|k_{w_\b}(p,q)k_{w_\b}(q,\eta)|\int_0^se^{-\nu(q)(s-\tau)} |h(\tau,\hat{x}-\hat{q}(s-\tau),\eta)|d\tau d\eta dqds\nonumber\\
&\leq \int_0^te^{-\nu(p)(t-s)} \int_{|q|\leq 2N,|\eta|\leq3N}|k_{w_\b}(p,q)k_{w_\b}(q,\eta)|\int_{s-\k}^se^{-\nu(q)(s-\tau)} |h(\tau,
\hat {x}-\hat{q}(s-\tau),\eta)|d\tau d\eta dqds\nonumber\\
&~~~+\int_0^te^{-\nu(p)(t-s)} \int_{|q|\leq 2N,|\eta|\leq3N}|k_{w_\b}(p,q)k_{w_\b}(q,\eta)|\nonumber\\
&~~~~~~~~~~~~~~~~~~~~~~~~~~~~~~~~~~~~~~~~~~~~~~~\times\int_0^{s-\k}e^{-\nu(q)(s-\tau)} |h(\tau,\hat{x}-\hat{q}(s-\tau),\eta)|d\tau d\eta dqds.
\end{align}
Using \eqref{4.20-1}, we can control the first term on the RHS of \eqref{4.33} by
\begin{align}\label{4.34}
C\k\sup_{0\leq s\leq t}\|h(s)\|_{L^\infty}\int_0^te^{-\nu(p)(t-s)} \nu(p)ds\leq C\k\sup_{0\leq s\leq t}\|h(s)\|_{L^\infty}.
\end{align}
Now we estimate the second term on the RHS of \eqref{4.33}.  Since $k_{w_\b}(p,q)$ may have  singularity of $\f{1}{|p-q|}$, we  choose a smooth compact support function $k_N(p,q)$ such that
\begin{align}\label{4.35}
\sup_{|p|\leq 3N}\int_{|q|\leq 3N}\Big|k_{w_\b}(p,q)-k_N(p,q)\Big|dq\leq \f{C}{N^7}.
\end{align}
Noting
\begin{align}\label{4.36}
k_{w_\b}(p,q)k_{w_\b}(q,\eta)&=\Big(k_{w_\b}(p,q)-k_N(p,q)\Big)k_{w_\b}(q,\eta)\nonumber\\
&~~~+\Big(k_{w_\b}(q,\eta)-k_N(q,\eta)\Big)k_N(p,q)+k_N(p,q)k_N(q,\eta),
\end{align}
and using \eqref{4.35} and \eqref{4.36},  the second term on the RHS of \eqref{4.33} is bounded by
\begin{align}\label{4.37}
&\f{C}{N^7}\sup_{0\leq s\leq t}\|h(s)\|_{L^\infty}\int_0^te^{-\nu(p)(t-s)}\int_0^{s-\k}e^{-\nu(q)(s-\tau)} d\tau ds \nonumber\\
&~~~~~~~~~~~~~~~~\times \left\{\sup_{|q|\leq2N}\int_{|\eta|\leq3N}|k_{w_\b}(q,\eta)|d\eta+\sup_{|p|\leq N}\int_{|q|\leq2N}|k_N(p,q)|dq\right\}\nonumber\\
&+\int_0^te^{-\nu(p)(t-s)} \int_{|q|\leq 2N,|\eta|\leq3N}|k_N(p,q)k_N(q,\eta)|\nonumber\\
&~~~~~~~~~~~~~~~~~~~~~\times \int_0^{s-\k}e^{-\nu(q)(s-\tau)} |h(\tau,\hat{x}-\hat{q}(s-\tau),\eta)|d\tau dqd\eta ds\nonumber\\
&\leq \f{C}{N}\sup_{0\leq s\leq t}\|h(s)\|_{L^\infty}+C_{N}\int_0^t\int_0^{s-\k}e^{-c_N(t-s)}e^{-c_N(s-\tau)}  \nonumber\\
&~~~~~~~~~~~~~~~~~~~~~~~~~~~~~~~~~~~~~\times\int_{|q|\leq 2N,|\eta|\leq3N} |h(\tau,\hat{x}-\hat{q}(s-\tau),\eta)|dqd\eta d\tau ds,
\end{align}
where we have used the facts that   $k_N(p,q)k_N(q,\eta)$ is bounded and
\begin{equation*}\label{4.38}
\nu(p)\geq c_N~~\mbox{for}~|p|\leq N,~~\mbox{and}~~\nu(q)\geq c_N~~\mbox{for}~|q|\leq 2N.
\end{equation*}
It follows from \eqref{2.52}  that
\begin{align}\label{4.39}
	&\int_{|q|\leq 2N,|\eta|\leq3N} |h(\tau,\hat{x}-\hat{q}(s-\tau),\eta)|d\eta dq\nonumber\\
	&\leq C_{N}\int_{|q|\leq 2N,|\eta|\leq3N} \f{|F(\tau,\hat{x}-\hat{q}(s-\tau),\eta)-J(\eta)|}{\sqrt{J(\eta)}} I_{\{|F(\tau,\hat{x}-\hat{q}(s-\tau),\eta)-J(\eta)|\leq J(\eta)\}}d\eta dq\nonumber\\
	&~~~+C_{N}\int_{|q|\leq 2N,|\eta|\leq3N} |F(\tau,\hat{x}-\hat{q}(s-\tau),\eta)-J(\eta)| I_{\{|F(\tau,\hat{x}-\hat{q}(s-\tau),\eta)-J(\eta)|\geq J(\eta)\}}d\eta dq \nonumber\\
	&\leq C_{N}\f{1+(s-\tau)^{\f32}}{(s-\tau)^{\f32}}\left\{\int_{\Omega}\int_{|\eta|\leq3N} \f{|F(\tau,z,\eta)-J(\eta)|^2}{J(\eta)} I_{\{|F(\tau,z,\eta)-J(\eta)|\leq J(\eta)\}}d\eta dz\right\}^{\f12}\nonumber\\
	 &~~~~~+C_{N}\f{1+(s-\tau)^{3}}{(s-\tau)^{3}}\int_{\Omega}\int_{|\eta|\leq3N} |F(\tau,z,\eta)-J(\eta)| I_{\{|F(\tau,z,\eta)-J(\eta)|\geq J(\eta)\}}d\eta dz\nonumber\\
	&\leq C_{N}\f{1+(s-\tau)^{\f32}}{(s-\tau)^{\f32}}\sqrt{\mathcal{E}(F_0)}+C_{N}\f{1+(s-\tau)^{3}}{(s-\tau)^{3}}\mathcal{E}(F_0),
	\end{align}
where we have made a change of variable $z:=\hat{x}-\hat{q}(s-\tau)$ and used
\begin{equation}\label{4.39-1}
dz=(s-\tau)^{3}d\hat{q}=\f{(s-\tau)^{3}}{(1+|q|^2)^{\f52}}dq.
\end{equation}
Using \eqref{4.39}, we can bound the second term on the RHS of \eqref{4.37} as follows
\begin{align}\label{4.40}
&C_{N}\int_0^t\int_0^{s-\k}e^{-c_N(t-s)}e^{-c_N(s-\tau)}\int_{|q|\leq 2N,|\eta|\leq3N} |h(\tau,\hat{x}-\hat{q}(s-\tau),\eta)|d\eta d\tau dqds\nonumber\\
&\leq C_{N}\k^{-\f32}\sqrt{\mathcal{E}(F_0)}+C_{N}\k^{-3}\mathcal{E}(F_0).
\end{align}
Combining \eqref{4.40}, \eqref{4.37}, \eqref{4.34} and \eqref{4.33} gives
\begin{align}\label{4.41}
&\int_0^te^{-\nu(p)(t-s)} \int_{|q|\leq 2N,|\eta|\leq3N}|k_{w_\b}(p,q)k_{w_\b}(q,\eta)|\int_0^se^{-\nu(q)(s-\tau)} |h(\tau,\hat{x}-\hat{q}(s-\tau),\eta)|d\tau d\eta dqds\nonumber\\
&\leq C(\k+\f{1}{N})\sup_{0\leq s\leq t}\|h(s)\|_{L^\infty}+C_{N}\k^{-\f32}\sqrt{\mathcal{E}(F_0)}+C_{N}\k^{-3}\mathcal{E}(F_0).
\end{align}
Thus collecting  all the above estimates,  we get, for any $\k>0$ and large $N\geq1$, that
\begin{align}\label{4.42}
\sup_{0\leq s\leq t}\|h(s)\|_{L^\infty}&\leq C_{N}\Big(\k^{-\f32}\sqrt{\mathcal{E}(F_0)}+\k^{-3}\mathcal{E}(F_0)+ \|h_0\|_{L^\infty}\Big)+C\Big(\k+\f{1}{N^{\xi}}\Big)\cdot\sup_{0\leq s\leq t}\|h(s)\|_{L^\infty}\nonumber\\
&~~~~+C\sup_{0\leq s\leq t,z\in\Omega}\left\{ \|h(s)\|^{\f{9d+1}{5d}}_{L^\infty}\Big(\int_{\mathbb{R}^3}|f(s,z,\eta)|d\eta\Big)^{\f{d-1}{5d}}\right\}.
\end{align}
First choosing  $\k$ small, then letting $N$ sufficiently large so that $C\Big(\k+\f{1}{N^\xi}\Big)\leq \f12$, one obtains, for $\b\geq1$, that
\begin{align}\label{4.43}
\sup_{0\leq s\leq t}\|h(s)\|_{L^\infty}&\leq C\Big\{\|h_0\|_{L^\infty}+\sqrt{\mathcal{E}(F_0)}+\mathcal{E}(F_0)\Big\}\nonumber\\
&~~~~~~~~~+C\sup_{0\leq s\leq t,z\in\Omega}\left\{ \|h(s)\|^{\f{9d+1}{5d}}_{L^\infty}\Big(\int_{\mathbb{R}^3}|f(s,z,\eta)|d\eta\Big)^{\f{d-1}{5d}}\right\}.
\end{align}
Using Theorem \ref{thm7.1}, one has, for $\b>14$, that
\begin{equation}\label{4.43-1}
\sup_{0\leq s\leq t_1,z\in\Omega}\left\{ \|h(s)\|^{\f{9d+1}{5d}}_{L^\infty}\Big(\int_{\mathbb{R}^3}|f(s,z,\eta)|d\eta\Big)^{\f{d-1}{5d}}\right\}\leq C\sup_{0\leq s\leq t_1}\|h(s)\|^2_{L^\infty}\leq C\|h_0\|^2_{L^\infty}.
\end{equation}
Substituting \eqref{4.43-1} into \eqref{4.43}  implies, for $\b>14$, that
\begin{align}\label{4.45}
\sup_{0\leq s\leq t}\|h(s)\|_{L^\infty}&\leq C_2\Big\{\|h_0\|_{L^\infty}+\|h_0\|^2_{L^\infty}+\sqrt{\mathcal{E}(F_0)}+\mathcal{E}(F_0)\Big\}\nonumber\\
&~~~~~~~+C_2\sup_{t_1\leq s\leq t,z\in\Omega}\left\{ \|h(s)\|^{\f{9d+1}{5d}}_{L^\infty}\Big(\int_{\mathbb{R}^3}|f(s,z,\eta)|d\eta\Big)^{\f{d-1}{5d}}\right\},
\end{align}
which yields immediately \eqref{4.13}, where the positive constant $C_2>0$ depends only on $a,b,\g,\b,d$.  Therefore the proof of Lemma \ref{lem4.2} is completed. $\hfill\Box$

\subsection{$L^\infty_xL^1_p$-Estimates}

In this subsection, we will concentrate on the estimate of
$$\int_{\mathbb{R}^3}|f(t,x,p)|dp.$$
We observe that if $\mathcal{E}(F_0)+\|f_0\|_{L^1_xL^\infty_p}$ is small,  $\int_{\mathbb{R}^3}|f(t,x,p)|dp$ is also small for $t\geq t_1$ even though it may be not small initially.  Indeed, we have the following lemma which plays a key role in this paper.
\begin{lemma}\label{lem5.1}
	Under the hypothesis {\bf H)},   it holds, for $\b>5$, that
	\begin{align}\label{5.2}
	\int_{\mathbb{R}^3}|f(t,x,p)|dp&\leq \int_{\mathbb{R}^3}e^{-\nu(p)t}|f_0(x-\hat{p}t,p)|dp+C_{\k,N}\Big[\sqrt{\mathcal{E}(F_0)}+\mathcal{E}(F_0)\Big]\nonumber\\
	&~~+C\Big(\k+\f{1}{N}\Big)\cdot\sup_{0\leq s\leq t}\Big\{\|h(s)\|_{L^\infty}+\|h(s)\|^2_{L^\infty}\Big\}\nonumber\\
	 &~~+C_{\k,N}\Big(\sqrt{\mathcal{E}(F_0)}+\mathcal{E}(F_0)\Big)^{1-\f1d}\cdot\sup_{0\leq s\leq t}\|h(s)\|^{1+\f1d}_{L^\infty}\nonumber\\
	 &~~+C_{\k,N}\Big(\sqrt{\mathcal{E}(F_0)}+\mathcal{E}(F_0)\Big)^{\f15(1-\f1d)}\sup_{0\leq s\leq t}\|h(s)\|^{\f{9d+1}{5d}}_{L^\infty},
	\end{align}
	where $\k>0$ and  $N\geq1$ will be chosen later.
\end{lemma}

\noindent{\bf Proof.} It follows from \eqref{1.27}  that
\begin{align}\label{5.3}
&\int_{\mathbb{R}^3}|f(t,x,p)|dp\nonumber\\
&\leq \int_{\mathbb{R}^3}e^{-\nu(p)t}|f_0(x-\hat{p}t,p)|dp+\int_0^t\int_{\mathbb{R}^3}e^{-\nu(p)(t-s)} \Big|(Kf)(s,x-\hat{p}(t-s),p)\Big|dpds\nonumber\\
&~~~~+\int_0^t\int_{\mathbb{R}^3}e^{-\nu(p)(t-s)} \Big|\Gamma(f,f)(s,x-\hat{p}(t-s),p)\Big|dpds\nonumber\\
&:= \int_{\mathbb{R}^3}e^{-\nu(p)t}|f_0(x-\hat{p}t,p)|dp+H_1+H_2.
\end{align}
For $H_1$, we    notice that
\begin{align}\label{5.5}
H_1
&= \int_{t-\k}^t\int_{\mathbb{R}^3}e^{-\nu(p)(t-s)} \Big|\int_{\mathbb{R}^3}k(p,q)f(s,x-\hat{p}(t-s),q)dq\Big|dpds\nonumber\\
&~~~~+\int_0^{t-\k}\int_{\mathbb{R}^3}e^{-\nu(p)(t-s)} \Big|\int_{\mathbb{R}^3}k(p,q)f(s,x-\hat{p}(t-s),q)dq\Big|dpds\nonumber\\
&:= H_{11}+H_{12}.
\end{align}
It is straightforward to obtain, for $\b> 3$, that
\begin{align}\label{5.6}
H_{11}&\leq \sup_{0\leq s\leq t}\|h(s)\|_{L^\infty}\int_{t-\k}^t\int_{\mathbb{R}^3}e^{-\nu(p)(t-s)}(1+|p|^2)^{-\f\b2} \int_{\mathbb{R}^3}\Big|k_{w_\b}(p,q)\Big| dqdp ds\nonumber\\
&\leq C\k \sup_{0\leq s\leq t}\|h(s)\|_{L^\infty},
\end{align}
where we have used the first part of \eqref{4.20-1} in the last inequality. For  $H_{12}$, one has that
\begin{align}\label{5.7}
H_{12}&\leq\int_0^{t-\k}\int_{|p|\geq N}e^{-\nu(p)(t-s)} \Big|\int_{\mathbb{R}^3_{q}}k(p,q)f(s,x-\hat{p}(t-s),q)dq\Big|dpds\nonumber\\
&~~~+\int_0^{t-\k}\int_{|p|\leq N}e^{-\nu(p)(t-s)} \Big|\int_{|q|\geq 2N}k(p,q)f(s,x-\hat{p}(t-s),q)dq\Big|dpds\nonumber\\
&~~~+\int_0^{t-\k}\int_{|p|\leq N}e^{-\nu(p)(t-s)} \Big|\int_{|q|\leq 2N}k(p,q)f(s,x-\hat{p}(t-s),q)dq\Big|dpds\nonumber\\
&:= H_{121}+H_{122}+H_{123}.
\end{align}
It follows from \eqref{4.20-1} and \eqref{4.31}, for $\b>3$,  that
\begin{align}\label{5.8}
H_{121}&\leq \sup_{0\leq s\leq t}\|h(s)\|_{L^\infty}\cdot\int_0^{t-\k}\int_{|p|\geq N}e^{-\nu(p)(t-s)}w_{\b}(p)^{-1} \Big|\int_{\mathbb{R}^3_{q}}k_{w_\b}(p,q)dq\Big|dpds\nonumber\\
&\leq C\sup_{0\leq s\leq t}\|h(s)\|_{L^\infty}\int_0^{t-\k}\int_{|p|\geq N}e^{-\nu(p)(t-s)}(1+|p|)^{-\xi-\b}\nu(p) dpds\nonumber\\
&\leq \f{C}{N^{\b-3}}\sup_{0\leq s\leq t}\|h(s)\|_{L^\infty},
\end{align}
and
\begin{align}\label{5.9}
H_{122}&\leq e^{-\f{N}{20}}\sup_{0\leq s\leq t}\|h(s)\|_{L^\infty}\cdot\int_0^{t-\k}\int_{|p|\leq N}e^{-\nu(p)(t-s)}w_{\b}(p)^{-1}  \Big|\int_{|q|\geq2N}k_{w_\b}(p,q)e^{\f{|p-q|}{20}}dq\Big|dpds\nonumber\\
&\leq C e^{-\f{N}{20}}\sup_{0\leq s\leq t}\|h(s)\|_{L^\infty}.
\end{align}
Noting $k_{w_{\b}}(p,q)$ may have singularity of $\f1{|p-q|}$,  then using   \eqref{4.35}, \eqref{4.36} and  \eqref{4.39}, we obtain
\begin{align}\label{5.10}
H_{123}&\leq \int_0^{t-\k}\int_{|p|\leq N}e^{-\nu(p)(t-s)}w_{\b}(p)^{-1} \int_{|q|\leq 2N}|k_{w_{\b}}(p,q)-k_N(p,q)|\cdot|h(s,x-\hat{p}(t-s),q)|dqdpds\nonumber\\
&~~~+\int_0^{t-\k}\int_{|p|\leq N}e^{-\nu(p)(t-s)} \int_{|q|\leq 2N}\Big|k_N(p,q)h(s,x-\hat{p}(t-s),q)\Big|dqdpds\nonumber\\
&\leq \f{C}{N}\sup_{0\leq s\leq t}\|h(s)\|_{L^\infty}+C_N\int_0^{t-\k}e^{-c_N(t-s)}\int_{|p|\leq N,|q|\leq 2N} \Big|h(s,x-\hat{p}(t-s),q)\Big|dqdpds\nonumber\\
&\leq \f{C}{N}\sup_{0\leq s\leq t}\|h(s)\|_{L^\infty}+C_N\Big[\k^{-\f32}\sqrt{\mathcal{E}(F_0)}+\k^{-3}\mathcal{E}(F_0)\Big].
\end{align}
Hence combining \eqref{5.5}-\eqref{5.10}, one gets, for $\b>3$, that
\begin{align}\label{5.11}
H_1\leq C\left(\k+\f1N+\f{1}{N^{\b-3}}\right)\sup_{0\leq s\leq t}\|h(s)\|_{L^\infty}+C_N\Big[\k^{-\f32}\sqrt{\mathcal{E}(F_0)}+\k^{-3}\mathcal{E}(F_0)\Big].
\end{align}

\

Next we estimate $H_2$. Firstly,  we note that
\begin{align}\label{5.12}
H_2&\leq \int_0^t\int_{\mathbb{R}^3}e^{-\nu(p)(t-s)} \int_{\mathbb{R}^3}\int_{\mathbb{S}^2}v_{\phi}\s(g,\t)\sqrt{J(q)}\nonumber\\
&~~~~~~~~~~~~~~~~~~~~~~~~~\times\Big|f(s,x-\hat{p}(t-s),p)f(s,x-\hat{p}(t-s),q)\Big|d\omega dqdpds\nonumber\\
&~~~+\int_0^t\int_{\mathbb{R}^3}e^{-\nu(p)(t-s)} \int_{\mathbb{R}^3}\int_{\mathbb{S}^2}v_{\phi}\s(g,\t)\sqrt{J(q)}\nonumber\\
&~~~~~~~~~~~~~~~~~~~~~~~~~\times\Big|f(s,x-\hat{p}(t-s),p')f(s,x-\hat{p}(t-s),q')\Big|d\omega  dqdpds\nonumber\\
&:= H_{21}+H_{22}.
\end{align}
For $H_{21}$, one has, for $\b>5$, that
\begin{align}\label{5.13}
H_{21}&\leq C\int_{t-\k}^t\int_{\mathbb{R}^3}e^{-\nu(p)(t-s)}\nu(p)w_{\b}(p)^{-1}\|h(s)\|^2_{L^\infty} dpds\nonumber\\
&~~+C\int_0^{t-\k}\int_{\mathbb{R}^3}\int_{\mathbb{R}^3}\int_{\mathbb{S}^2}\|h(s)\|_{L^\infty}e^{-\nu(p)(t-s)} w_{\b}(p)^{-1}\nonumber\\
&~~~~~~~~~~~~~~~~~~~~~~~~~~~~~~~~~~~~~~\times v_{\phi}\s(g,\t)\sqrt{J(q)}|f(s,x-\hat{p}(t-s),q)|d\omega dqdpds\nonumber\\
&\leq C\k\sup_{0\leq s\leq t}\|h(s)\|^2_{L^\infty}+C\int_0^{t-\k}\Big\{\int_{|p|\geq N}\int_{\mathbb{R}^3_{q}}\int_{\mathbb{S}^2}+\int_{\mathbb{R}^3_p}\int_{|q|\geq N}\int_{\mathbb{S}^2}\Big\}\{\cdots\}dqdpds\nonumber\\
&~~~~~~~~~~~~~~+C\int_0^{t-\k}\int_{|p|\leq N}\int_{|q|\leq N}\int_{\mathbb{S}^2}\{\cdots\}d\omega dqdpds\nonumber\\
&\leq C\left(\k+\f{1}{N^{\b-3}}\right)\sup_{0\leq s\leq t}\|h(s)\|^2_{L^\infty}+C\int_0^{t-\k}\int_{|p|\leq N}\int_{|q|\leq N}\int_{\mathbb{S}^2}\{\cdots\}d\omega dqdpds.
\end{align}
For  the last term of above, it holds that
\begin{align}\label{5.14}
&C\int_0^{t-\k}\int_{|p|\leq N}\int_{|q|\leq N}\int_{\mathbb{S}^2}\{\cdots\}d\omega dqdpds\nonumber\\
&\leq C\int_0^{t-\k}\int_{|p|\leq N}\|h(s)\|_{L^\infty}e^{-c_N(t-s)}w_{\b}(p)^{-1}\nonumber\\
&~~~~~~~~~~~~~~~~~~~~~~~~\times\int_{|q|\leq N}\int_{\mathbb{S}^2}  v_{\phi}\s(g,\t)\sqrt{J(q)}|f(s,x-\hat{p}(t-s),q)|d\omega dqdpds\nonumber\\
&\leq C\int_0^{t-\k}e^{-c_N(t-s)}\|h(s)\|_{L^\infty}\left(\int_{|p|\leq N}\int_{|q|\leq N}\int_{\mathbb{S}^2}|v_{\phi}\s(g,\t)|^{d}\sqrt{J(q)}d\omega dq dp\right)^{\f1d}\nonumber\\
&~~~~~~~~~~~~~~~~~~~~~~~~~~~~~~~~~~~\times\left(\int_{|p|\leq N}\int_{|q|\leq N}|f(s,x-\hat{p}(t-s),q)|^{\f{d}{d-1}}dqdp\right)^{1-\f1d}ds\nonumber\\
&\leq C_N\int_0^{t-\k}e^{-c_N(t-s)}\|h(s)\|^{1+\f1d}_{L^\infty}\left(\int_{|p|\leq N}\int_{|q|\leq N}|f(s,x-\hat{p}(t-s),q)|dqdp\right)^{1-\f1d}ds\nonumber\\
&\leq C_N\k^{-3}\Big[\sqrt{\mathcal{E}(F_0)}+\mathcal{E}(F_0)\Big]^{1-\f1d}\sup_{0\leq s\leq t}\|h(s)\|^{1+\f1d}_{L^\infty}.
\end{align}
where we have used
\begin{equation*}\label{5.15}
\int_{|p|\leq N,|q|\leq 3N} |f(s,x-\hat{p}(t-s),q)|dqdp\leq C_N\f{1+(t-s)^{\f32}}{(t-s)^{\f32}}\sqrt{\mathcal{E}(F_0)}+C_N\f{1+(t-s)^{3}}{(t-s)^{3}}\mathcal{E}(F_0).
\end{equation*}
Hence, from \eqref{5.13} and \eqref{5.14}, one obtains, for $\b>5$,  that
\begin{equation}\label{5.16}
H_{21}\leq C\left(\k+\f1{N^{\b-3}}\right)\sup_{0\leq s\leq t}\|h(s)\|^2_{L^\infty}+C_N\k^{-3}\Big(\sqrt{\mathcal{E}(F_0)}+\mathcal{E}(F_0)\Big)^{1-\f1d}\sup_{0\leq s\leq t}\|h(s)\|^{1+\f1d}_{L^\infty}.
\end{equation}
Now we estimate $H_{22}$. For  $\b>5$,   it holds that
\begin{align}\label{5.17}
H_{22}&\leq C\int_{t-\k}^t\int_{\mathbb{R}^3}e^{-\nu(p)(t-s)}\nu(p)w_{\b}(p)^{-1}\|h(s)\|^2_{L^\infty} dpds\nonumber\\
&+\int_0^{t-\k}\int_{\mathbb{R}^3}\int_{\mathbb{R}^3}\int_{\mathbb{S}^2}\|h(s)\|_{L^\infty}e^{-\nu(p)(t-s)} w_{\b}(p)^{-1}  v_{\phi}\s(g,\t)\sqrt{J(q)}|h(s,x-\hat{p}(t-s),p')|d\omega dq dpds\nonumber\\
&\leq C\left(\k+\f{1}{N^{\b-3}}\right)\sup_{0\leq s\leq t}\|h(s)\|^2_{L^\infty}+C_N\int_0^{t-\k}e^{-c_N(t-s)} \|h(s)\|_{L^\infty}\nonumber\\
&~~~~~~~~~~~~~~~~~~~~~~~~~~\times\int_{|p|\leq N}\int_{|q|\leq N}\int_{\mathbb{S}^2}v_{\phi}\s(g,\t)\sqrt{J(q)}|h(s,x-\hat{p}(t-s),p')|d\omega dq dpds.
\end{align}
To estimate the last term on RHS of \eqref{5.17}, we note that
{\small\begin{align}\label{5.17-1}
&\int_{|p|\leq N}\int_{|q|\leq N}\int_{\mathbb{S}^2}v_{\phi}\s(g,\t)\sqrt{J(q)}|h(s,x-\hat{p}(t-s),p')|d\omega dq dp\nonumber\\
&\leq \left(\int_{|p|\leq N}\int_{|q|\leq N}\int_{\mathbb{S}^2}|v_{\phi}\s(g,\t)|^{d}\sqrt{J(q)}d\omega dq dp\right)^{\f1d}\left(\int_{|p|\leq N}\int_{|q|\leq N}\int_{\mathbb{S}^2}\sqrt{J(q)}|h(s,\hat{x},p')|^{\f{d}{d-1}}d\omega dq dp\right)^{1-\f1d}\nonumber\\
&\leq C_N \left(\int_{|p|\leq N}\int_{|q|\leq N}\int_{\mathbb{S}^2}\sqrt{J(q)}|f(s,\hat{x},p')|^{\f{d}{d-1}}d\omega dq dp\right)^{1-\f1d},
\end{align}}
where we have denoted  $\hat{x}=x-\hat{p}(t-s)$  and used \eqref{A.19} in the last inequality.
Noting  $d\omega=\f{\sqrt{s}}{2g}\d^{(4)}(p^\mu+q^\mu-{p^\mu}'-{q^\mu}')\f{dp'dq'}{p_0'q_0'}$, one has that
{\small\begin{align}\label{5.17-2}
&\int_{|q|\leq N}\int_{\mathbb{S}^2}\sqrt{J(q)}|f(s,\hat{x},p')|^{\f{d}{d-1}}d\omega dq\leq C_N\int_{|q|\leq N}\int_{\mathbb{S}^2}\sqrt{J(q)}\sqrt{J(p')}|f(s,\hat{x},p')|^{\f{d}{d-1}}d\omega dq\nonumber\\
&\leq C_N\int_{|q|\leq N}\int_{\mathbb{S}^2}\sqrt{J(q)}\sqrt{J(p')}\Big[|f(s,\hat{x},p')|^{\f{d}{d-1}}I_{\{|F(s,\hat{x},p')-J(p')|<J(p')\}}\nonumber\\
&~~~~~~~~~~~~~~~~~~~~~~~~~~~~~~~~~~~~~~~~~~+|F(s,\hat{x},p')-J(p')|^{\f{d}{d-1}}I_{\{|F(s,\hat{x},p')-J(p')|\geq J(p')\}}\Big]d\omega dq\nonumber\\
&\leq C_N\int_{|q|\leq N}\int_{\mathbb{R}^3}\int_{\mathbb{R}^3}\f{\sqrt{s}}{g}\d^{(4)}(p^\mu+q^\mu-{p^\mu}'-{q^\mu}')J(q)^{\f14}\sqrt{J(p')}\Big[|f(s,\hat{x},p')|^{\f{d}{d-1}}I_{\{|F(s,\hat{x},p')-J(p')|<J(p')\}}\nonumber\\
&~~~~~~~~~~~~~~~~~~~~~~~~~~~~~~~~~~~~~~~~~~~~~~~+|F(s,\hat{x},p')-J(p')|^{\f{d}{d-1}}I_{\{|F(s,\hat{x},p')-J(p')|\geq J(p')\}}\Big]\f{dp'dq'dq}{p_0'q_0'q_0}\nonumber\\
&\leq C_N\int_{\mathbb{R}^3}\int_{\mathbb{R}^3}\int_{\mathbb{R}^3}\f{\sqrt{\bar s}}{\bar g}\d^{(4)}(p^\mu+{p^\mu}'-q^\mu-{q^\mu}') J(p')^{\f14}\sqrt{J(q)}\Big[|f(s,\hat{x},q)|^{\f{d}{d-1}}I_{\{|F(s,\hat{x},q)-J(q)|<J(q)\}}\nonumber\\
&~~~~~~~~~~~~~~~~~~~~~~~~~~~~~~~~~~~~~~~~~~~~~~~~~~+|F(s,\hat{x},q)-J(q)|^{\f{d}{d-1}}I_{\{|F(s,\hat{x},q)-J(q)|\geq J(q)\}}\Big]\f{dp'dq'dq}{p_0'q_0'q_0}\nonumber\\
&\leq C_N\int_{\mathbb{R}^3} \tilde{A}(p,q)\sqrt{J(q)}\Big[|f(s,\hat{x},q)|^{\f{d}{d-1}}I_{\{|F(s,\hat{x},q)-J(q)|<J(q)\}}\nonumber\\
&~~~~~~~~~~~~~~~~~~~~~~~~~~~~~~~~~~~~~~~~~~~~~~~~~~~~~~~~~+|F(s,\hat{x},q)-J(q)|^{\f{d}{d-1}}I_{\{|F(s,\hat{x},q)-J(q)|\geq J(q)\}}\Big]  dq,
\end{align}}
where $\tilde{A}(p,q)$ is defined in \eqref{4.10} above.  Then applying the bound of  $\tilde{A}(p,q)$ in \eqref{4.10}, we have
\begin{align}
&\int_{|q|\leq N}\int_{\mathbb{S}^2}\sqrt{J(q)}|f(s,\hat{x},p')|^{\f{d}{d-1}}d\omega dq\nonumber\\
&\leq  C_N\int_{\mathbb{R}^3}\f{p_0q_0\sqrt{J(q)}}{|p-q|^2(1+|p-q|)}\Big[|f(s,\hat{x},q)|^{\f{d}{d-1}}I_{\{|F(s,\hat{x},q)-J(q)|<J(q)\}}\nonumber\\
&~~~~~~~~~~~~~~~~~~~~~~~~~~~~~~~~~~~~~~~~~~~~~~+|F(s,\hat{x},q)-J(q)|^{\f{d}{d-1}}I_{\{|F(s,\hat{x},q)-J(q)|\geq J(q)\}}\Big] dq\nonumber\\
&\leq C_N\Big(\int_{\mathbb{R}^3}\sqrt{J(q)}|f(s,\hat{x},q)|^{\f{5d}{d-1}}I_{\{|F(s,\hat{x},q)-J(q)|<J(q)\}}dq\nonumber\\
&~~~~~~~~~~~~~~~~~~~~~~~~~~~~~~~~~~~~+\int_{\mathbb{R}^3}|F(s,\hat{x},q)-J(q)|^{\f{5d}{d-1}}I_{\{|F(s,\hat{x},q)-J(q)|\geq J(q)\}} dq\Big)^{\f15}\nonumber\\
&\leq C_N\|f(s)\|_{L^\infty}^{\f{4d+1}{5(d-1)}}\Big(\int_{\mathbb{R}^3}\sqrt{J(q)}|f(s,\hat{x},q)|I_{\{|F(s,\hat{x},q)-J(q)|<J(q)\}}dq\nonumber\\
&~~~~~~~~~~~~~~~~~~~~~~~~~~~~~~~~~~~~~~~~~~~+\int_{\mathbb{R}^3}|F(s,\hat{x},q)-J(q)|I_{\{|F(s,\hat{x},q)-J(q)|\geq J(q)\}} dq\Big)^{\f15},\nonumber
\end{align}
which, together with \eqref{5.17-1} and \eqref{2.52}, yields that
\begin{align}\label{5.17-3}
&\int_{|p|\leq N}\int_{|q|\leq N}\int_{\mathbb{S}^2}v_{\phi}\s(g,\t)\sqrt{J(q)}|h(s,x-\hat{p}(t-s),p')|d\omega dq dp\nonumber\\
&\leq C_{\k,N}\|f(s)\|_{L^\infty}^{\f{4d+1}{5d}} \Big[\sqrt{\mathcal{E}(F_0)}+\mathcal{E}(F_0)\Big]^{\f{d-1}{5d}}.
\end{align}
Substituting \eqref{5.17-3} into \eqref{5.17}, we get, for $\b>5$, that
\begin{equation}\label{5.20}
H_{22}\leq C\Big(\k+\f{1}{N^{\b-3}}\Big)\sup_{0\leq s\leq t}\|h(s)\|^2_{L^\infty}+C_{\k,N}\|h(s)\|_{L^\infty}^{\f{9d+1}{5d}} \Big[\sqrt{\mathcal{E}(F_0)}+\mathcal{E}(F_0)\Big]^{\f{d-1}{5d}}.
\end{equation}
Combining \eqref{5.20} and \eqref{5.16}, one has
\begin{align}\label{5.20-1}
H_2&\leq C\Big(\k+\f1N+\f{1}{N^{\b-3}}\Big)\sup_{0\leq s\leq t}\|h(s)\|^2_{L^\infty}
+C_{\k,N}\Big[\sqrt{\mathcal{E}(F_0)}+\mathcal{E}(F_0)\Big]^{1-\f1d}\sup_{0\leq s\leq t}\|h(s)\|^{1+\f1d}_{L^\infty}\nonumber\\
&~~~~+C_{\k,N}\|h(s)\|_{L^\infty}^{\f{9d+1}{5d}} \Big[\sqrt{\mathcal{E}(F_0)}+\mathcal{E}(F_0)\Big]^{\f{1}{5}(1-\f1d)}.
\end{align}
Submitting \eqref{5.11} and \eqref{5.20-1}  into \eqref{5.3}, one proves \eqref{5.2}  for $\b> 5$. Therefore the proof of Lemma \ref{lem5.1} is completed. $\hfill\Box$

\

\subsection{Proof of Theorem \ref{thm1.1}}

Let $\b>14$,  we make the {\it a priori} assumption
\begin{align}\label{5.22}
\|h(t)\|_{L^\infty}\leq 4R_0:= 4C_2\Big\{2\bar{M}^2+\sqrt{\mathcal{E}(F_0)}+\mathcal{E}(F_0)\Big\},
\end{align}
where the positive constant $C_2\geq1$ is defined in Lemma \ref{lem4.2}. Then it follows from Lemma \ref{lem4.2} and the a priori assumption \eqref{5.22} that
\begin{align}\label{5.23}
\|h(t)\|_{L^\infty}\leq R_0+C_2(4R_0)^{\f{9d+1}{5d}}\cdot\sup_{t_1\leq s\leq t,z\in\Omega} \Big(\int_{\mathbb{R}^3}|f(s,z,q)|dq\Big)^{\f{d-1}{5d}}.
\end{align}
For $t\geq t_1$,   it holds that
\begin{align}\label{5.40}
\int_{\mathbb{R}^3}e^{-\nu(p)t}|f_0(x-\hat{p}t,p)|dp&\leq C\|w_\b f_0\|_{L^\infty} \tilde{N}^{-\b+3}+\int_{|p|\leq \tilde{N}}|f_0(x-\hat{p}t,p)|dp\nonumber\\
&\leq C\|w_\b f_0\|_{L^\infty} \tilde{N}^{-\b+3}+\int_{|p|\leq \tilde{N}}|f_0(x-\hat{p}t,p)|(1+|p|^2)^{\f52}d\hat{p}\nonumber\\
&\leq C\|w_\b f_0\|_{L^\infty} \tilde{N}^{-\b+3}+C\tilde{N}^5\int_{|\hat{p}|\leq 1}\|f_0(x-\hat{p}t,\cdot)\|_{L^\infty}d\hat{p}\nonumber\\
&\leq C\|w_\b f_0\|_{L^\infty} \tilde{N}^{-\b+3}+C\tilde{N}^5(1+t_1^{-3})\|f_0\|_{L^1_xL^\infty_p}\nonumber\\
&\leq C\|w_\b f_0\|_{L^\infty} \tilde{N}^{-\b+3}+C\tilde{N}^5\bar{M}^3\|f_0\|_{L^1_xL^\infty_p}\nonumber\\
&\leq C \bar{M}^{\f{3\b-4}{\b+2}} \|f_0\|_{L^1_xL^\infty_p}^{\f{\b-3}{\b+2}}.
\end{align}
where we have chosen $\tilde{N}=\bar{M}^{-\f{2}{\b+2}}\|f_0\|_{L^1_xL^\infty_p}^{-\f1{\b+2}}$ above.
Then it follows from  Lemma \ref{lem5.1}, \eqref{5.40} and the a priori assumption \eqref{5.22} that
\begin{align}\label{5.24}
&\sup_{t_1\leq s\leq t,z\in\mathbb{R}^3}\int_{\mathbb{R}^3}|f(t,z,p)|dp\nonumber\\
&\leq C\bar{M}^{\f{3\b-4}{\b+2}} \|f_0\|_{L^1_xL^\infty_p}^{\f{\b-3}{\b+2}}+C_{\k,N}\Big[\sqrt{\mathcal{E}(F_0)}+\mathcal{E}(F_0)\Big]+C\Big[\k+\f{1}{N}\Big](4R_0)^2\nonumber\\
&~~~+C_{\k,N}\Big(\sqrt{\mathcal{E}(F_0)}+\mathcal{E}(F_0)\Big)^{\f15(1-\f1d)}(4R_0)^2.
\end{align}
Noting $\b>14, d>1$, firstly choosing $\k$ sufficiently small, then $N\geq1$ large enough, finally letting $\mathcal{E}(F_0)+\|f_0\|_{L^1_xL^\infty_p}\leq \epsilon_0$ with $\epsilon_0$ small depending only on $a,b,\b,\g$ and $\bar{M}$ such that
\begin{align}\label{5.26}
16C_2R_0^{\f{4d+1}{5d}}\cdot\sup_{t_1\leq s\leq t,z\in\Omega} \Big(\int_{\mathbb{R}^3}|f(s,z,p)|dp\Big)^{\f{d-1}{5d}}\leq \f34,
\end{align}
which, together with \eqref{5.23}, yields immediately that
\begin{equation}\label{5.27}
\|h(t)\|_{L^\infty}\leq \f74R_0~~\mbox{for}~~t\geq0.
\end{equation}
That is, the a priori assumption \eqref{5.22} is closed and therefore  the proof of Theorem \ref{thm1.1} is completed. $\hfill\Box$



\section{Decay Estimates in Torus $\mathbb{T}^3$}

In this section, we try to obtain the decay rates for the global  solutions  obtained in Theorem \ref{thm1.1} in the case $\Omega=\mathbb{T}^3$.  Consider the following linearized Boltzmann equation
\begin{align}\label{6.1}
\varphi_t+\hat{p}\cdot\nabla_x\varphi+\nu(p)\varphi-K\varphi=0,~~\varphi(0,x,p)=\varphi_0(x,p).
\end{align}
Denote the semigroup for \eqref{6.1} by  $S(t)$, then we have
\begin{equation*}\label{6.1-3}
\varphi(t)=S(t)\varphi_0.
\end{equation*}
Let $\varphi(t,x,p)$ be the solution of linearized equation \eqref{6.1}  and define
\begin{equation*}\label{6.1-2}
\Phi(t,x,p):=w_\b(p)\varphi(t,x,p),
\end{equation*}
then it holds
\begin{align}\label{6.1-1}
\Phi_{t}+\hat{p}\cdot\nabla_x\Phi+\nu(p) \Phi-K_{w_\b}\Phi=0,~~\Phi(0,x,p)=\Phi_{0}(x,p).
\end{align}
For later use, we denote the semigroup for \eqref{6.1-1} by $U(t)$, then one has
\begin{equation*}\label{6.5}
\Phi(t)=U(t)\Phi_{0}.
\end{equation*}

From the $H$-theorem, $L$ is a nonnegative, and the null space of $L$ is given by the five dimensional space \cite{Glassey}:
\begin{align}\label{6.5-1}
\mathcal{N}=span\left\{\sqrt{J}, p\sqrt{J}, p_0\sqrt{J} \right\}.
\end{align}
We define $\mathbf{P}$ as its $p$-projection in $L^2(\mathbb{R}^3)$ to the null space $\mathcal{N}$. Then we decompose $f(t,x,p)$ uniquely as
\begin{equation*}\label{6.5-2}
f(t,x,p)=\mathbf{P}f+\mathbf{(I-P)}f.
\end{equation*}
Furthermore, we expand $\mathbf{P}f$ as a linear combination of the basis in \eqref{6.5-1},
\begin{equation*}\label{6.5-3}
\mathbf{P}f=\left\{a(t,x)+\sum_{i=1}^3b_i(t,x)p_i+c(t,x)p_0 \right\}\sqrt{J}.
\end{equation*}

Using \eqref{2.47}, \eqref{2.48}  and Lemma 3.5.1 of \cite{Glassey}, it is easy to check  that $K$ is a  compact operator, from $L^2(\mathbb{R}^3_p)$ to $L^2(\mathbb{R}^3_p)$ . Then the following lemma holds:
\begin{lemma}\label{lem6.0}
Under the hypothesis {\bf H)}, there exists a positive constant $\d_0>0$ such that
\begin{equation*}\label{6.5-4}
\langle f,Lf \rangle\geq \d_0 |\mathbf{(I-P)}f|^2_{\nu},
\end{equation*}
and $Lf=0$ if and only if $f=\mathbf{P}f$.
\end{lemma}

Guo \cite{Guo2} firstly established the following lemma for the Newtonian Boltzmann equation. Later, Strain \cite{Strain1} extended it to the relativistic Boltzmann equation in the case of soft potentials.  Indeed, the lemma also holds for relativistic Boltzmann equation in the case of hard potentials, the proof is almost the same as in \cite{Strain}, we omit the details of proof for simplicity of presentation.
\begin{lemma}\label{lem6.0-1}
Let  $f(t,x,p)$ be any solution to the linearized relativistic Boltzmann equation \eqref{6.1} in the sense of distribution, and $f$ satisfies \eqref{1.17}-\eqref{1.19} with $(M_0,\tilde{M}_0,E_0)=(0,0,0)$,  then there exists a positive constant $M>0$ such that
\begin{equation*}
\int_0^1\|\mathbf{(I-P)}f(s)\|_{\nu}^2ds\geq M \int_0^1\|\mathbf{P}f(s)\|^2_{\nu}ds.
\end{equation*}
\end{lemma}

\subsection{Decay Estimate for Hard Potentials}
In this subsection, we consider the decay estimate for the case of hard potentials on torus $\mathbb{T}^3$, i.e.,  $\g>-\f43, a\in[0,2]\cap[0,\min\{2+\g,4+3\g\}),b\in[0,2)$.

\

Utilizing Lemmas \ref{lem6.0} and  \ref{lem6.0-1} and by similar arguments as in Theorem 5 of \cite{Guo2}, we have the following $L^2$-exponential decay estimate for the linearized equation \eqref{6.1}.
\begin{lemma}\label{lem6.1}
Let $\Omega=\mathbb{T}^3$,  $\varphi(t,x,p)$ be any solution to the linearized Boltzmann equation \eqref{6.1} and satisfies the conservations of mass \eqref{1.17}, momentum \eqref{1.18} and energy \eqref{1.19} with $(M_0,\tilde{M}_0,E_0)=(0,0,0)\in \mathbb{R}\times\mathbb{R}^3\times\mathbb{R}$. Then there exist  positive constants $\l_1>0$ and $C>0$ such that
	\begin{align}\label{6.2}
	\|S(t)\varphi_0\|_{L^2}=\|\varphi(t)\|_{L^2}\leq Ce^{-\l_1 t}\|\varphi_0\|_{L^2},~~\mbox{for}~~t\geq0.
	\end{align}
\end{lemma}

Utilizing  Lemma \ref{lem6.1}, we can obtain the following $L^\infty$ decay estimate for the linearized Boltzmann equation.
\begin{lemma}\label{lem6.2}
	Let $\Omega=\mathbb{T}^3$ and $\beta>\f32$. Let $\varphi(t,x,v)$ be any solution to the linearized Boltzmann equation \eqref{6.1} and satisfies the conservations of mass \eqref{1.17}, momentum \eqref{1.18} and energy \eqref{1.19} with $(M_0,\tilde{M}_0,E_0)=(0,0,0)\in \mathbb{R}\times\mathbb{R}^3\times\mathbb{R}$. Then there exist  positive constants $0<\l_2\leq \l_1$ and $C>0$ such that
	\begin{align}\label{6.6}
	\|U(t)\Phi_0\|_{L^\infty}=\|\Phi(t)\|_{L^\infty}\leq Ce^{-\l_2 t}\|w_\b \varphi_0\|_{L^\infty},~~\mbox{for}~~t\geq0.
	\end{align}
\end{lemma}

\noindent{\bf Proof.} The mild form of \eqref{6.1-1} is given by
\begin{equation*}\label{6.4-1}
\Phi(t,x,p)=e^{-\nu(v)t}\Phi_0(x-\hat{p}t,p)
+\int_0^te^{-\nu(p)(t-s)} (K_{w_{\b}} \Phi)(s,x-\hat{p}(t-s),p)ds.
\end{equation*}
Noting  $\nu(p)\geq \nu_0>0$, by similar arguments as in Lemma \ref{lem4.2}, we can obtain
\begin{align}\label{6.4-2}
\|\Phi(t)\|_{L^\infty}&\leq C(1+t)e^{-\nu_0t}\|\Phi_0\|_{L^\infty}+C(\k+\f1{N^{\xi}})e^{-\f{\nu_0 t}{2}}\sup_{0\leq s\leq t}\left(e^{\f{\nu_0}2s}\|\Phi(s)\|_{L^\infty}\right)\nonumber\\
&~~+C_{N}\int_0^t\int_0^{s-\k}e^{-c_N(t-s)}e^{-c_N(s-\tau)}  \int_{|q|\leq 2N,|\eta|\leq3N} |\Phi(\tau,\hat{x}-\hat{q}(s-\tau),\eta)|dqd\eta d\tau ds
\end{align}
A direct calculation shows that
\begin{align}
\int_{|q|\leq 2N,|\eta|\leq3N} |\Phi(\tau,\hat{x}-\hat{q}(s-\tau),\eta)|dqd\eta&\leq C_N\int_{|q|\leq 2N,|\eta|\leq3N} |\varphi(\tau,\hat{x}-\hat{q}(s-\tau),\eta)|dqd\eta\nonumber\\
&\leq C_{\k,N}\|\varphi(s)\|_{L^2},\nonumber
\end{align}
which, together with \eqref{6.4-2}, yields that
\begin{align}
\sup_{0\leq s\leq t}\left(e^{\f{\nu_0 t}{2}}\|\Phi(t)\|_{L^\infty}\right)&\leq C(1+t)\|\Phi_0\|_{L^\infty}+C(\k+\f1{N^{\xi}})\sup_{0\leq s\leq t}\left\{e^{\f{\nu_0}2s}\|\Phi(s)\|_{L^\infty}\right\}\nonumber\\
&~~~+C_{\k,N}\int_0^t\|\varphi(s)\|_{L^2} ds.\nonumber
\end{align}
We first choose $T_0$ large so that $2C(1+T_0)e^{-\f{\nu_0}{2}T_0}=e^{-\f{\nu_0}{4} T_0}$, then choose $N$ large and $\k$ sufficiently small, thus we have
\begin{equation}\label{6.4-5}
\|\Phi(T_0)\|_{L^\infty}\leq e^{-\f{\nu_0}{4}T_0}\|\Phi_0\|_{L^\infty}+C_{T_0}\int_0^{T_0}\|\varphi(s)\|_{L^2}ds.
\end{equation}
Thus, \eqref{6.6} follows from \eqref{6.2}, \eqref{6.4-5} and   Lemma 19 of  Guo \cite{Guo2}. Here we omit the details of proof for simplicity of presentation. Therefore  the proof of this lemma is completed.  $\hfill\Box$

\

Based on the above preparations, we apply  Lemma \ref{lem6.2} to prove Theorem \ref{thm1.2}.\\[1.5mm]
\noindent{\bf Proof of Theorem \ref{thm1.2}:}  Using the semigroup $U(t)$ of \eqref{6.1-1} and  the Duhamel Principle,   we obtain the formula  of solutions to the weighted relativistic Boltzmann equation \eqref{4.1} as
\begin{equation*}\label{6.8}
h(t)=U(t)h_0+\int_0^tU(t-s)\Big\{w_{\b}\Gamma(f,f)(s)\Big\}ds.
\end{equation*}
Then it follows from \eqref{6.6} that
\begin{align}\label{6.9}
\|h(t)\|_{L^\infty}
&\leq Ce^{-\l_2 t}\|h_0\|_{L^\infty}+\Big\|\int_0^tU(t-s)\Big\{w_{\b}\Gamma(f,f)(s)\Big\}ds\Big\|_{L^\infty}.
\end{align}
To bound the last term on the RHS of \eqref{6.9}, we notice that
\begin{align}\label{6.10}
&\int_0^tU(t-s)\Big\{w_{\b}\Gamma(f,f)(s)\Big\}ds=\int_0^te^{-\nu(p)(t-s)}\Big\{w_{\b}\Gamma(f,f)(s)\Big\}ds\nonumber\\
&~~~~~~~~~~~~~~~~~~~~~~~~~~~~~~~+\int_0^t\int_s^te^{-\nu(p)(t-s_1)}K_{w_{\b}}\Big\{U(s_1-s)w_{\b}\Gamma(f,f)(s)\Big\}ds_1ds.
\end{align}
For the first term on the RHS of \eqref{6.10}, it follows from \eqref{4.6} that
\begin{align}\label{6.11}
&\Big|\int_0^te^{-\nu(p)(t-s)}\Big\{w_{\b}(p)\Gamma(f,f)(s)\Big\}ds\Big|\nonumber\\
&\leq C\int_0^te^{-\nu(p)(t-s)}\nu(p)\|h(s)\|_{L^\infty}^{\f{9d+1}{5d}}\sup_{z\in\Omega}\Big(\int_{\mathbb{R}^3}|f(s,z,\eta)|d\eta\Big)^{\f{d-1}{5d}}ds\nonumber\\
&\leq C\int_0^te^{-\nu(p)(t-s)}\nu(p) e^{-\f{\l_2}2s}{\small\sup_{0\leq s\leq t, z\in\Omega}
	 \left\{\Big[e^{\f{\l_2}2s}\|h(s)\|_{L^\infty}\Big]\cdot\|h(s)\|^{\f{4d+1}{5d}}_{L^\infty}\Big(\int_{\mathbb{R}^3}|f(s,z,\eta)|d\eta\Big)^{\f{d-1}{5d}}\right\}}ds\nonumber\\
&\leq Ce^{-\f{\l_2}{2} t}\sup_{0\leq s\leq t, z\in\Omega}
\left\{\Big[e^{\f{\l_2}2s}\|h(s)\|_{L^\infty}\Big]\cdot\|h(s)\|^{\f{4d+1}{5d}}_{L^\infty}\Big(\int_{\mathbb{R}^3}|f(s,z,\eta)|d\eta\Big)^{\f{d-1}{5d}}\right\}.
\end{align}

For the second term on the RHS of \eqref{6.10},  motivated by \cite{Guo2}, we define a new semigroup $\tilde{U}(t)$ solving
\begin{equation*}\label{6.12}
\Big\{\partial_t+\hat{p}\cdot\nabla_x+\nu(p)-K_{\tilde{w}}\Big\}\{\tilde{U}(t)\tilde{h}_0\}=0,~~\tilde{U}(0)\tilde{h}_0=\tilde{h}_0,
\end{equation*}
with $\tilde{w}(p)=\f{w_{\b}(p)}{(1+|p|^2)^{\f{a}4}}$. A direct calculation shows that $(1+|p|^2)^{\f{a}4}\tilde{U}(t)$
also solves the  weighted linearized Boltzmann equation \eqref{6.1-1}. The uniqueness in $L^\infty$ with $\tilde{h}_0=\f{h_0}{(1+|p|^2)^{\f{a}4}}$ yields that
\begin{equation*}\label{6.13}
U(t)h_0\equiv (1+|p|^2)^{\f{a}4}\tilde{U}(t)\Big\{\f{h_0}{(1+|p|^2)^{\f{a}4}}\Big\}.
\end{equation*}
We note that  \eqref{6.6} also holds for semigroup $\tilde{U}(t)$.
Then it follows from \eqref{6.6} and \eqref{4.6} that
{\small\begin{align}\label{6.14}
&\Big|\int_0^t\int_s^te^{-\nu(p)(t-s_1)}K_{w_{\b}}\Big\{U(s_1-s)w_{\b}\Gamma(f,f)(s)\Big\}ds_1ds\Big|\nonumber\\
&\leq \int_0^t\int_s^te^{-\nu(p)(t-s_1)}\Big|\int_{\mathbb{R}^3_{q}}k_{w_{\b}}(p,q)(1+|q|^2)^{\f{a}4}\Big\{\tilde{U}(s_1-s)\f{w_{\b}}{(1+|q|^2)^{\f{a}4}}\Gamma(f,f)(s)\Big\}\Big|dqds_1ds\nonumber\\
&\leq \int_0^t\int_s^te^{-\nu(p)(t-s_1)}\Big|\int_{\mathbb{R}^3_{q}}k_{w_{\b}}(p,q)(1+|q|^2)^{\f{a}4}dq\Big|\cdot\Big\|\Big\{\tilde{U}(s_1-s)\f{w_{\b}}{(1+|p|^2)^{\f{a}4}}\Gamma(f,f)(s)\Big\}\Big\|_{L^\infty}dqds_1ds\nonumber\\
&\leq \int_0^t\int_s^te^{-\nu(p)(t-s_1)}e^{-\l_2(s_1-s)}\nu(p)\Big\|\f{w_{\b}(p)}{(1+|p|^2)^{\f{a}4}}\Gamma(f,f)(s)\Big\|_{L^\infty}ds_1 ds\nonumber\\
&\leq C\int_0^t\int_s^te^{-\nu(p)(t-s_1)}e^{-\l_2(s_1-s)}\nu(p)\Big\|\f{\nu(p)}{(1+|p|^2)^{\f{a}4}}\Big\|_{L^\infty}\|h(s)\|_{L^\infty}^{\f{9d+1}{5d}}\sup_{z\in\Omega}\Big(\int_{\mathbb{R}^3}|f(s,z,\eta)|d\eta\Big)^{\f{d-1}{5d}}ds_1ds\nonumber\\
&\leq C\int_0^t\int_s^te^{-\nu(p)(t-s_1)}e^{-\l_2(s_1-s)-\f{\l_2}{2}s}\nu(p)ds_1ds\sup_{0\leq s\leq t, z\in\Omega}
\left\{\Big[e^{\f{\l_2}2s}\|h(s)\|_{L^\infty}\Big]\cdot\|h(s)\|^{\f{4p+1}{5p}}_{L^\infty}\Big(\int_{\mathbb{R}^3}|f(s,z,\eta)|d\eta\Big)^{\f{d-1}{5d}}\right\}\nonumber\\
&\leq Ce^{-\f{\l_2}{2} t}\sup_{0\leq s\leq t, z\in\Omega}
\left\{\Big[e^{\f{\l_2}2s}\|h(s)\|_{L^\infty}\Big]\cdot\|h(s)\|^{\f{4p+1}{5p}}_{L^\infty}\Big(\int_{\mathbb{R}^3}|f(s,z,\eta)|d\eta\Big)^{\f{d-1}{5d}}\right\},
\end{align}}
where we have used the following bound in the last inequality 
\begin{align}
&\int_0^t\int_s^te^{-\nu(p)(t-s_1)}e^{-\l_2(s_1-s)}e^{-\f{\l_2}2s}\nu(p)ds_1 ds
=\int_0^te^{-\nu(p)t}e^{\f{\l_2}2s}ds\int_s^te^{[\nu(p)-\l_2]s_1}\nu(p)ds_1 \nonumber\\
&\leq C\int_0^te^{-\nu(p)t}e^{\f{\l_2}2s}ds\int_s^te^{[\nu(p)-\l_2]s_1}d(\nu(p)-\l_2)s_1\nonumber\\
&\leq C\int_0^te^{-\nu(p)t}e^{\f{\l_2}2s}e^{[\nu(p)-\l_2]t}ds\leq Ce^{-\f{\l_2}2t}.\nonumber
\end{align}

Combining \eqref{6.9}-\eqref{6.11}, \eqref{6.14} and using \eqref{5.27}, one obtains that
\begin{align}\label{6.15}
&\sup_{0\leq s\leq t}\Big\{e^{\f{\l_2}{2}s}\|h(s)\|_{L^\infty}\Big\}\nonumber\\
&\leq C\|h_0\|_{L^\infty}+C\sup_{0\leq s\leq t, z\in\Omega}
\left\{\Big[e^{\f{\l_2}2s}\|h(s)\|_{L^\infty}\Big]\cdot\|h(s)\|^{\f{4d+1}{5d}}_{L^\infty}\Big(\int_{\mathbb{R}^3}|f(s,z,\eta)|d\eta\Big)^{\f{d-1}{5d}}\right\}\nonumber\\
&\leq C\Big\{\|h_0\|_{L^\infty}+\sup_{0\leq s\leq 1}\|h(s)\|^2_{L^\infty}\Big\}\nonumber\\
&~~~~~~~~~~~~+C\sup_{1\leq s\leq t}\Big[e^{\f{\l_2}2s}\|h(s)\|_{L^\infty}\Big]
\cdot\sup_{1\leq s\leq t,z\in\Omega}\left\{\|h(s)\|^{\f{4d+1}{5d}}_{L^\infty}\Big(\int_{\mathbb{R}^3}|f(s,z,\eta)|d\eta\Big)^{\f{d-1}{5d}}\right\}\nonumber\\
&\leq C_3\bar{M}^4+C_3\sup_{1\leq s\leq t}\Big[e^{\f{\l_2}2s}\|h(s)\|_{L^\infty}\Big]
\cdot\sup_{1\leq s\leq t,z\in\Omega}\left\{\|h(s)\|^{\f{4d+1}{5d}}_{L^\infty}\Big(\int_{\mathbb{R}^3}|f(s,z,\eta)|d\eta\Big)^{\f{d-1}{5d}}\right\}.
\end{align}
Using \eqref{5.24} and by  similar arguments as in \eqref{5.26},  if $\epsilon_0$ is small, one gets
\begin{align}\label{6.20}
&C_3\sup_{1\leq s\leq t,z\in\Omega}\left\{\|h(s)\|^{\f{4d+1}{5d}}_{L^\infty}\Big(\int_{\mathbb{R}^3}|f(s,z,\eta)|d\eta\Big)^{\f{d-1}{5d}}\right\}\leq \f12.
\end{align}
Substituting \eqref{6.20} into \eqref{6.15} gives that
\begin{align}\label{6.21}
e^{\f{\l_2}{2}t}\|h(t)\|_{L^\infty}
\leq 2C_3\bar{M}^4,~~t\geq0.
\end{align}
Choosing
$$\l_0=\f{\l_2}{2}~~~\mbox{and}~~~\tilde{C}_2=2C_3\bar{M}^4,$$
we obtain \eqref{1.29} from \eqref{6.21}.  Therefore the proof of Theorem \ref{thm1.2} is completed. $\hfill\Box$

\subsection{Decay Estimate for Soft Potentials}

In this subsection, we consider the decay estimate for soft potentials on torus, i.e. $b\in(0,2), \g>-\min\{\f43,4-2b\}$. Firstly, we give the following estimate on the linearized operator $L$. The proof of following lemma will be given in the appendix since it is similar to the ones of \cite{Strain}.
\begin{lemma}\label{lem8.1}
There exists a constant $R\geq1$, such that the following inequality holds
\begin{align}\label{6.22}
\langle  w_{\vartheta}^2Lf,f\rangle\geq \f12|w_\vartheta f|^2_{\nu}-C_{R}|I_{\leq R}f|_{L^2}^2,
\end{align}
where the positive constant $C_{R}>0$ depends only on $R$.
\end{lemma}
Using Lemma \ref{lem6.0-1} and  Lemma \ref{lem8.1}, and by same arguments as in the proof of Theorem 3.7 of  \cite{Strain}, we can obtain the following $L^2$-decay estimate. Here we omit the details of the proof for simplicity of presentation.
\begin{lemma}\label{lem8.2}
Let $\Omega=\mathbb{T}^3$ and  $\varphi(t,x,p)$ be any solution to the linearized Boltzmann equation \eqref{6.1} and satisfies  \eqref{1.17}-\eqref{1.19} with $(M_0,\tilde{M}_0,E_0)=(0,0,0)\in \mathbb{R}\times\mathbb{R}^3\times\mathbb{R}$. Then, for any $\a,k\geq 0$, there exist  positive constants $C_{\a,k}>0$ such that
\begin{align}\label{6.23}
\|w_\a\varphi(t)\|_{L^2}\leq C_{\a,k}(1+t)^{-k}\|w_{\a+k\f{b}{2}}\varphi_0\|_{L^2}.
\end{align}
\end{lemma}

Based on Lemma \ref{lem8.2}, we  have the following $L^\infty$-decay estimate for the solutions to the  linearized relativistic Boltzmann equation \eqref{6.1}.
\begin{lemma}\label{lem8.3}
	Under the assumptions of  Lemma \ref{lem8.2}, given $\vartheta>\f32$ and $k\in[0,1+\f{\xi_1}b]$, then it holds  that
	\begin{align}\label{8.5}
	\|S(t)\varphi_0\|_{L^\infty}=\|\varphi(t)\|_{L^\infty}\leq C_{\vartheta,k}(1+t)^{-k}\|w_{\vartheta+k\f{b}{2}}\varphi_0\|_{L^\infty},
	\end{align}
where $\x_1>0$ is defined in Lemma \ref{lem2.4} and  the positive constant $C_{\vartheta,k}$ depending only on $\vartheta,k$.
\end{lemma}

\noindent{\bf Proof.}
It is noted that
\begin{align}\label{8.6}
\varphi(t,x,p)= e^{-\nu(p)t}\varphi_0(x-\hat{p}t,p)+\int_0^t e^{-\nu(p)(t-s)}K\varphi(s,x-\hat{p}(t-s),p)ds.
\end{align}
As in Vidav \cite{Vidav,Strain}, we iterate \eqref{8.6} again  to obtain that
\begin{align}\label{8.7}
\varphi(t,x,p)
&= e^{-\nu(p)t}\varphi_0(x-\hat{p}t,p)+\int_0^t e^{-\nu(p)(t-s)}\int_{\mathbb{R}^3}K(p,q) e^{-\nu(q)s}\varphi_0(\hat{x}-\hat{q}s,q)
dqds\nonumber\\
&~+\int_0^t e^{-\nu(p)(t-s)}\int_{\mathbb{R}^3}K(p,q)\int_0^se^{-\nu(q)(s-\tau)}\int_{\mathbb{R}^3}K(q,\eta)\varphi(\tau,\hat{x}-\hat{q}(s-\tau),\eta)d\eta d\tau dqds\nonumber\\
&:=L_1+L_2+L_3,
\end{align}
where we have denoted $\hat{x}=x-\hat{p}(t-s)$. Firstly, it is straightforward to get that
\begin{align}\label{8.8}
|L_1|\leq C(1+t)^{-k}\|\nu(p)^{-k}\varphi_0\|_{L^\infty}\leq C(1+t)^{-k}\|w_{k\f{b}{2}}\varphi_0\|_{L^\infty}.
\end{align}
For $L_2$, we decompose it to be
\begin{align}\label{8.9}
L_2&=\int_0^t e^{-\nu(p)(t-s)}\int_{2|p|\geq |q|}K(p,q) e^{-\nu(q)s}\varphi_0(\hat{x}-\hat{q}s,q)dqds\nonumber\\
&~~+\int_0^t e^{-\nu(p)(t-s)}\int_{|q|\geq2|p|}K(p,q) e^{-\nu(q)s}\varphi_0(\hat{x}-\hat{q}s,q)dqds:=L_{21}+L_{22}.
\end{align}
A direct calculation shows that
\begin{align}\label{8.10}
|L_{21}|&\leq \|w_{k\f{b}2}\varphi_0\|_{L^\infty}\int_0^t e^{-c\nu(p)t}\int_{2|p|\geq |q|}|K(p,q)|w_{k\f{b}2}(q)^{-1}dqds\nonumber\\
&\leq C\|w_{k\f{b}2}\varphi_0\|_{L^\infty}\int_0^te^{-c\nu(p)t}p_0^{-\xi_1} \cdot \nu(p)^{1+k}ds\leq C(1+t)^{k}\|w_{k\f{b}2}\varphi_0\|_{L^\infty}.
\end{align}
To estimate $L_{22}$, for $|q|\geq2|p|$,  we note $$e^{-\f1{20}|p-q|}\leq   e^{-\f1{20}(|q|-|p|)}\leq e^{-\f1{40}|q|}\leq C,$$
which yields immediately that
\begin{align}\label{8.11}
|L_{22}|&\leq \|\varphi_0\|_{L^\infty}\int_0^t e^{-\nu(q)t}\int_{|q|\geq2|p|}\Big|K(p,q)e^{\f1{20}|p-q|}\Big| e^{-\f1{20}|p-q|}dqds\nonumber\\
&\leq  \|\varphi_0\|_{L^\infty}\int_0^t e^{-\nu(q)t} e^{-\f1{40}|q|}\int_{|q|\geq2|p|}\Big|K(p,q)e^{\f1{20}|p-q|}\Big|dqds\leq C(1+t)^{-k}\|\varphi_0\|_{L^\infty}.
\end{align}
Thus, it follows from \eqref{8.9}-\eqref{8.11} that
\begin{align}\label{8.12}
|L_2|\leq C(1+t)^{-k}\|w_{k\f{b}2}\varphi_0\|_{L^\infty}.
\end{align}

We now focus on the term $L_3$.  As in the previous, we divide the proof into three cases.\\[1mm]
\noindent{\it Case 1.} For $|p|\geq N$, it follows from \eqref{2.47} that
\begin{align}\label{8.15-2}
L_{3}&\leq\int_0^t e^{-\nu(p)(t-s)}\int_{\mathbb{R}^3}|K(p,q)|
\int_{\mathbb{R}^3}|K(q,\eta)|\int_0^se^{-\nu(q)(s-\tau)}\varphi(\tau,\hat{y},\eta)d\tau d\eta dq  ds \nonumber\\
&\leq C\sup_{0\leq \tau\leq t}\{(1+\tau)^k\|\varphi(\tau)\|_{L^\infty}\}\int_0^t e^{-\nu(p)(t-s)}\int_{\mathbb{R}^3}K(p,q)  q_0^{-\f{b}2-\xi_1}\int_0^s\f{e^{-\nu(q)(s-\tau)}}{(1+\tau)^{k}}d\tau dqds\nonumber\\
&\leq C\sup_{0\leq \tau\leq t}\{(1+\tau)^k\|\varphi(\tau)\|_{L^\infty}\}\int_0^t e^{-\nu(p)(t-s)}p_0^{-\f{b}2-\xi_1}\int_0^s\f{(1+\tau)^{k}}{(1+s-\tau)^{1+\f{2\xi_1}{b}}}d\tau\nonumber\\
&\leq C\sup_{0\leq \tau\leq t}\{(1+\tau)^k\|\varphi(\tau)\|_{L^\infty}\}p_0^{-\f{\xi_1}b} \int_0^t (1+t-s)^{-1-\f{\xi_1}{b}}(1+s)^{-k}ds\nonumber\\
&\leq \f{C}{N^{\f{\xi_1}{b}}}(1+t)^{-k}\sup_{0\leq \tau\leq t}\{(1+\tau)^k\|\varphi(\tau)\|_{L^\infty}\},
\end{align}
where we have denoted $\hat{y}:=\hat{x}-\hat{q}(s-\tau)$.\\

\noindent{\it Case 2.} For $|p|\leq N, |q|\geq 2N$(or $|p|\leq N, |q|\leq 2N, |\eta|\geq3N$), we have $|p-q|\geq N$(or $|q-\eta|\geq N$). Then it follows from \eqref{4.31} that
\begin{align}\label{8.15-3}
&\int_0^t e^{-\nu(p)(t-s)}\left\{\int_{|p|\leq N, |q|\geq 2N}+\int_{|q|\leq 2N, |\eta|\geq3N}\right\}|K(p,q)K(q,\eta)|\int_0^se^{-\nu(q)(s-\tau)}\varphi(\tau,\hat{y},\eta)d\tau d\eta dq  ds \nonumber\\
&\leq e^{-\f{N}{20}}\sup_{0\leq\tau\leq t}\{(1+\tau)^k\|\varphi(\tau)\|_{L^\infty}\}\int_0^t e^{-\nu(p)(t-s)}\int_{\mathbb{R}^3}|K(p,q)e^{\f{|p-q|}{20}}|dq\nonumber\\
&~~~~~~~~~~~~~~~~~~~~~~~~~~~~~~~~~~~~~~~~~\times\int_{\mathbb{R}^3}|K(q,\eta)e^{-\f{|q-\eta|}{20}}|d\eta d\tau\int_0^se^{-\nu(q)(s-\tau)}(1+\tau)^{-k}d\tau  ds \nonumber\\
&\leq Ce^{-\f{N}{20}}(1+t)^{-k}\sup_{0\leq\tau\leq t}\{(1+\tau)^k\|\varphi(\tau)\|_{L^\infty}\}
.
\end{align}

\noindent{\it Case 3.} For $|p|\leq N, |q|\leq 2N, |\eta|\leq 3N$, we have
\begin{align}\label{8.15-4}
&\int_0^t e^{-\nu(p)(t-s)}\int_{|q|\leq 2N, |\eta|\leq3N}|K(p,q)K(q,\eta)|\int_0^se^{-\nu(q)(s-\tau)}\varphi(\tau,\hat{y},\eta)d\tau d\eta dq  ds \nonumber\\
&\leq \int_0^t e^{-\nu(p)(t-s)}\int_{|q|\leq 2N, |\eta|\leq3N}|K(p,q)K(q,\eta)|\int_{s-\k}^se^{-\nu(q)(s-\tau)}\varphi(\tau,\hat{y},\eta)d\tau d\eta dq  ds \nonumber\\
&~~~+\int_0^t e^{-c_N(t-s)}\int_{|q|\leq 2N, |\eta|\leq3N}|K(p,q)K(q,\eta)|\int_0^{s-\k}e^{-c_N(s-\tau)}\varphi(\tau,\hat{y},\eta)d\tau d\eta dq  ds,
\end{align}
where we have used  $\nu(p)\geq c_N\cong N^{-\f{b}2}~~\mbox{for}~|p|\leq N,~~\mbox{and}~~\nu(q)\geq c_N\cong N^{-\f{b}2}~~\mbox{for}~|q|\leq 2N$.
A directly calculation shows that the first term on the RHS of \eqref{8.15-4} can be estimated by
\begin{align}\label{8.15-5}
&C\k\sup_{0\leq\tau\leq t}\{(1+\tau)^k\|\varphi(\tau)\|_{L^\infty}\}\int_0^t e^{-\nu(p)(t-s)}(1+s)^{-k} p_0^{-\f{b}2-\xi_1}ds\nonumber\\
&\leq C\k(1+t)^{-k}\sup_{0\leq\tau\leq t}\{(1+\tau)^k\|\varphi(\tau)\|_{L^\infty}\}
\end{align}
Next we estimate the second term on the RHS of \eqref{8.15-4}.  Since $k(p,q)$ may have singularity of $|p-q|^{-1}$,  as previously, we  choose a smooth compact support function $\tilde{k}_N(p,q)$ such that
\begin{align}\label{8.15-6}
\sup_{|p|\leq 3N}\int_{|q|\leq 3N}\Big|k(p,q)-\tilde{k}_N(p,q)\Big|dq\leq C N^{-2b-2}.
\end{align}
Noting
\begin{align}\label{8.12-11}
k(p,q)k(q,\eta)&=\Big(k(p,q)-\tilde{k}_N(p,q)\Big)k(q,\eta)\nonumber\\
&\quad+\Big(k(q,\eta)-\tilde{k}_N(q,\eta)\Big)\tilde{k}_N(p,q)+\tilde{k}_N(p,q)\tilde{k}_N(q,\eta),
\end{align}
and  using \eqref{2.47}, \eqref{8.15-6}  and \eqref{8.12-11}, we can bound the second term on the RHS of \eqref{8.15-4} by
\begin{align}\label{8.12-12}
&N^{-2b-2}\sup_{0\leq \tau\leq t}\Big\{(1+\tau)^{-k}\|\varphi(\tau)\|_{L^\infty}\Big\}\int_0^te^{-c_N(t-s)}\int_0^{s-\k}e^{-c_N(s-\tau)} (1+\tau)^{-k}d\tau ds \nonumber\\
&~+\int_0^te^{-c_N(t-s)} \int_{|q|\leq 2N,|\eta|\leq3N}|\tilde{k}_N(p,q)\tilde{k}_N(q,\eta)|\int_0^{s-\k}e^{-c_N(s-\tau)} |\varphi(\tau,\hat{x}-\hat{q}(s-\tau),\eta)|d\tau d\eta dqds\nonumber\\
&\leq \f{C}{N}(1+t)^{-k}\sup_{0\leq s\leq t}\Big\{(1+\tau)^{-k}\|\varphi(\tau)\|_{L^\infty}\Big\}\nonumber\\
&\quad+C_{N}\int_0^t\int_0^{s-\k}e^{-c_N(t-s)} e^{-c_N(s-\tau)} \int_{|q|\leq 2N,|\eta|\leq3N} |\varphi(\tau,\hat{x}-\hat{q}(s-\tau),\eta)|d\eta dqd\tau ds,
\end{align}
where we have used the facts that   $\tilde{k}_N(p,q)\tilde{k}_N(q,\eta)$ is bounded.
As in Section 4, using the changing of variables, one can obtain that
\begin{align}\label{8.12-14}
& C_{N}\int_0^t\int_0^{s-\k}e^{-c_N(t-s)} e^{-c_N(s-\tau)} \int_{|q|\leq 2N,|\eta|\leq3N} |\varphi(\tau,\hat{x}-\hat{q}(s-\tau),\eta)|d\eta dqd\tau ds
\nonumber\\
&\leq C_{N}\int_0^t\int_0^{s-\k}e^{-c_N(t-s)} e^{-c_N(s-\tau)} \|\varphi(\tau)\|_{L^2}d\tau\nonumber\\
&\leq C_{N}(1+t)^{-k}\sup_{0\leq\tau\leq t}\Big\{(1+\tau)^k\|\varphi(\tau)\|_{L^2}\Big\}
\leq C_{N}(1+t)^{-k}\|w_{k\f{b}2}\varphi_0\|_{L^2}\nonumber\\
&\leq C_{N}(1+t)^{-k}\|w_{\vartheta+k\f{b}2}\varphi_0\|_{L^\infty},
\end{align}
where we have used $\vartheta>\f32$ and  \eqref{6.23} with $\a=0$. Thus combining \eqref{8.15-4}-\eqref{8.12-14}, one gets that
\begin{align}\label{8.12-15}
&\int_0^t e^{-\nu(p)(t-s)}\int_{|q|\leq 2N, |\eta|\leq3N}|K(p,q)K(q,\eta)|\int_0^se^{-\nu(q)(s-\tau)}\varphi(\tau,\hat{y},\eta)d\tau d\eta dq  ds \nonumber\\
&\leq  C\Big(\k +\f{1}{N}\Big)(1+t)^{-k}\cdot\sup_{0\leq \tau\leq t}\Big\{(1+\tau)^k\|\varphi(\tau)\|_{L^\infty}\Big\}+C_{N}(1+t)^{-k}\|w_{\vartheta+k\f{b}2}\varphi_0\|_{L^\infty}.
\end{align}
Therefore, it follows from \eqref{8.15-2},  \eqref{8.15-3} and \eqref{8.12-15} that
\begin{align}\label{8.16-1}
L_{3}&\leq  C\Big[\k +\f{1}{N^{\f{\xi_1}{b}}}+\f1N\Big](1+t)^{-k}\cdot\sup_{0\leq \tau\leq t}\Big\{(1+\tau)^k\|\varphi(\tau)\|_{L^\infty}\Big\}+C_{N}(1+t)^{-k}\|w_{\vartheta+k\f{b}2}\varphi_0\|_{L^\infty}.
\end{align}
which, together with \eqref{8.8}  and  \eqref{8.12}, yield that
\begin{align}
&\sup_{0\leq s\leq t}\Big\{(1+s)^k\|\varphi(s)\|_{L^\infty}\Big\}\nonumber\\
&\leq C\Big[\k +\f{1}{N^{\f{\xi_1}{b}}}+\f1N\Big]\sup_{0\leq \tau\leq t}\Big\{(1+\tau)^k\|\varphi(\tau)\|_{L^\infty}\Big\}+C_{N}(1+t)^{-k}\|w_{\vartheta+k\f{b}2}\varphi_0\|_{L^\infty}.\nonumber
\end{align}
By first choosing  $\k$ small, then letting $N$ sufficiently large so that $C\Big[\k +\f{1}{N^{\f{\xi_1}{b}}}+\f1N\Big]\leq \f12$, one obtains that
\begin{equation*}\label{8.15}
\|\varphi(t)\|_{L^\infty}\leq C(1+t)^{-k}\|w_{\vartheta+k\f{b}2}\varphi_0\|_{L^\infty},
\end{equation*}
for all $t\geq 0$. This yields immediately \eqref{8.5}. Therefore we complete the proof of this lemma.    $\hfill\Box$

\

Based on the above preparations, we begin to prove Theorem \ref{thm1.3}.

\noindent{\bf Proof of Theorem \ref{thm1.3}:}  Using the semigroup $S(t)$ for the linearized  relativistic Boltzmann equation \eqref{6.1}, by the Duhamel Principle,   we have the solution formula  for the nonlinear relativistic  Boltzmann equation \eqref{1.22} as
\begin{align}
f(t)=S(t)f_0+\int_0^tS(t-s)\Big\{\Gamma(f,f)(s)\Big\}ds.\notag
\end{align}
From now on, we take $k:=1+\f{\xi_1}b>1$. Then it follows from \eqref{8.5} that
\begin{align}\label{8.17}
\|f(t)\|_{L^\infty}
&\leq C(1+t)^{-k}\|w_{\vartheta+k\f{b}{2}}f_0\|_{L^\infty}+C\int_0^t(1+t-s)^{-k}\left\|w_{\vartheta+k\f{b}{2}}\Big\{\Gamma(f,f)(s)\Big\}\right\|_{L^\infty}ds\nonumber\\
&\leq C(1+t)^{-k}\|w_{\vartheta+k\f{b}{2}}f_0\|_{L^\infty}+C\int_0^t(1+t-s)^{-k}\left\|w_{\vartheta+k\f{b}{2}}\Big\{\Gamma(f,f)(s)\Big\}\right\|_{L^\infty}ds.
\end{align}
Noting $0<\xi_1\leq \f14$ and taking $\vartheta=\f74$, it follows from \eqref{4.6}  that
\begin{align}\label{8.20}
&\Big|(1+|p|)^{\vartheta+k\f{b}2}\Big\{\Gamma(f,f)(s,x-\hat{p}(t-s),p)\Big\}\Big|\nonumber\\
&\leq C\|w_{\vartheta+\f{\xi_1}2}f(s)\|_{L^\infty}\|w_{1}f(s)\|^{\f{4d+1}{5d}}_{L^\infty}\sup_{z\in\Omega}\Big(\int_{\mathbb{R}^3}|f(s,z,\eta)|d\eta\Big)^{\f{d-1}{5d}}\nonumber\\
&\leq C\|w_{\b}f(s)\|^{\f{\vartheta+\f{\xi_1}2}\b}_{L^\infty}\|f(s)\|^{1-\f{\vartheta+\f{\xi_1}2}\b}_{L^\infty}
\|w_{\b}f(s)\|^{\f1{\b}\f{4d+1}{5d}}_{L^\infty}\|f(s)\|^{\f{4d+1}{5d}(1-\f1{\b})}_{L^\infty}\sup_{z\in\Omega}\Big(\int_{\mathbb{R}^3}|f(s,z,\eta)|d\eta\Big)^{\f{d-1}{5d}}\nonumber\\
&\leq C\|f(s)\|^{1-\f{\vartheta+\f{\xi_1}2}\b+\f{4d+1}{5d}(1-\f1{\b})}_{L^\infty}\|w_{\b}f(s)\|^{\f{\vartheta+\f{\xi_1}2}\b+\f{4d+1}{5d}\f1{\b}}_{L^\infty}\sup_{z\in\Omega}\Big(\int_{\mathbb{R}^3}|f(s,z,\eta)|d\eta\Big)^{\f{d-1}{5d}}\nonumber\\
&\leq C\|f(s)\|_{L^\infty}\|w_{\b}f(s)\|^{\f{4d+1}{5d}}_{L^\infty}\sup_{z\in\Omega}\Big(\int_{\mathbb{R}^3}|f(s,z,\eta)|d\eta\Big)^{\f{d-1}{5d}},
\end{align}
where we have used $\f{4d+1}{5d}(1-\f1{\b})-\f{\vartheta+\f{\xi_1}2}\b\geq 0$ due to  $\b>14$. Then it follows from \eqref{8.20} that
\begin{align}
&C\int_0^t(1+t-s)^{-k}\left\|(1+|p|)^{\vartheta+k\f{b}2}\Big\{\Gamma(f,f)(s)\Big\}\right\|_{L^\infty}ds\nonumber\\
&\leq C\int_0^t(1+t-s)^{-k}(1+s)^{-k}\nonumber\\
&\qquad\times \sup_{0\leq s\leq t,z\in\Omega}\left\{\Big[(1+s)^k\|f(s)\|_{L^\infty}\Big]\|w_{\b}f(s)\|^{\f{4d+1}{5d}}_{L^\infty}\Big(\int_{\mathbb{R}^3}|f(s,z,\eta)|d\eta\Big)^{\f{d-1}{5d}}\right\}ds\nonumber\\
&\leq C(1+t)^{-k}\sup_{0\leq s\leq t,z\in\Omega}\left\{\Big[(1+s)^k\|f(s)\|_{L^\infty}\Big]\|w_{\b}f(s)\|^{\f{4d+1}{5d}}_{L^\infty}\Big(\int_{\mathbb{R}^3}|f(s,z,\eta)|d\eta\Big)^{\f{d-1}{5d}}\right\},\nonumber
\end{align}
which together with \eqref{5.27} and \eqref{8.17}, yield that for $\b>14$,
\begin{align}\label{8.25}
&\sup_{0\leq s\leq t}\Big\{(1+s)^k\|f(s)\|_{L^\infty}\Big\}\nonumber\\
&\leq C\|w_\b f_0\|_{L^\infty}+\sup_{0\leq s\leq t,z\in\Omega}\left\{\Big[(1+s)^k\|f(s)\|_{L^\infty}\Big]\|w_{\b}f(s)\|^{\f{4d+1}{5d}}_{L^\infty}\Big(\int_{\mathbb{R}^3}|f(s,z,\eta)|d\eta\Big)^{\f{d-1}{5d}}\right\}\nonumber\\
&\leq C_4\bar{M}^4+C_4\sup_{1\leq s\leq t}\Big[(1+s)^r\|f(s)\|_{L^\infty}\Big]\sup_{1\leq s\leq t,z\in\Omega}\left\{\|w_{\b}f(s)\|^{\f{4d+1}{5d}}_{L^\infty}\Big(\int_{\mathbb{R}^3}|f(s,z,\eta)|d\eta\Big)^{\f{d-1}{5d}}\right\}.
\end{align}
Then, using \eqref{5.24} and  similar arguments as in \eqref{5.26},  if $\v_0$ is small enough, one can obtain that
\begin{align}\label{8.29}
&C_4\sup_{1\leq s\leq t,y\in\Omega}\left\{\|f(s)\|^{\f{4d+1}{5d}}_{L^\infty}\Big(\int_{\mathbb{R}^3}|f(s,y,\eta)|d\eta\Big)^{\f{d-1}{5d}}\right\}\leq \f12.
\end{align}
Substituting \eqref{8.29} into \eqref{8.25}, one proves that for $\b>14$,
\begin{align}\label{8.30}
\|f(t)\|_{L^\infty}\leq 2C_4\bar{M}^4(1+t)^{-k},\quad \forall\,t\geq0.
\end{align}
Taking
$$\tilde{C}_3=2C_4\bar{M}^4,$$
then we obtain \eqref{1.30} from \eqref{8.30}.  Therefore the proof of Theorem \ref{thm1.3} is completed. $\hfill\Box$


\section{Appendix}

\begin{lemma}[Glassey\& Strauss \cite{Glassey1}]\label{lemA.1}
Let $R>r\geq0$, and define
\begin{equation*}\label{A.2}
J_1(R,r):=\int_0^\infty z e^{-R\sqrt{1+z^2}} I_0(rz)dz,~~J_2(R,r):=\int_0^\infty \f{z}{\sqrt{1+z^2}} e^{-R\sqrt{1+z^2}} I_0(rz)dz.
\end{equation*}
Then it holds that
\begin{align}\label{A.3}
J_1(R,r)=\f{R}{R^2-r^2}\Big[1+\f1{\sqrt{R^2-r^2}}\Big]e^{-\sqrt{R^2-r^2}},~~~J_2(R,r)=\f1{\sqrt{R^2-r^2}}e^{-\sqrt{R^2-r^2}}.
\end{align}
\end{lemma}

\

The following Lemmas \ref{lemA.2} and \ref{lemA.3} can be regarded as a refine version of \cite{Glassey1}. These refined Lemmas play important roles in the proof of Theorem \ref{thm7.1} and Lemma \ref{lem4.1}.
\begin{lemma}\label{lemA.2}
For $l,j$ defined in \eqref{A.4}, we have, for $\a\in[-2,2]$, that
\begin{align}\label{A.5}
K_\a(l,j):=\int_0^\infty z(1+z^2)^{\f\a4} e^{-l\sqrt{1+z^2}} I_0(jz)dz\leq C\f{l^{1+\f\a2}}{(l^2-j^2)^{1+\f{\a}4}}e^{-\sqrt{l^2-j^2}}.
\end{align}
\end{lemma}
\noindent{\bf Proof.} From \eqref{A.3} and \eqref{2.9}, we know that
\begin{align}\label{A.5-1}
K_{-2}(l,j)=\f{1}{\sqrt{l^2-j^2}}e^{-\sqrt{l^2-j^2}},~~~K_0(l,j)\leq \f{2l}{l^2-j^2}e^{-\sqrt{l^2-j^2}},
\end{align}
It follows from \cite{Glassey1}(Page 323) that
\begin{align}\label{A.5-2}
K_2(l,j)&=(l^2-j^2)^{-\f52}\Big\{[(l^2-j^2)+3(l^2-j^2)^{\f12}+3]l^2-(l^2-j^2)-(l^2-j^2)^{\f32}\Big\}e^{-\sqrt{l^2-j^2}}\nonumber\\
&\leq \f{Cl^2}{(l^2-j^2)^{\f32}}e^{-\sqrt{l^2-j^2}}.
\end{align}
Combining \eqref{A.5-1}, \eqref{A.5-2} and Holder inequality, one can prove \eqref{A.5} directly. The details are omitted here. $\hfill\Box$

\

\begin{lemma}\label{lemA.3}
Let $0<\g<2$ and $l,j$ be as in \eqref{A.4}. For any given small $\v>0$, it holds that
\begin{align}\label{A.6}
I_\g(l,j):= \int_0^1z^{1-\g} e^{-l\sqrt{1+z^2}}I_0(jz)dz\leq \f{C}{(l^2-j^2)^{\f1{2r}}} e^{-\sqrt{l^2-j^2}},
\end{align}	
where we choose $r:=1+\v+\f{\g}{2-\g}$ and  $\f1{2r}=\f12-\f{\g}4-\f{(2-\g)^2\v}{4(2+2\v-\g\v)}$.
\end{lemma}
\noindent{\bf Proof.}  It follows from the H\"{o}lder inequality that
\begin{align}\label{A.7}
I_\g(l,j)&\leq \Big(\int_0^1 z^{1-\f{r}{r-1}\gamma}dz\Big)^{1-\f1r}\Big(\int_0^1z e^{-rl\sqrt{1+z^2}}I_0(jz)^rdz\Big)^{\f1r}\nonumber\\
&\leq C_\g\Big(\int_0^1z e^{-rl\sqrt{1+z^2}}I_0(jz)^rdz\Big)^{\f1r},
\end{align}
where we have used the fact $1-\f{r}{r-1}\g>-1$ in the last inequality. Noting
\begin{equation*}\label{A.8}
I_0(jz)^r\leq \Big(\f{1}{2\pi}\int_{0}^{2\pi} e^{jz\cos\varphi}d\varphi\Big)^r\leq C_r\f{1}{2\pi}\int_{0}^{2\pi} e^{rjz\cos\varphi}d\varphi\leq C_r I_0(rjz),
\end{equation*}
which, together with \eqref{A.7}, yields that
\begin{align}
I_\g(l,j)&\leq C_\g\Big(\int_0^1z e^{-rl\sqrt{1+z^2}}I_0(rjz)dz\Big)^{\f1r}
\leq C_\g\Big(\int_0^{\infty}z (1+z^2)^{-\f12}e^{-rl\sqrt{1+z^2}}I_0(rjz)dz\Big)^{\f1r}\nonumber\\
&\leq C_\g \Big(\f{1}{r\sqrt{l^2-j^2}}e^{-r\sqrt{l^2-j^2}}\Big)^{\f1r}\leq C_\g\f{1}{(l^2-j^2)^{\f1{2r}}} e^{-\sqrt{l^2-j^2}},\nonumber
\end{align}
where we have used \eqref{A.5} with $\a=-2$. Therefore, the proof of Lemma \ref{lemA.3} is completed. $\hfill\Box$

\

\begin{lemma}\label{lemA.4}
For any fixed $\vartheta\in\mathbb{R}$, it holds that
\begin{align}\label{A.9}
(p_0+q_0)^{\vartheta}e^{-c|p-q|}\leq C (p_0q_0)^{\f{\vartheta}2}e^{-\f{c}2|p-q|},
\end{align}
and
\begin{align}\label{A.10}
 (p_0q_0)^\vartheta e^{-c|p-q|}\leq C p_0^{2\vartheta}e^{-\f{c}2|p-q|}.
\end{align}
\end{lemma}
\noindent{\bf Proof.} We divide it into three cases.\\
\noindent{\it Case 1:} $\f12|p|\leq |q|\leq 2|p|$.  For this case, \eqref{A.9} and \eqref{A.10} are easy to obtain.\\[1mm]
\noindent{\it Case 2:} $\f12|p|\geq |q|$. A directly calculation shows that
\begin{equation*}\label{A.11}
|p-q|\geq \f12|p-q|+\f12(|p|-|q|)\geq \f12|p-q|+\f14|p|,
\end{equation*}
which yields immediately that
\begin{align}
(p_0+q_0)^\vartheta e^{-c|p-q|}
&\lesssim
\begin{cases}
p_0^\vartheta e^{-\f{c}4|p|} e^{-\f{c}2|p-q|},~~\mbox{for}~\vartheta\geq0\\
q_0^{\vartheta}e^{-\f{c}4|p|} e^{-\f{c}2|p-q|},~~\mbox{for}~\vartheta\leq0
\end{cases}
\lesssim (p_0q_0)^{\f{\vartheta}2}e^{-\f{c}2|p-q|},\nonumber
\end{align}
and
\begin{align}
 (p_0q_0)^\vartheta e^{-c|p-q|}
 &\lesssim
 \begin{cases}
 p_0^{2\vartheta} e^{-\f{c}4|p|} e^{-\f{c}2|p-q|},~~\mbox{for}~\vartheta\geq0\\
e^{-\f{c}4|p|} e^{-\f{c}2|p-q|},~~\mbox{for}~\vartheta\leq0
 \end{cases}
 \lesssim p_0^{2\vartheta}e^{-\f{c}2|p-q|}.\nonumber
\end{align}
\noindent{\it Case 3:} $\f12|q|\geq |p|$. A directly calculation shows that
\begin{equation*}\label{A.14}
|p-q|\geq \f12|p-q|+\f12(|q|-|p|)\geq \f12|p-q|+\f14|q|,
\end{equation*}
which yields immediately that
\begin{align}
(p_0+q_0)^\vartheta e^{-c|p-q|}
&\lesssim
\begin{cases}
q_0^\vartheta e^{-\f{c}4|q|} e^{-\f{c}2|p-q|},~~\mbox{for}~\vartheta\geq0\\
p_0^{\vartheta}e^{-\f{c}4|q|} e^{-\f{c}2|p-q|},~~\mbox{for}~\vartheta\leq0
\end{cases}
\lesssim (p_0q_0)^{\f{\vartheta}2}e^{-\f{c}2|p-q|},\nonumber
\end{align}
and
\begin{align}
(p_0q_0)^\vartheta e^{-c|p-q|}
&\lesssim
\begin{cases}
q_0^{2\vartheta} e^{-\f{c}4|q|} e^{-\f{c}2|p-q|},~~\mbox{for}~\vartheta\geq0\\
p_0^{2\vartheta} e^{-\f{c}4|q|} e^{-\f{c}2|p-q|},~~\mbox{for}~\vartheta\leq0
\end{cases}
\lesssim p_0^{2\vartheta}e^{-\f{c}2|p-q|}.\nonumber
\end{align}
Combining all the above estimates, we complete the proof of Lemma \ref{lemA.4}. $\hfill\Box$

\begin{lemma}
Let $0\leq \a<3$, it holds that
\begin{align}\label{A.17}
\int_{\mathbb{R}^3}\f{e^{-|p-q|}}{[|p\times q|+|p-q|]^\a} dq\leq
\begin{cases}
C_\a (1+|p|)^{-\a},~~\mbox{for}~0\leq \a<2,\\
C(1+|p|)^{-2}\ln(1+|p|),~~\mbox{for}~\a=2,\\
C_\a (1+|p|)^{-2},~~\mbox{for}~\a>2.\\
\end{cases}
\end{align}
\end{lemma}
\noindent{\bf Proof.} We need only to calculate the case for $|p|\geq2$. Using the Polar coordinates, we have
\begin{align}
\int_{\mathbb{R}^3}\f{e^{-|p-q|}}{[|p\times q|+|p-q|]^\a} dq
&\lesssim \int_0^\infty e^{-r}r^{2-\a} dr\int_0^{\pi}\f{\sin\t}{[1+|p|\sin\t]^{\a}}d\t\nonumber\\
&\lesssim \int_0^{\pi}\f{\sin\t}{[1+|p|\sin\t]^{\a}}d\t
\lesssim \int_0^{\f\pi2}\f{\sin\t}{[1+|p|\sin\t]^{\a}}d\t\nonumber\\
&\lesssim -\int_0^{\f\pi2}\f{d\cos\t}{[1+|p|^2-|p|^2\cos^2\t]^{\f\a2}}
\lesssim \int_0^1\f{dy}{[1+|p|^2-|p|^2y^2]^{\f\a2}}\nonumber\\
&\cong \f{1}{|p|^\a} \int_0^1 \f{dy}{[1+\f1{|p|^2}-y^2]^{\f\a2}}
\cong \f{1}{|p|^\a} \int_0^1 \f{dy}{[\sqrt{1+\f1{|p|^2}}-y]^{\f\a2}}.\nonumber
\end{align}
For $0\leq \a<2$, a direct calculation shows that
\begin{align}
\int_0^1 \f{dy}{[\sqrt{1+\f1{|p|^2}}-y]^{\f\a2}}
=-\f{1}{1-\f\a2}\Big[\sqrt{1+\f1{|p|^2}}-y\Big]^{1-\f\a2}\Big|_0^1 \lesssim 1.\nonumber
\end{align}
For $\a=2$, one can obtains
\begin{align}
\int_0^1 \f{dy}{\sqrt{1+\f1{|p|^2}}-y}&=- \ln(\sqrt{1+\f1{|p|^2}}-y)\Big|_0^1=\ln\sqrt{1+\f1{|p|^2}}-\ln(\sqrt{1+\f1{|p|^2}}-1)
\nonumber\\
&\lesssim  \ln|p|.\nonumber
\end{align}
For $2< \a<3$, a direct calculation shows that
\begin{align}
\int_0^1 \f{dy}{[\sqrt{1+\f1{|p|^2}}-y]^{\f\a2}}
&=-\f{1}{1-\f\a2}\Big[\sqrt{1+\f1{|p|^2}}-y\Big]^{1-\f\a2}\Big|_0^1 \nonumber\\
&\lesssim \Big[\sqrt{1+\f1{|p|^2}}-1\Big]^{1-\f\a2}\lesssim |p|^{-2+\a}.\nonumber
\end{align}
Combining the above estimates, we complete the proof of \eqref{A.17}.   $\hfill\Box$

\begin{lemma}\label{lemA.6}
For $0\leq \a_1<3, ~\a_2\geq0$ and $\b>3+\a_1+\a_2$,  it holds that
\begin{align}\label{A.18}
\int_{\mathbb{R}^3} \f{(1+|q|)^{-\b}}{|p-q|^{\a_1}(1+|p-q|)^{\a_2}}dq\leq C(1+|p|)^{-\a_1-\a_2}.
\end{align}
\end{lemma}
Here we omit the proof of Lemma \ref{lemA.6} since the proof  is standard.

\

\begin{lemma}\label{lemA.7}
Let \eqref{1.30}, \eqref{1.31} hold, and $\g>-2, a\in[0, 2+\g]$, $b\in[0,\min\{4,4+\g\})$.   For any fixed $c>0$ and $1\leq d<\min\left\{ \f2{\max\{-\g,1\}},\f3{\max\{b-1,1\}}\right\}$, it holds that
\begin{equation}\label{A.19}
\int_{\mathbb{R}^3\times \mathbb{S}^2} |v_{\phi}\s(g,\t)|^{d} e^{-cq_0}d\omega dq
\cong
\begin{cases}
(p_0^{\f{a}2})^{d},~\mbox{for hard potentials},\\
(p_0^{-\f{b}2})^d,~\mbox{for soft potentials},
\end{cases}
\end{equation}
which immediately yields 
\begin{align}\label{n1.20-1}
\nu(p)\cong
\begin{cases}
p_0^{\f{a}2},~\mbox{for hard potentials},\\
p_0^{-\f{b}2},~\mbox{for soft potentials}.
\end{cases}
\end{align}
\end{lemma}
{\bf Proof.} Firstly, we calculate the the upper bound.  It is noted that
\begin{align}\label{A.20}
\s(g,\t)\lesssim
\begin{cases}
\s_a(g,\t)+\s_{b}(g,\t),~\mbox{for hard potentials},\\
\s_{b}(g,\t),~\mbox{for hard potentials},
\end{cases}
\end{align}
where $\s_a(\cdot,\cdot)$ and $\s_b(\cdot,\cdot)$ are defined in \eqref{2.15-1} above.
Thus we need only to calculate
\begin{align}\label{A.21}
\int_{\mathbb{R}^3\times \mathbb{S}^2} |v_{\phi}\s_a(g,\t)|^{d} e^{-cq_0}d\omega dq~~\mbox{and}~~\int_{\mathbb{R}^3\times \mathbb{S}^2} |v_{\phi}\s_b(g,\t)|^{d} e^{-cq_0}d\omega dq.
\end{align}
Since $g\lesssim \sqrt{p_0q_0},~s\lesssim p_0q_0$, we have
\begin{equation*}
|v_{\phi}\s_a(g,\t)|^{d}=\Big|\f{\sqrt{s}}{p_0q_0}g^{1+a} \sin^{\g}\t\Big|^{d}
\lesssim (p_0q_0)^{\f{da}2}\sin^{d\g}\t,
\end{equation*}
which implies that
\begin{align}\label{A.22}
\int_{\mathbb{R}^3\times \mathbb{S}^2} |v_{\phi}\s_a(g,\t)|^{d} e^{-cq_0}d\omega dq&\lesssim \int_{\mathbb{R}^3} (p_0q_0)^{\f{da}2}e^{-cq_0}dq\int_0^{\pi}\sin^{1+d\g}\t d\t \lesssim p_0^{\f{da}2},
\end{align}
where we have used the fact $1+\g d>-1$.\\

To estimate the second part of \eqref{A.21},
we divide the proof into two cases.\\[1mm]
\noindent{\it Case 1.} For $0\leq b<1$, it follows from $g\lesssim \sqrt{p_0q_0},~s\lesssim p_0q_0$ that
\begin{equation*}
|v_{\phi}\s_b(g,\t)|^{d}=\Big|\f{\sqrt{s}}{p_0q_0}g^{1-b} \sin^{\g}\t\Big|^{d}
\lesssim (p_0q_0)^{\f{-b}{2}}\sin^{d\g}\t,
\end{equation*}
which yields that
\begin{align}\label{A.23}
\int_{\mathbb{R}^3\times \mathbb{S}^2} |v_{\phi}\s_b(g,\t)|^{d} e^{-cq_0}d\omega dq\leq \int_{\mathbb{R}^3} (p_0q_0)^{\f{-db}{2}} e^{-cq_0}dq \int_{0}^{\pi}\sin^{1+d\g}\t d\t  \lesssim p_0^{-\f{db}2}.
\end{align}

\noindent{\it Case 2.} For $1\leq b<\min\{4,4+\g\}$,  using $g^{-1}\lesssim \f{\sqrt{p_0q_0}}{|p-q|},~s\lesssim p_0q_0$, we have
\begin{equation*}
|v_{\phi}\s_b(g,\t)|^{d}=\Big|\f{\sqrt{s}}{p_0q_0}g^{1-b} \sin^{\g}\t\Big|^{d}
\lesssim \f{(p_0q_0)^{\f{d(b-2)}{2}}}{|p-q|^{d(b-1)}}\sin^{d\g}\t,
\end{equation*}
which yields that
\begin{align}\label{A.24}
\int_{\mathbb{R}^3\times \mathbb{S}^2} |v_{\phi}\s_b(g,\t)|^{d} e^{-cq_0}d\omega dq\leq \int_{\mathbb{R}^3} \f{(p_0q_0)^{\f{d(b-2)}{2}}}{|p-q|^{d(b-1)}}e^{-cq_0}dq \int_{0}^{\pi}\sin^{1+d\g}\t d\t  \lesssim p_0^{-\f{db}2},
\end{align}
where we have used that fact $0\leq d(b-1)<3$. Thus it follows from \eqref{A.23} and \eqref{A.24} that
\begin{align}\label{A.25}
\int_{\mathbb{R}^3\times \mathbb{S}^2} |v_{\phi}\s_b(g,\t)|^{d} e^{-cq_0}d\omega dq \leq p_0^{-\f{db}2}.
\end{align}
Combining \eqref{A.20}, \eqref{A.22} and \eqref{A.25}, we obtain
\begin{equation*}\label{A.19-1}
\int_{\mathbb{R}^3\times \mathbb{S}^2} |v_{\phi}\s(g,\t)|^{d} e^{-cq_0}d\omega dq
\lesssim
\begin{cases}
(p_0^{\f{a}2})^{d},~\mbox{for hard potentials},\\
(p_0^{-\f{b}2})^d,~\mbox{for soft potentials},
\end{cases}
\end{equation*}

Next we calculate the lower bound. It follows from \eqref{2.7} and $s\lesssim p_0q_0$ that 
\begin{align}
v_{\phi}\f{g}{\sqrt{s}}g^a \gtrsim \f{|p-q|^{2+a}}{(p_0q_0)^{2+\f{a}2}}~~~\mbox{and}~~~v_{\phi}\f{g}{\sqrt{s}}g^{-b}\gtrsim
\begin{cases}
\f{|p-q|^{2-b}}{(p_0q_0)^{2-\f{b}2}},~&\mbox{for}~0\leq b\leq 2,\\
(p_0q_0)^{-\f{b}2},~&\mbox{for}~b>2,\nonumber
\end{cases}
\end{align}
which yield immediately that 
 \begin{equation*}\label{A.19-2}
 \int_{\mathbb{R}^3\times \mathbb{S}^2} |v_{\phi}\s(g,\t)|^{d} e^{-cq_0}d\omega dq
\gtrsim
 \begin{cases}
 (p_0^{\f{a}2})^{d},~\mbox{for hard potentials},\\
 (p_0^{-\f{b}2})^d,~\mbox{for soft potentials}.
 \end{cases}
 \end{equation*}
Therefore we complete the proof of Lemma \ref{lemA.7}.  $\hfill\Box$

\

\noindent{\bf Proof of Lemma \ref{lem8.1}:} It follows from \eqref{1.23} that
\begin{align}\label{A.26}
\langle w_{\vartheta}^2Lf,f\rangle=|w_{\vartheta}f|_{\nu}^2-\langle w_{\vartheta}Kf,w_{\vartheta}f\rangle.
\end{align}
Using \eqref{2.47} and H\"{o}lder inequality, one has
\begin{align}\label{A.27}
|\langle w_{\vartheta}Kf,w_{\vartheta}f\rangle|
&=\Big|\int_{\mathbb{R}^3}w_{\vartheta}(p)f(p)\int_{\mathbb{R}^3}k_{w_{\vartheta}}(p,q)w_{\vartheta}(q)f(q)dqdp\Big|\nonumber\\
&\leq \Big(\int_{\mathbb{R}^3}\nu(p)|w_{\vartheta}(p)f(p)|^2dp\Big)^{\f12}
\Big(\int_{\mathbb{R}^3}\nu(p)^{-1}\Big|\int_{\mathbb{R}^3}k_{w_{\vartheta}}(p,q)w_{\vartheta}(q)f(q)dq\Big|^2dp\Big)^{\f12}\nonumber\\
&\leq C|w_{\vartheta}f|_{\nu}
\Big(\int_{\mathbb{R}^3}\int_{\mathbb{R}^3}p_0^{-\xi_1}|k_{w_{\vartheta}}(p,q)|\cdot|w_{\vartheta}(q)f(q)|^2dqdp\Big)^{\f12}\nonumber\\
&\leq \f14|w_{\vartheta}f|_{\nu}^2+\int_{\mathbb{R}^3}\nu(q)q_0^{-2\xi_1}|w_{\vartheta}(q)f(q)|^2dq\nonumber\\
&\leq \f14|w_{\vartheta}f|_{\nu}^2+C_R|I_{\leq R}f|_{L^2}^2+CR^{-2\xi_1}\int_{|q|\geq R}\nu(q)|w_{\vartheta}(q)f(q)|^2dq\nonumber\\
&\leq \f12|w_{\vartheta}f|_{\nu}^2+C_R|I_{\leq R}f|_{L^2}^2,
\end{align}
where we have chosen  $R\gg1$ so that $CR^{-2\xi_1}\leq \f14$. Substituting \eqref{A.27} into \eqref{A.26}, we get \eqref{6.22}. Therefore, the proof of Lemma \ref{lem8.1} is completed. $\hfill\Box$

\

\noindent{\bf Acknowledgments.} 
Yong Wang is partly supported by NSFC Grant No.  11771429, 11671237 and  11688101. The author greatly appreciates the referees for their invaluable comments, which helped to improve the presentation of the paper. 

\end{document}